\documentclass[a4paper,10pt,oneside]{article} 
\usepackage{geometry}
\usepackage{standalone}
\usepackage[T1] {fontenc}
\usepackage [utf8] {inputenc}
\usepackage[english] {babel}
\usepackage{xcolor}
\usepackage{graphicx}
\usepackage{authblk}
\usepackage{bbm}
\usepackage{amsmath,amssymb,amsthm, xcolor, listings}
\usepackage{mathrsfs}
\usepackage{tikz}
\usetikzlibrary{patterns}
\usetikzlibrary{decorations.pathreplacing}
\usepackage{paralist}
\usepackage{fancyhdr}
\usepackage{mathtools}
\usepackage{titlesec}
\usepackage{enumitem}
\usepackage{float,placeins}
\usepackage{url}
\usepackage{hyperref}
\usepackage{bm}
\usepackage{mathdots}
\usepackage{microtype} 
\usepackage{eurosym}

\newtheorem{theorem}{Theorem}[section]
\newtheorem{proposition}{Proposition}[section]
\newtheorem{lemma}{Lemma}[section]
\newtheorem{cor}{Corollary}[section]

\theoremstyle{definition}
\newtheorem{definition}{Definition}[section]
\newtheorem{remark}{Remark}[section]

\numberwithin{equation}{section}
\newtheorem{assumption}{Assumption}[section]

\newcommand{\mycircleblack}{\tikz\draw[black, fill=black!88!white] (0,0) circle (3pt);}
\newcommand{\mycircleintgray}{\tikz\draw[black, fill=black!50!white] (0,0) circle (3pt);}
\newcommand{\mycirclegray}{\tikz\draw[black, fill=black!28!white] (0,0) circle (3pt);}

\newcommand{\mycirclelightgray}{\tikz\draw[black, fill=black!12!white] (0,0) circle (3pt);}
\newcommand{\mycircleerreprimo}{\tikz\draw[black, fill=black!70!white] (0,0) circle (3pt);}

\newcommand{\mycirclewhite}{\tikz\draw[black, fill=white] (0,0) circle (3pt);} 

\begin{document}
\title{Metastability for the degenerate Potts Model with\\ negative external magnetic field under Glauber dynamics.}

\author[
         {}\hspace{0.5pt}\protect\hyperlink{hyp:email1}{1},\protect\hyperlink{hyp:affil1}{a},
         \protect\hyperlink{hyp:corresponding}{$\dagger$}
        ]
        {\protect\hypertarget{hyp:author1}{Gianmarco Bet}}

\author[
         {}\hspace{0.5pt}\protect\hyperlink{hyp:email2}{2},\protect\hyperlink{hyp:affil3}{c}
        ]
        {\protect\hypertarget{hyp:author2}{Anna Gallo}}

\author[
         {}\hspace{0.5pt}\protect\hyperlink{hyp:email3}{3},\protect\hyperlink{hyp:affil1}{a},\protect\hyperlink{hyp:affil2}{b}
        ]
        {\protect\hypertarget{hyp:author3}{Francesca R.~Nardi}}

\affil[ ]{
          \small\parbox{365pt}{
             \parbox{5pt}{\textsuperscript{\protect\hypertarget{hyp:affil2}{a}}}Università degli Studi di Firenze (Italy),
            \enspace
             \parbox{5pt}{\textsuperscript{\protect\hypertarget{hyp:affil1}{b}}}Eindhoven University of Technology (The Netherlands),
             \enspace
             \parbox{5pt}{\textsuperscript{\protect\hypertarget{hyp:affil3}{c}}}IMT School for Advanced Studies Lucca (Italy)
            }
          }

\affil[ ]{
          \small\parbox{365pt}{
             \parbox{5pt}{\textsuperscript{\protect\hypertarget{hyp:email1}{1}}}\texttt{\footnotesize\href{mailto:gianmarco.bet@unifi.it}{gianmarco.bet@unifi.it}},
             \parbox{5pt}{\textsuperscript{\protect\hypertarget{hyp:email2}{2}}}\texttt{\footnotesize\href{mailto:anna.gallo@imtlucca.it}{anna.gallo@imtlucca.it}},
             \parbox{5pt}{\textsuperscript{\protect\hypertarget{hyp:email3}{3}}}\texttt{\footnotesize\href{mailto:francescaromana.nardi@unifi.it}{francescaromana.nardi@unifi.it}}
            }
          }

\affil[ ]{
          \small\parbox{365pt}{
             \parbox{5pt}{\textsuperscript{\protect\hypertarget{hyp:corresponding}{$\dagger$}}}Corresponding author
            }
          }

\date{\today}

\maketitle
\begin{abstract}
We consider the ferromagnetic $q$-state Potts model on a finite grid with non-zero external field and periodic boundary conditions. The system evolves according to Glauber-type dynamics described by the Metropolis algorithm, and we focus on the low temperature asymptotic regime. We analyze the case of negative external magnetic field. In this scenario there are $q-1$ stable configurations and a unique metastable state. We describe the asymptotic behavior of the first hitting time from the metastable state to the set of the stable states as $\beta\to\infty$ in probability, in expectation, and in distribution. We also identify the exponent of the mixing time and find an upper and a lower bound for the spectral gap. We identify the minimal gates for the transition from the metastable state to the set of the stable states and for the transition from the metastable state to a fixed stable state. Furthermore, we identify the tube of typical trajectories for these two transitions. The detailed description of the energy landscape that we develop allows us to give precise asymptotics for the expected transition time from the unique metastable state to the set of the stable configurations.  

\\
 \emph{Acknowledgment:} The research of Francesca R.~Nardi was partially supported by the NWO Gravitation Grant 024.002.003--NETWORKS and by the PRIN Grant 20155PAWZB ``Large Scale Random Structures''. The authors are grateful to Simone Baldassarri, Vanessa Jacquier and Cristian Spitoni for the detailed and fruitful discussions.  
\end{abstract}

\section{Introduction}\label{intro}
Metastability is a phenomenon that is observed when a physical system is close to a first--order phase transition. When a physical system lies close to its phase coexistence line, it may remain stuck for a long time in a state which is different from the equilibrium state. The former is known as the \textit{metastable state}. After a long (random) time, the system may perform a sudden transition from the metastable state to the stable state. When the system lies exactly on the phase coexistence line, it is of interest to understand precisely the \textit{tunneling transition} between two or more stable states. Many models for metastable behavior have been developed throughout the years. In these models 
a suitable stochastic dynamics is chosen and typically three main issues are investigated. The first is the study of the \textit{first hitting time} of the stable state(s) for the process started in the metastable state. The second issue is the study of the \textit{critical configurations} visited by the process  with probability close to one during the transition from the metastable state to the stable state(s). The final issue is the study of the \textit{tube of typical paths} of the process during the transition from the metastable state to the stable state(s). When a system lies on the phase coexistence line the same three issues above are investigated for the transition between any two stable states.
 
In this paper we study the metastable behavior of the $q$-state Potts model with non-zero external magnetic field on a finite two-dimensional discrete torus $\Lambda$. 
Each site $i$ of $\Lambda$ lies a spin with value $\sigma(i)\in\{1,\dots,q\}$, hence the $q$-state Potts model is an extension of the classical Ising model from $q=2$ to an arbitrary number $q$ of spins with $q>2$. To each configuration $\sigma$ is associated an energy $H(\sigma)$ that depends on the ferromagnetic interaction between nearest-neighbor spins, and on an external magnetic field $h$ which favors to a specific spin value. We focus on the regime of large inverse temperature $\beta\to\infty$. The stochastic evolution is described by a \textit{Glauber-type dynamics}, which is a Markov chain, given by the Metropolis algorithm, that only allows single spin flip updates.  This dynamics is reversible with respect to the \textit{Gibbs measure} $\mu_\beta$, see \eqref{gibbs}.

Our analysis focuses on the case of negative external magnetic field. In this scenario there are one metastable state and $q-1$ stable states. Without loss of generality, in the metastable configuration all spins are equal to $1$. The remaining constant configurations are stable states. We focus our attention on the transition from the metastable state to the set of stable configurations and from the metastable state to some fixed stable state.  When there is more than one stable state, these transitions are quite different because there may be intermediate transitions between different stable states. 

The goal of this paper is to investigate all the three issues of metastability introduced above for the $q$-state Potts model with negative external magnetic field. We focus on two classes of transitions: from the metastable state to the set of stable states (briefly denoted $\mathbf 1\to\mathcal X^s_\text{neg}$) and from the metastable state to any \textit{fixed} stable state (briefly denoted $\mathbf 1\to\mathbf s$). For both transitions, we investigate transition time, the minimal gates and the tube of typical trajectories. Finally, we identify the prefactor of the expected transition time.

Let us now briefly describe our approach. First we prove that the only metastable configuration is the configuration with all spins equal to $1$. For the transition $\mathbf 1\to\mathcal X^s_\text{neg}$, we are able to obtain the expected value and distribution of the transition time. This is more complicated for the transition $\mathbf 1\to\mathbf s$. Indeed, in this case with probability strictly positive the optimal path visits a stable state different from $\bold s$ before hitting $\bold s$. We prove that the energy barrier between two stable states is strictly larger than the energy barrier between a stable state and any other (non--stable) state. In view of this, we prove that the lower and the upper asymptotic bounds for the transition time have different exponents, see Remark \ref{labelintrorecall}. Moreover, we characterize the behavior of the \textit{mixing time} in the low-temperature regime and give an estimate of the \textit{spectral gap}, see \eqref{mixingtimedef} and \eqref{spectralgapdef} for the formal definitions. Next, we identify the set of all minimal gates. In particular, we prove that this set is given by those configurations in which all spins are $1$ except those, which are $s\in\{2,\dots,q\}$, in a quasi-square with a unit protuberance on one of the longest sides. The process hits the set of the stable configurations in any stable state with the same probability, thus it follows a uniform distribution over $\{2,\dots,q\}$. Using the so-called \textit{potential theoretic approach}, we give sharp estimates on the expected transition time by computing the so-called \textit{prefactor} explicitly. This requires a detailed knowledge of the critical configurations and the configurations connected to them. Finally, we give a geometric characterizationof the configurations that belong to the tube of typical paths for both transitions. 

\paragraph{Literature on the Potts model} All grouped citations here and henceforth are in chronological order of publication. The Potts model is one of the most studied statistical physics models, as the vast literature on the subject, both on the mathematics side and the physics side, attests. The study of the equilibrium properties of the Potts model and their dependence on $q$, have been investigated on the square lattice $\mathbb Z^d$ in  \cite{baxter1973potts,baxter1982critical}, on the triangular lattice in \cite{baxter1978triangular,enting1982triangular} and on the Bethe lattice in \cite{ananikyan1995phase,de1991metastability,di1987potts}. The mean-field version of the Potts model has been studied in \cite{costeniuc2005complete,ellis1990limit,ellis1992limit,gandolfo2010limit,wang1994solutions}. Furthermore, the tunneling behaviour for the Potts model with zero external magnetic field has been studied in \cite{nardi2019tunneling,bet2021critical,kim2021metastability}. In this energy landscape there are $q$ stable states and there is not any relevant metastable state. 
In \cite{nardi2019tunneling}, the authors derive the asymptotic behavior of the first hitting time for the transition between stable configurations, and give results in probability, in expectation and in distribution. They also characterize the behavior of the mixing time and give a lower and an upper bound for the spectral gap. In \cite{bet2021critical}, the authors study the tunneling from a stable state to the other stable configurations and between two stable states. In both cases, they geometrically identify the union of all minimal gates and the tube of typical trajectories. Finally, in \cite{kim2021metastability}, the authors study the model in dimensions two and three. They give a description of the so-called \textit{gateway configurations} in order to compute the prefactor. These gateway configurations are quite different from the minimal gates in \cite{bet2021critical}.
The $q$-Potts model with positive external magnetic field has been studied in \cite{bet2021metastabilitypos}. In this scenario there are $q-1$ multiple degenerate metastable states and a unique stable configuration. The authors answer all the three issues of the metastability introduced above for the transition from any metastable to the stable state. 

\paragraph{Literature on metastability} In this paper we adopt the framework known as \textit{pathwise approach}, which was initiated in 1984 by Cassandro, Galves, Olivieri, Vares in \cite{cassandro1984metastable} and it was further developed in \cite{olivieri1995markov,olivieri1996markov,olivieri2005large} and independently in \cite{catoni1997exit}. The pathwise approach requires a detailed knowledge of the energy landscape to give quantitative answers to the three issues of metastability in the form of ad hoc large deviations estimates. This approach was further developed in \cite{manzo2004essential,cirillo2013relaxation,cirillo2015metastability,nardi2016hitting,fernandez2015asymptotically,fernandez2016conditioned} by separating the study of the transition time and of critical configurations from that of the tube of typical trajectories. Indeed, it was recognized that the latter requires more detailed model-dependant inputs. 
The pathwise approach has been used in \cite{arous1996metastability,cirillo1998metastability,cirillo1996metastability,kotecky1994shapes,nardi1996low,neves1991critical,neves1992behavior,olivieri2005large} to tackle the three issues for Ising-like models with Glauber dynamics. Moreover, it was also used in \cite{hollander2000metastability,den2003droplet,gaudilliere2005nucleation,apollonio2021metastability,nardi2016hitting,zocca2019tunneling} to study the transition time and the gates for Ising-like and hard-core models with Kawasaki and Glauber dynamics. Finally, this method was applied to probabilistic cellular automata (parallel dynamics) in \cite{cirillo2003metastability,cirillo2008competitive,cirillo2008metastability,procacci2016probabilistic,dai2015fast}.
The so-called \emph{potential-theoretical approach} exploits a suitable Dirichlet form and spectral properties of the transition matrix to give sharp asymptotics for the hitting time. More precisely, this method estimates the leading order of the expected value of the transition time including its \textit{prefactor}, see \cite{bovier2002metastability,bovier2004metastability,bovier2016metastability,cirillo2017sum}. The potential theoretical approach was applied to find the prefactor for Ising-like models and the hard-core model in \cite{bashiri2019on,boviermanzo2002metastability,cirillo2017sum,bovier2006sharp,den2012metastability,jovanovski2017metastability,den2018metastability} for Glauber and Kawasaki dynamics and in \cite{nardi2012sharp,bet2020effect} for parallel dynamics.
Recently, other approaches have been formulated in \cite{beltran2010tunneling,beltran2012tunneling,gaudillierelandim2014} and in \cite{bianchi2016metastable} and they are particularly adapted to estimate the pre-factor when dealing with the tunnelling between two or more stable states.
\paragraph{Outline} In Section \ref{moddescr} we define the ferromagnetic $q$-state Potts model and the associated Hamiltonian. We state our main results in Section \ref{mainres}. In Section \ref{secenergylandscape} we analyse the energy landscape and give the proofs of some useful model-dependent results that are used throughout all the next sections. In Subsections \ref{proofgates} and \ref{prooftube} we give the explicit proofs of the main results on the critical configurations and on the tube of typical paths, respectively. Finally, in Section \ref{secprefactor} we compute the prefactor and refine the estimate on the expected transition time. 
\noindent We omit thoose proofs which are technically straightforward, but nevertheless lengthy. We refer the interested reader to \cite{bet2021metastabilityneg}.
\section{Model description}\label{moddescr}
In the $q$-state Potts model each spin lies on a vertex of a finite two-dimensional rectangular lattice $\Lambda=(V,E)$, where $V=\{0,\dots,K-1\}\times\{0,\dots,L-1\}$ is the vertex set and $E$ is the edge set, namely the set of the pairs of vertices whose spins interact with each other. We identify
each pair of vertices lying on opposite sides of the rectangular lattice, so that we obtain a two-dimensional torus. 
We denote by $S$ the set of spin values, i.e., $S:= \{1,\dots,q\}$ and assume $q>2$. To each vertex $v\in V$ is associated a spin value $\sigma(v)\in S$, and $\mathcal X := S^V$ denotes the set of spin
configurations.  

We denote by $\textbf{1},\dots,\textbf{q} \in \mathcal X$ those configurations in which all the vertices have spin value $1,\dots,q$, respectively.  

\noindent To each configuration $\sigma\in\mathcal X$ we associate the energy $H(\sigma)$ given by
\begin{align}\label{hamiltonianneg}
H(\sigma)
=-J \sum_{(v,w)\in E} \mathbbm{1}_{\{\sigma(v)=\sigma(w)\}}+h\sum_{u\in V} \mathbbm{1}_{\{\sigma(u)=1\}}, 
\end{align}
where $J$ is the \textit{coupling} or \textit{interation constant} and $h$ is the \textit{negative external magnetic field}. We call $h$ \textit{negative} since there is a minus in front of $H$. In this paper we consider the ferromagnetic Potts model and set  $J=1$. 

The \textit{Gibbs measure} for the $q$-state Potts model on $\Lambda$ is a probability distribution on the state space $\mathcal X$ given by
\begin{align}\label{gibbs}
\mu_\beta(\sigma):=\frac{e^{-\beta H_\text{neg}(\sigma)}}{Z},
\end{align}
where $\beta>0$ is the inverse temperature and where $Z:=\sum_{\sigma'\in\mathcal X}e^{-\beta H(\sigma')}$.

The spin system evolves according to a Glauber-type dynamics. This dynamics is described by a single-spin update Markov chain $\{X_t^\beta\}_{t\in\mathbb{N}}$ on the state space $\mathcal X$ with the following transition probabilities: for $\sigma, \sigma' \in \mathcal X$,
\begin{align}\label{metropolisTP}
P_\beta(\sigma,\sigma'):=
\begin{cases}
Q(\sigma,\sigma')e^{-\beta [H_\text{neg}(\sigma')-H_\text{neg}(\sigma)]^+}, &\text{if}\ \sigma \neq \sigma',\\
1-\sum_{\eta \neq \sigma} P_\beta (\sigma, \eta), &\text{if}\ \sigma=\sigma',
\end{cases}
\end{align}
where $[n]^+:=\max\{0,n\}$ is the positive part of $n$ and 
\begin{align}\label{Qmatrix}
Q(\sigma,\sigma'):=
\begin{cases}
\frac{1}{q|V|}, &\text{if}\ |\{v\in V: \sigma(v) \neq \sigma'(v)\}|=1,\\
0, &\text{if}\ |\{v\in V: \sigma(v) \neq \sigma'(v)\}|>1,
\end{cases}
\end{align}
for any $\sigma, \sigma' \in \mathcal X$. $Q$ is the so-called \textit{connectivity matrix} and it is symmetric and irreducible, i.e., for all $\sigma, \sigma' \in \mathcal X$, there exists a finite sequence of configurations $\omega_1,\dots,\omega_n \in \mathcal X$ such that $\omega_1=\sigma$, $\omega_n=\sigma'$ and $Q(\omega_i,\omega_{i+1})>0$ for $i=1,\dots,n-1$.  Hence, the resulting stochastic dynamics defined by \eqref{metropolisTP} is reversible with respect to the Gibbs measure \eqref{gibbs}. The triplet $(\mathcal X,H,Q)$ is called the \textit{energy landscape}. 
\noindent The dynamics defined above belongs to the class of Metropolis dynamics. In particular, at each step the update of vertex $v$ depends on the neighboring spins of $v$ and on the following energy difference
\begin{align}\label{energydifferenceneg}
H_{\text{neg}}&(\sigma^{v,s})-H_{\text{neg}}(\sigma)=\begin{cases}
\sum_{w \sim v} (\mathbbm{1}_{\{\sigma(v)=\sigma(w)\}}-\mathbbm{1}_{\{\sigma(w)=s\}})-h, &\text{if}\ \sigma(v)=1,\ s\neq1,\\
\sum_{w \sim v} (\mathbbm{1}_{\{\sigma(v)=\sigma(w)\}}-\mathbbm{1}_{\{\sigma(w)=s\}}), &\text{if}\ \sigma(v)\neq1,\ s\neq1,\\
\sum_{w \sim v} (\mathbbm{1}_{\{\sigma(v)=\sigma(w)\}}-\mathbbm{1}_{\{\sigma(w)=s\}})+h, &\text{if}\ \sigma(v)\neq1,\ s=1,
\end{cases}
\end{align}
where $\sigma^{v,s}$ is the configuration obtained from $\sigma$ by updating the spin in the vertex $v$ to $s$, i.e., $\sigma^{v,s}(w)=\sigma(w)$ if $w\neq v$, $\sigma^{v,s}(w)=s$ if $w=v$.


\section{Main results on the $q$-state Potts model with negative external magnetic field}\label{mainres}
In this section we state our main results. Note that we give the proof of the main results by considering the condition $L\ge K\ge3\ell^*$, where $\ell^*:=\left\lceil \frac{2}{h} \right\rceil$ is the \textit{critical length}. It is possible to extend the results to the case $K>L$ by interchanging the role of rows and columns in the proof.

In order to state our main results on the Potts model with Hamiltonian as in \eqref{hamiltonianneg}, we have the following assumption. 
\begin{assumption}\label{remarkconditionneg}
We assume that the following conditions are verified:

\text{(i)} the magnetic field $h$ is such that $0<h<1$;

\text{(ii)} $2/h$ is not integer.

\end{assumption}
\subsection{Energy landscape}\label{subsecenerland}
The first result that we give is the identification of the set of the global minima of the Hamiltonian \ref{hamiltonianneg}. This follows by simple algebraic calculations. 
\begin{proposition}[Identification of $\mathcal X^s_\text{neg}$]\label{stablesetnegprop}
If the external magnetic field is negative, then the set of the global minima $\mathcal X^s_\text{neg}$ of the Hamiltonian \eqref{hamiltonianneg} is given by $\mathcal X^s_{\emph{neg}}=\{\bold 2,\dots,\bold q\}$.
\end{proposition}

Next, we prove that the $q$-state Potts model with Hamiltonian $H_{\text{neg}}$ defined in \eqref{hamiltonianneg} has only one metastable state and we give an estimate of the stability level of this configuration. Formally, we call \textit{path} a finite sequence $\omega$ of configurations $\omega_0,\dots,\omega_n \in \mathcal X$, $n \in \mathbb{N}$, such that $Q(\omega_i,\omega_{i+1})>0$ for $i=0,\dots,n-1$. Let $\Omega_{\sigma,\sigma'}$ be the set of all paths between $\sigma$ and $\sigma'$.
 Given a path $\omega=(\omega_0,\dots,\omega_n)$, we define the \textit{height} of $\omega$ as
\begin{align}\label{height}
\Phi_\omega:=\max_{i=0,\dots,n} H(\omega_i).
\end{align}

\noindent For any pair $\sigma, \sigma' \in \mathcal X$, the \textit{communication height} $\Phi(\sigma,\sigma')$ between $\sigma$ and $\sigma'$ is the minimal energy across all paths $\omega:\sigma \to \sigma'$, i.e., 
\begin{align}\label{comheight}
\Phi(\sigma,\sigma'):=\min_{\omega:\sigma \to \sigma'} \Phi_\omega = \min_{\omega:\sigma \to \sigma'} \max_{\eta \in \omega} H(\eta). 
\end{align}

\noindent We define the set of \textit{optimal paths} between $\sigma, \sigma' \in\mathcal X$ as
\begin{align}\label{optpaths}
\Omega_{\sigma,\sigma'}^{opt}:=\{\omega\in\Omega_{\sigma,\sigma'}:\ \max_{\eta\in\omega} H(\eta)=\Phi(\sigma,\sigma')\}.
\end{align}

\noindent For any $\sigma \in \mathcal X$, let $\mathcal{I}_\sigma:=\{\eta \in \mathcal X:\ H(\eta)<H(\sigma)\}$ 
be the set of states with energy strictly smaller than $H(\sigma)$.
We define \textit{stability level} of $\sigma$ the energy barrier
\begin{align}\label{stabilitylevel}
V_\sigma:=\Phi(\sigma,\mathcal{I}_\sigma)-H(\sigma).
\end{align}
If $\mathcal{I}_\sigma = \varnothing$, we set $V_\sigma:=\infty$.
Finally, we define the set of \textit{metastable states} as 
\begin{align}\label{metastableset}
\mathcal X^m:=\{\eta \in \mathcal X: V_\eta=\max_{\sigma\in\mathcal X\backslash\mathcal X^s} V_\sigma\}.
\end{align}
Furthermore, for any $\sigma\in\mathcal X$ and any $\varnothing\neq\mathcal A\subset\mathcal X$, we set $\Gamma(\sigma,\mathcal A):=\Phi(\sigma,\mathcal A)-H(\sigma)$.

\noindent We refer to Figure \ref{figurenegative} for an illustration of the $4$-Potts model.
\begin{figure}[h!]
\centering
\begin{tikzpicture}[scale=0.6,transform shape]
\draw[black!55!white] (1,0.7)--(1.37,0.32);\draw[black!55!white] (1.15,0.725)--(1.51,0.35);
\fill [pattern=horizontal lines,pattern color=black!55!white] (1,.7)--(1.35,.34)--(1.49,.37)-- (1.15,.715)--(1,.7);

\draw[black!55!white] (1,.18)--(-1.4,.2);\draw[black!55!white] (1,.25)--(-1.2,.28);
\fill [pattern=vertical lines,pattern color=black!55!white] (1,.18)--(-1.4,.2)--(-1,.28)--(0.9,.25)--(1,.18);

\draw[black!55!white] (-1.4,.4)--(-0.5,.8);\draw[black!55!white] (-1.25,.35)--(-0.45,.74);
\fill [pattern=horizontal lines,pattern color=black!55!white]  (-1.35,.43)--(-0.5,.8)--(-0.45,.74)--(-1.1,.43)--(-1.25,.43);

\draw[black!40!white] (-2.1,0.3) ellipse (0.945cm and 0.15cm);
\draw[black!40!white] (1.8,0.2) ellipse (0.945cm and 0.15cm);
\draw[black!40!white] (0.3,0.8) ellipse (0.945cm and 0.15cm);
\fill[white] (0.3,0.8) ellipse (0.945cm and 0.15cm);
\fill[white]  (1.8,0.2) ellipse (0.945cm and 0.15cm);
\fill[white]  (-2.1,0.3) ellipse (0.945cm and 0.15cm);
\fill[pattern=north west lines,pattern color=black!40!white,thin] (0.3,0.8) ellipse (0.945cm and 0.15cm);
\fill[pattern=north west lines,pattern color=black!40!white,thin]  (1.8,0.2) ellipse (0.945cm and 0.15cm);
\fill[pattern=north west lines,pattern color=black!40!white,thin]  (-2.1,0.3) ellipse (0.945cm and 0.15cm);
\draw (0.3,2.85) ellipse (1.3cm and 0.15cm);
\fill [black!10!white] (0.3,2.85) ellipse (1.3cm and 0.15cm);
\draw (0.3,-1.4) parabola (-1,2.85);
\draw[dotted] (0.3,-1.4) parabola (1.6,2.85);
\draw[thick,white] (0.3,-1.4) parabola (0.966,-0.28);
\draw (0.3,-1.4) parabola (0.966,-0.28);
\draw (1.6,2.85)--(1.535,2.43);
\draw[white,thick] (-0.955,2.55)--(-0.82,1.73);
\draw[dotted] (-0.955,2.55)--(-0.82,1.73);
\fill (0.3,-1.4) circle (1.3pt); \draw (0.3,-1.4) node[below] {{$\bold 3$}};
\fill[black!3!white] (0,0.5) parabola (-1,2.5) (0,0.5) parabola (1,2.5) (0,0.5)--(-1,2.5)--(1,2.5)--(0,0.5);

\draw (0,0.5) parabola (-1,2.5);\draw (0,0.5) parabola (1,2.5);
\draw[white,thick] (0,0.5) parabola (1,2.5);\draw (0,0.5) parabola (0.64,1.32);
\draw[dotted] (0,0.5) parabola (1,2.5);
\draw[white,thick] (-1,2.5)--(-0.87,1.9);

\draw[dotted] (-1,2.5)--(-0.87,1.9);
\draw(1,2.5)--(0.98,2.4);
\draw (0,2.5) ellipse (1cm and 0.2cm);
\fill [black!10!white] (0,2.5) ellipse (1cm and 0.2cm);
\fill (0,0.5) circle (1.3pt)  node[below] {{$\bold 1$}};
\draw (1.8,-2) parabola (0.5,2.24);\draw (1.8,-2) parabola (3.1,2.24);\draw (1.8,2.24) ellipse (1.3cm and 0.2cm);
\fill [black!10!white](1.8,2.24) ellipse (1.3cm and 0.2cm);
\fill (1.8,-2) circle (1.3pt); \draw (1.8,-2) node[below] {{$\bold 4$}};
\draw (-2.1,-1.9) parabola (-0.8,2.34);\draw (-2.1,-1.9) parabola (-3.4,2.34);\draw (-2.1,2.34) ellipse (1.3cm and 0.2cm);
\fill [black!10!white](-2.1,2.34) ellipse (1.3cm and 0.2cm);
\fill (-2.1,-1.9) circle (1.3pt); \draw (-2.1,-1.9) node[below] {{$\bold 2$}};
\draw(-4.45,2.4) node {{$\Phi_{\text{neg}}(\bold 1,\mathcal X^s_{\text{neg}})$}};
\draw[black!55!white] (-4.5,.25) node {{$\Phi_{\text{neg}}(\bold s,\mathcal X^s_{\text{neg}}\backslash\{\bold s\})$}};
\fill (-0.9,2.4) circle (1.3pt);
\fill (0.75,2.36) circle (1.3pt);
\fill (0.15,2.7) circle (1.3pt);
\draw[->,ultra thin] (-0.9,2.4)--(-1.7,3);\draw (-1.8,2.9) node[above] {$\bar B_{\ell^*-1,\ell^*}^1(1,2)$};
\draw[->,ultra thin] (0.75,2.4)--(3,3);\draw (3,2.9) node[above,right] {$\bar B_{\ell^*-1,\ell^*}^1(1,4)$};
\draw[->,ultra thin] (0.15,2.7)--(1,3.7);\draw (1.1,3.6) node[above] {$\bar B_{\ell^*-1,\ell^*}^1(1,3)$};
\end{tikzpicture}
\caption{\label{figurenegative} Schematic picture of the energy landscape below $\Phi_\text{neg}(\bold 1,\mathcal X^s_\text{neg})$ of the $4$-state Potts model with negative external magnetic field with $S=\{1,2,3,4\}$, $\mathcal X^s_{\text{neg}}=\{\bold 2,\bold 3,\bold 4\}$. We have not represented the cycles (valleys) that contain configurations with stability level smaller than or equal to $2$ (see Proposition \ref{ricorrenzaneg}).
} 
\end{figure}
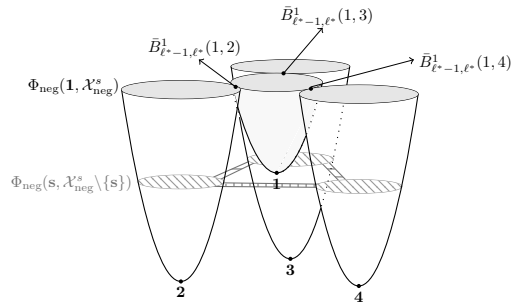\FloatBarrier
\begin{theorem}[Identification of $\mathcal X^m_\text{neg}$]\label{teometastableneg}
If the external magnetic field is negative, then
$\mathcal X^m_{\emph{neg}}=\{\bold 1\}$ 
and
\begin{align}\label{estimatestablevelmetaneg}
\Gamma^m_{\emph{neg}}=\Gamma_{\emph{neg}}(\bold 1,\mathcal X^s_{\emph{neg}})=4\ell^*-h(\ell^*(\ell^*-1)+1).
\end{align}
\end{theorem}
\textit{Proof.} 
To prove this, we apply \cite[Theorem 2.4]{cirillo2013relaxation}. The first assumption on the identification of the communication height follows by Proposition \ref{refpathneg} and Proposition \ref{lowerboundneg}. The second assumption, the estimate of the stability level of any $\sigma\in\mathcal X\backslash\{\bold 1,\dots,\bold q\}$, is proved in Proposition \ref{ricorrenzaneg}. $\qed$

In the following proposition, which we prove in Subsection \ref{submetastableneg}, we give a uniform estimate of the stability level for any configuration $\eta\in\mathcal X\backslash \{\bold 1,\dots,\bold q\}$.
\begin{proposition}[Estimate on the stability level]\label{ricorrenzaneg} 
If the external magnetic field is negative, then for any $\eta\in\mathcal X\backslash \{\bold 1,\dots,\bold q\},$ $V^\emph{neg}_\eta\le 2< \Gamma_{\emph{neg}}(\bold 1,\mathcal X^s_{\emph{neg}}).$
\end{proposition}
We define \textit{metastable set at level $V$} the set of all the configurations with stability level larger than $V$, i.e.,
\begin{align}\label{metasetV}
\mathcal X_V:=\{\sigma\in\mathcal X:V_\sigma>V\}.
\end{align} 
Moreover, given a non-empty subset $\mathcal A \subset \mathcal X$ and a configuration $\sigma \in \mathcal X$, we define
\begin{align}\label{firsthittingtime}
\tau_\mathcal A^\sigma := \text{inf}\{t>0: \ X_t^\beta \in \mathcal A\}
\end{align}
as the \textit{first hitting time} of the subset $\mathcal A$ for the Markov chain $\{X_t^\beta\}_{t \in \mathbb{N}}$ starting from $\sigma$ at time $t=0$. 
Exploiting the estimate of the stability level in Proposition \ref{ricorrenzaneg}, we obtain the following result on a recurrence property to metastable and stable states, i.e., $\{\bold 1,\dots,\bold q\}$. 
\begin{theorem}[Recurrence property]\label{teorecprop}
If the external magnetic field is negative, then for any  $\sigma\in\mathcal X$ and for any $\epsilon>0$ there exists $k>0$ such that for $\beta$ sufficiently large 
\begin{align}
\mathbb P(\tau^\sigma_{\{\bold 1,\dots,\bold q\}}>e^{\beta(2+\epsilon)})\le e^{-e^{k\beta}}.
\end{align}
\end{theorem}
\textit{Proof}. Apply \cite[Theorem 3.1]{manzo2004essential} with $V=2$ and use \eqref{metasetV} and Proposition \ref{ricorrenzaneg} to get $\mathcal X_2=\{\bold 1,\dots,\bold q\}=\mathcal X^s_{\text{neg}}\cup\mathcal X^m_{\text{neg}}$, where the last equality follows by Proposition \ref{stablesetnegprop} and Theorem \ref{teometastableneg}. $\qed$

From Theorem \ref{teorecprop} follows that the function $\beta\to f(\beta) := \mathbb P(\tau^\sigma_{\{\bold 1,\dots,\bold q\}}>e^{\beta(2+\epsilon)})$ satisfies $\lim_{\beta\to\infty} \frac{\log f(\beta)}{\beta}=-\infty$ and such a function is known as \textit{super-exponentially small}.

From Proposition \ref{stablesetnegprop}, we have that when $q>2$ the energy landscape $(\mathcal X,H_\text{neg},Q)$ has multiple stable states. We are interested in studying the transition from the metastable state $\bold 1$ to $\mathcal X^s_\text{neg}$ and also the transition from $\bold 1$ to a fixed stable configuration $\bold s\in\mathcal X^s_\text{neg}$. To this end, it is useful to compare the communication energy between two different stable states and the communication energy between the metastable state and a stable configuration. Furthermore, for any $\bold s\in\mathcal X^s_\text{neg}$ in order to find the asymptotic upper bound in probability for $\tau_\bold s^\bold 1$, we estimate the maximum energy barrier that the process started from $\bold r\in\mathcal X^s_\text{neg}\backslash\{\bold s\}$ has to overcome so as to reach $\bold s$ (in Theorem \ref{timenegative}). These are the goals of the following theorem, for whose proof we refer to \cite[Theorem 4.3]{bet2021metastabilityneg}. 

In order to state the next result, we need some further definitions. 
A non-empty subset $\mathcal C\subset\mathcal X$ is called \textit{cycle} if it is either a singleton or a connected set such that
\begin{align}\label{cycle}
\max_{\sigma\in\mathcal C} H(\sigma)<H(\mathscr{F}(\partial\mathcal C)).
\end{align}
When $\mathcal C$ is a singleton, it is said to be a \textit{trivial cycle}. Let $\mathscr C(\mathcal X)$ be the set of cycles of $\mathcal X$. 

\noindent The \textit{depth} of a cycle $\mathcal C$ is given by
\begin{align}\label{defdepcycle}
\Gamma(\mathcal C):=H(\mathscr F(\partial\mathcal C))-H(\mathscr F(\mathcal C)).
\end{align}
If $\mathcal C$ is a trivial cycle we set $\Gamma(\mathcal C)=0$.

\noindent Given a non-empty set $\mathcal A\subset\mathcal X$, we denote by $\mathcal M(\mathcal A)$ the \textit{collection of maximal cycles} $\mathcal A$, i.e., 
$\mathcal M(\mathcal A):=\{\mathcal C\in\mathscr C(\mathcal X)|\ \mathcal C$  maximal by inclusion under constraint $\mathcal C\subseteq\mathcal A\}.$
For any $\mathcal A\subset\mathcal X$, we define the \textit{maximum depth of $\mathcal A$} as the maximum depth of a cycle contained in $\mathcal A$, i.e., 
\begin{align}\label{defgammatildenzb}
\widetilde{\Gamma}(\mathcal A):=\max_{\mathcal C\in\mathcal M(\mathcal A)} \Gamma(\mathcal C).
\end{align}
In \cite[Lemma 3.6]{nardi2016hitting} the authors give an alternative characterization of \eqref{defgammatildenzb} as the maximum initial energy barrier that the process started from a configuration $\eta\in\mathcal A$ possibly has to overcome to exit from $\mathcal A$, i.e., $\widetilde{\Gamma}(\mathcal A)=\max_{\eta\in\mathcal A} \Gamma(\eta,\mathcal X\backslash\mathcal A)$. 

\noindent Finally, for any $\sigma\in\mathcal X$, if $\mathcal{A}$ is a non-empty target set, we define the \textit{initial cycle} for the transition from $\sigma$ to $\mathcal A$ as
$\mathcal C_{\mathcal{A}}^\sigma(\Gamma):=\{\sigma\}\cup\{\eta\in\mathcal X:\ \Phi(\sigma,\eta)-H(\sigma)<\Gamma=\Phi(\sigma,\mathcal A)-H(\sigma)\}.$
Note that if $\sigma\notin\mathcal{A}$, then $C_{\mathcal{A}}^\sigma(\Gamma)\cap{\mathcal{A}}=\varnothing$.
\begin{theorem}\label{theoremcomparisonneg}
Consider the $q$-state Potts model on a $K\times L$ grid $\Lambda$, with periodic boundary conditions and with negative external magnetic field. For any $\bold s\in\mathcal X^s_{\emph{neg}}$, we have
\begin{align}\label{comparisonalignneg}
&\Phi_{\emph{neg}}(\bold 1,\mathcal X^s_{\emph{neg}})>\Phi_{\emph{neg}}(\bold s,\mathcal X^s_{\emph{neg}}\backslash\{\bold s\}),\\
\label{comparisongammaneg}
&\Gamma_{\emph{neg}}(\bold 1,\mathcal X^s_{\emph{neg}}) < \Gamma_{\emph{neg}}(\bold s, \mathcal X^s_{\emph{neg}}\backslash\{\bold s\}),\\
\label{gammatildenegtime}
&\widetilde{\Gamma}_\emph{neg}(\mathcal X\backslash\{\bold s\})=\Gamma_{\emph{neg}}(\bold r, \mathcal X^s_{\emph{neg}}\backslash\{\bold r\}),\ \text{with $\bold r\in\mathcal X^s_\emph{neg}$.}
\end{align}
\end{theorem}
We refer to Figure \ref{figurenegative} and  Figure \ref{sidenegative} for a schematic representation of the energy landscape  and of the quantities of Theorem \ref{theoremcomparisonneg} for the $4$-state Potts model with negative magnetic field.
\begin{figure}
\centering
\begin{tikzpicture}[scale=0.5,transform shape]
\fill[black!8!white] (0,0.7) parabola (-1,2.6) (0,0.7) parabola (1,2.6) (0,0.7)--(-1,2.6)--(1,2.6)--(0,0.7);
\draw[thick,dotted](-5.1,2.6)--(4.8,2.6) node[right]{$\Phi_{\text{neg}}(\bold 1,\mathcal X^s_{\text{neg}})$};
\draw(-5,-.5) node[black!50!white,left]{$\Phi_{\text{neg}}(\bold s,\mathcal X^s_{\text{neg}}\backslash\{\bold s\})$};
\draw [black!50!white,<-] (-4.75,0.32)--(-4.75,-2.4);
\draw[black!50!white,dashed] (-5.1,.35)--(4.8,.35);
\draw[black!50!white,dotted](-5.1,-2)--(4.8,-2);
\draw[black!50!white,<->] (2.7,-1.95)--(2.7,0.32); \draw (2.7,-0.65) node[black!50!white,right] {$\Gamma_{\text{neg}}(\bold s,\mathcal X^s_{\text{neg}}\backslash\{\bold s\})$};
\draw[<->] (0,0.75)--(0,2.52); \draw (-0.05,1.95) node[right] {$\Gamma^m_\text{neg}$};
\draw (0,0.7) parabola (-1,2.5);\draw (0,0.7) parabola (1,2.5);
\draw (2.7,-2) parabola (1.4,1.5);\draw (2.7,-2) parabola (4.1,2.6);
\draw  (1.1,2.6) parabola (1.4,1.5);\draw(1.1,2.6)parabola(1,2.5);
\draw (-2.7,-2) parabola (-1.4,1.5);\draw (-2.7,-2) parabola (-4.1,2.6);
\draw  (-1.1,2.6) parabola (-1.4,1.5);\draw(-1.1,2.6)parabola(-1,2.5);

\fill (-2.7,-2) circle (1.3pt) node[below,black] {$\bold 2$};
\fill (0,0.7) circle (1.3pt) node[below] {$\bold 1$};
\fill (2.7,-2) circle (1.3pt) node[below,black] {$\bold 3$};

\fill (1.1,2.6) circle (1.3pt); \draw(1.6,2.6) node[above] {$\bar B_{\ell^*-1,\ell^*}^1(1,3)$};
\fill (-1.1,2.6) circle (1.3pt);\draw (-1,2.6) node[above] {$\bar B_{\ell^*-1,\ell^*}^1(1,2)$};
\end{tikzpicture}\ \ \ \ \ \
\begin{tikzpicture}[scale=0.5,transform shape]
\draw (0,0) circle (0.7cm);
\fill[black!10!white] (0,0) circle (0.7cm);
\fill[white] (-1.5,-0.2) -- (-0.58,1.52)--(-0.47,1.33)--(-1.3,-0.23)--(-1.5,-0.2);\draw[black!45!white,thin] (-1.5,-0.2) -- (-0.58,1.52)(-1.3,-0.23)--(-0.47,1.33);
\fill[pattern=horizontal lines,pattern color=black!45!white] (-1.5,-0.2) -- (-0.58,1.52)--(-0.47,1.33)--(-1.3,-0.23)--(-1.5,-0.2);
\draw[black!45!white,thin] (1.5,-0.2) -- (0.58,1.52)(1.3,-0.23)--(0.47,1.33);
\fill[pattern=horizontal lines,pattern color=black!45!white] (1.5,-0.2) -- (0.58,1.52)--(0.47,1.33)--(1.3,-0.23)--(1.5,-0.2);
\draw[black!45!white,thin] (-1.16,-1.3) -- (1.16,-1.3) (-1,-1.13)--(1,-1.13);
\fill[pattern=vertical lines,pattern color=black!45!white](-1.16,-1.3) -- (1.16,-1.3)--(1,-1.13)-- (-1,-1.13)-- (-1.16,-1.3);
\draw(0,1.7) circle (1cm);\draw(1.5,-0.8) circle (1cm);\draw (-1.5,-0.8) circle (1cm);
\fill[black!10!white] (0,1.7) circle (1cm);
\fill[black!10!white] (1.5,-0.8) circle (1cm);
\fill[black!10!white] (-1.5,-0.8) circle (1cm);
\draw[black!45!white,dashed](0,1.7) circle (0.6cm);\draw[black!45!white,dashed](1.5,-0.8) circle (0.6cm);\draw[black!45!white,dashed] (-1.5,-0.8) circle (0.6cm);
\fill[pattern=north west lines,pattern color=black!45!white] (0,1.7) circle (0.6cm);
\fill[pattern=north west lines,pattern color=black!45!white](1.5,-0.8) circle (0.6cm);
\fill[pattern=north west lines,pattern color=black!45!white] (-1.5,-0.8) circle (0.6cm);
\fill (0,1.7) circle (1pt) node[below] {{$\bold 3$}};
\fill (1.5,-0.8) circle (1pt) node[below] {{$\bold 4$}};
\fill (-1.5,-0.8) circle (1pt) node[below] {{$\bold 2$}};
\fill (0,0) circle (1pt) node[below] {{$\bold 1$}};

\fill (-0.6,-0.38) circle (1.6pt);
\fill (0.6,-0.38) circle (1.6pt);
\fill (0,0.7) circle (1.6pt);
\draw[->,ultra thin] (-0.6,-0.38)--(-1.1,-2);\draw(-1.25,-1.9) node[below] {$\bar B_{\ell^*-1,\ell^*}^1(1,2)$};
\draw[->,ultra thin] (0.6,-0.38)--(1.1,-2);\draw(1.25,-1.9) node[below] {$\bar B_{\ell^*-1,\ell^*}^1(1,4)$};
\draw[->,ultra thin] (0,0.7)--(-1.4,0.7) node[left] {$\bar B_{\ell^*-1,\ell^*}^1(1,3)$};
\end{tikzpicture}
\caption{\label{sidenegative} On the left, we give a side view (vertical section) of the energy landscape depicted in Figure \ref{figurenegative}. We colour light gray the initial cycle $\mathcal C^\bold1_{\mathcal X^s_\text{neg}}(\Gamma^m_\text{neg})$. On the right, viewpoint from above of the energy landscape depicted in Figure \ref{figurenegative} cut to the energy level $\Phi_{\text{neg}}(\bold 1,\bold s)$, for some $\bold s\in\mathcal X^s_{\text{neg}}$.  The dashed part denotes the energy landscape whose energy value is smaller than $\Phi_{\text{neg}}(\bold 1,\bold s)$. The cycles whose bottom is a stable state are deeper than the cycle $\mathcal C^\bold 1_{\mathcal X^s_\text{neg}}(\Gamma_\text{neg}^m)$ of the metastable state, hence we depict them with circles whose diameter is larger than the one related to the metastable state $\bold 1$.}
\end{figure}
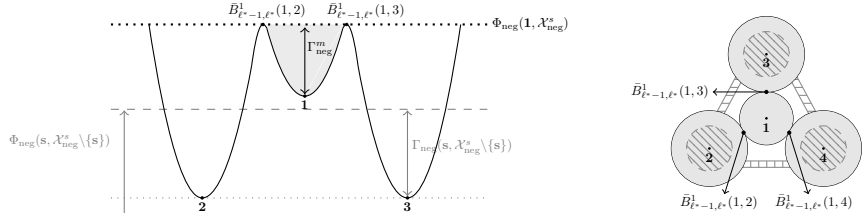\FloatBarrier

\subsection{Asymptotic behavior of $\tau_{\mathcal X^s_{\text{neg}}}^{\textbf{1}}$ and $\tau_{\bold s}^{\textbf{1}}$ and mixing time}\label{mainrestime}
In the following theorem we give asymptotic bounds in probability for both $\tau_{\mathcal X^s_{\text{neg}}}^{\textbf{1}}$ and $\tau_{\bold s}^{\textbf{1}}$, identify the order of magnitude of the expected value of $\tau_{\mathcal X^s_{\text{neg}}}^{\textbf{1}}$ and prove that the asymptotic rescaled distibution of $\tau_{\mathcal X^s_{\text{neg}}}^{\textbf{1}}$ is exponential. Furthermore, we also identify the mixing time and give an upper and a lower bound for the spectral gap. Formally, let $\{X_t^\beta\}_{t\in\mathbb N}$ be the Markov chain with transition probabilities \eqref{metropolisTP} and stationary distribution \eqref{gibbs}. For every $\epsilon\in(0,1)$, we define the \textit{mixing time} $t^{\text{mix}}_\beta(\epsilon)$ by
\begin{align}\label{mixingtimedef}
t^{\text{mix}}_\beta(\epsilon):=\min\{n\ge 0|\ \max_{\sigma\in\mathcal X}||P_\beta^n(\sigma,\cdot)-\mu_\beta(\cdot)||_{\text{TV}}\le\epsilon\},
\end{align}
where the total variance distance is defined by $||\nu-\nu'||_{\text{TV}}:=\frac 1 2 \sum_{\sigma\in\mathcal X}|\nu(\sigma)-\nu'(\sigma)|$ for every two probability distribution $\nu,\nu'$ on $\mathcal X$. Furthermore, we define \textit{spectral gap} as
\begin{align}\label{spectralgapdef}
\rho_\beta:=1-\lambda_\beta^{(2)},
\end{align}
where $1=\lambda_\beta^{(1)}>\lambda_\beta^{(2)}\ge\dots\ge\lambda_\beta^{(|\mathcal X|)}\ge-1$ are the eigenvalues of the matrix $P_\beta(\sigma,\eta))_{\sigma,\eta\in\mathcal X}$.
\begin{theorem}[Asymptotic behavior of $\tau_{\mathcal X^s_{\text{neg}}}^{\textbf{1}}$ and $\tau_{\bold s}^{\textbf{1}}$ and mixing time]\label{timenegative}
If the external magnetic field is negative, then for any $\bold s \in \mathcal X^s_{\emph{neg}}$, the following statements hold:
\begin{itemize}
\item[\emph{(a)}] for any $\epsilon>0$, $\lim_{\beta\to\infty}\mathbb P_\beta(e^{\beta(\Gamma^m_{\emph{neg}}-\epsilon)}<\tau_{\mathcal X^s_{\emph{neg}}}^{\bold 1}<e^{\beta(\Gamma^m_{\emph{neg}}+\epsilon)})=1$;
\item[\emph{(b)}] for any $\epsilon>0$, $\lim_{\beta\to\infty}\mathbb P_\beta(e^{\beta(\Gamma^m_{\emph{neg}}-\epsilon)}<\tau_{\bold s}^{\bold 1}<e^{\beta(\Gamma_{\emph{neg}}(\bold s,\mathcal X^s_\emph{neg}\backslash\{\bold s\})+\epsilon)})=1$;
\item[\emph{(c)}] $\lim_{\beta\to\infty} \frac{1}{\beta}\log\mathbb E[\tau_{\mathcal X^s_{\emph{neg}}}^\bold 1]=\Gamma^m_{\emph{neg}}$;
\item[\emph{(d)}] $\frac{\tau_{\mathcal X^s_{\emph{neg}}}^\bold 1}{\mathbb E[\tau_{\mathcal X^s_{\emph{neg}}}^\bold 1]}\xrightarrow{d}$ \emph{Exp(1)};
\item[\emph{(e)}] for every $\epsilon\in(0,1)$ and $\bold s\in\mathcal X^s_\emph{neg}$, $\lim_{\beta\to\infty}\frac 1 \beta \log t^{\text{mix}}_\beta(\epsilon)=\Gamma_\emph{neg}(\bold s,\mathcal X^s_\emph{neg}\backslash\{\bold s\})$ and there exist two constants $0<c_1\le c_2<\infty$ independent of $\beta$ such that, for any $\beta>0$, $c_1e^{-\beta\Gamma_\emph{neg}(\bold s,\mathcal X^s_\emph{neg}\backslash\{\bold s\})}\le\rho_\beta\le c_2e^{-\beta\Gamma_\emph{neg}(\bold s,\mathcal X^s_\emph{neg}\backslash\{\bold s\})}$.
\end{itemize}
\end{theorem}
\textit{Proof}. Item (a) holds in view of Theorem \ref{teometastableneg} and \cite[Theorem 4.1]{manzo2004essential}. The lower bound of item (b) follows by Theorem \ref{teometastableneg} and \cite[Propositions 3.4]{nardi2016hitting}, while the upper bound by \eqref{gammatildenegtime} and \cite[Propositions 3.7]{nardi2016hitting}. 
Item (c) follows from Theorem \ref{teometastableneg} and \cite[Theorem 4.9]{manzo2004essential}.
Lastly, item (d), i.e., the asymptotic exponentiality of $\tau_{\mathcal X^s_{\text{neg}}}^\bold 1$, follows from Theorem \ref{gammatildenegtime} and \cite[Theorem 4.15]{manzo2004essential}. For this last item, we refer also to \cite[Theorem 3.19, Example 3]{nardi2016hitting}. Item (e) follows by \eqref{gammatildenegtime} and by \cite[Proposition 3.24]{nardi2016hitting}. $\qed$

\begin{remark}\label{labelintrorecall}
Note that the lower and upper bounds for $\tau_\bold s^\bold 1$ in item (b) have different exponents. Indeed, the presence of a subset of the optimal paths, that the process follows with probability strictly positive, going from $\bold 1$ to $\bold s$ without crossing $\mathcal X^s_\text{neg}\backslash\{\bold s\}$, implies that the lower bound is sharp. Moreover, the presence of a subset of the optimal paths going from $\bold 1$ to $\bold s$ crossing $\mathcal X^s_\text{neg}\backslash\{\bold s\}$, ensures that  the process, with probability strictly positive, enters at least a cycle $\mathcal C^\bold r_{\bold s}(\Gamma_\text{neg}(\bold r,\bold s))$ for any given $\bold r\in\mathcal X^s_\text{neg}\backslash\{\bold s\}$ which is deeper than the initial cycle $\mathcal C^\bold 1_{\bold s}(\Gamma^m_\text{neg})$. This implies that the maximum depth of the cycles crossed by these paths is $\Gamma_\text{neg}(\bold r,\bold s)$, thus the upper is sharp. 
Finally, we remark that in \cite[Theorem 4.3]{bet2021metastabilitypos} items (a) and (b) coincide since in that scenario there is a unique stable state.
\end{remark}

\subsection{Minimal gates for the metastable transitions}\label{mainresgates}
We also identify the set of minimal gates for the transition $\bold 1\to\mathcal X^s_{\text{neg}}$ and also for the transition $\bold 1\to\bold s$ for some fixed $\bold s\in\mathcal X^s_\text{neg}$. To this end, we need some further definitions.
The set of \textit{minimal saddles} between $\sigma, \sigma' \in \mathcal X$ is defined as
\begin{align}\label{saddles}
\mathcal S(\sigma,\sigma'):=\{\xi\in\mathcal X:\exists\omega\in\Omega_{\sigma,\sigma'}^{opt},\ \xi\in\omega:\ \max_{\eta\in\omega} H(\eta)=H(\xi)\}.
\end{align}
We say that $\eta\in\mathcal S(\sigma,\sigma')$ is an \textit{essential saddle} if there exists $\omega\in\Omega_{\sigma,\sigma'}^{opt}$ such that either
\begin{itemize}
\item  $\{\text{arg max}_\omega H\}=\{\eta\}$ or
\item $\{\text{arg max}_\omega H\}\supset\{\eta\}$ and $\{\text{arg max}_{\omega'} H\}\not\subseteq\{\text{arg max}_\omega H\}\backslash \{\eta\}$ for all $\omega'\in\Omega_{\sigma,\sigma'}^{opt}$.
\end{itemize}
A saddle $\eta\in\mathcal S(\sigma,\sigma')$ that is not essential is said to be \textit{unessential}. 

\noindent Given $\sigma, \sigma' \in \mathcal X$, we say that $\mathcal W(\sigma,\sigma')$ is a \textit{gate} for the transition from $\sigma$ to $\sigma'$ if  $\mathcal W(\sigma,\sigma')\subseteq\mathcal S(\sigma,\sigma')$ and $\omega\cap\mathcal W(\sigma,\sigma')\neq\varnothing$ for all $\omega\in\Omega_{\sigma,\sigma'}^{opt}$.
We say that  $\mathcal W(\sigma,\sigma')$ is a \textit{minimal gate} for the transition from $\sigma$ to $\sigma'$ if it is a minimal (by inclusion) subset of $\mathcal S(\sigma,\sigma')$ that is visited by all optimal paths. More in detail, it is a gate and for any $\mathcal W'\subset\mathcal W(\sigma,\sigma')$ there exists $\omega'\in\Omega_{\sigma,\sigma'}^{opt}$ such that $\omega'\cap\mathcal W'=\varnothing$. We denote by $\mathcal{G}=\mathcal{G}(\sigma,\sigma')$ the union of all minimal gates for the transition $\sigma\to\sigma'$.

\noindent In our scenario, we define
\begin{align}\label{gate1neg}
&\mathcal W_{\text{neg}}(\bold 1,\mathcal X^s_{\text{neg}}):=\bigcup_{t=2}^q \bar B_{\ell^*-1,\ell^*}^1(1,t)
\ \ \text{and}\ \ \mathcal{W}_{\text{neg}}'(\bold 1,\mathcal X^s_{\text{neg}}):=\bigcup_{t=2}^q \bar B_{\ell^*,\ell^*-1}^1(1,t),
\end{align}
where $\bar B_{a,b}^l(r, s)$ denotes the set of those configurations in which all the vertices have spins $r$, except those, which have spins $s$, in a rectangle $a\times b$ with a bar $1\times l$ adjacent to one of the sides of length $b$, with $1\le l\le b-1$.

\noindent We refer to Figure \ref{figureexample}(b)--(c) for an example of configurations belonging respectively to $\mathcal{W}_{\text{neg}}'(\bold 1,\mathcal X^s_{\text{neg}})$ and to $\mathcal{W}_{\text{neg}}(\bold 1,\mathcal X^s_{\text{neg}})$. These sets are investigated in Subsection \ref{secmingatesneg} where we study the gate for the transition from the metastable state $\bold 1$ to $\mathcal X^s_{\text{neg}}$. 
Similarly to the Ising Model, see for instance \cite[Section 5.4]{manzo2004essential}, \cite[Section 7.4]{olivieri2005large}, if the Assumption \ref{remarkconditionneg}(ii) holds, then there exists only one minimal gate. 
\begin{theorem}[Minimal gates for the transition $\bold 1\to\mathcal X^s_{\text{neg}}$]\label{teogatenegset}
If the external magnetic field is negative, then $\mathcal W_{\emph{neg}}(\bold 1,\mathcal X^s_{\emph{neg}})$ is a minimal gate for the transition from the metastable state $\bold 1$ to $\mathcal X^s_{\emph{neg}}$. Moreover,
\begin{align}\label{mingateneg}
\mathcal G_{\emph{neg}}(\bold 1,\mathcal X^s_{\emph{neg}})=\mathcal W_{\emph{neg}}(\bold 1,\mathcal X^s_{\emph{neg}}).
\end{align}
\end{theorem}
The proof of the theorem above is given in Subsection \ref{proofgates}. We refer to Figures \ref{figurenegative} and 
\ref{sidenegative} for illustrations of the set $\mathcal W_{\text{neg}}(\bold 1,\mathcal X^s_{\text{neg}})$ when $q=4$.

Finally, in Theorem \ref{teogatenegtarget} we establish the set of all minimal gates for the transition from the metastable state $\bold 1$ to a fixed stable configuration $\bold s\in\mathcal X^s_{\text{neg}}$. In particular, starting from $\bold 1$ the process may visit some stable states in $\mathcal X^s_{\text{neg}}\backslash\{\bold s\}$ before hitting $\bold s$.  Thanks to Theorem \ref{theoremcomparisonneg}, we get that 
along any optimal path between two different stable states the process only visits states with energy value lower than $\Phi_{\text{neg}}(\bold 1,\mathcal X^s_{\text{neg}})$ and so it does not visit any other gate. See for instance Figure \ref{figurenegative} and Figure \ref{sidenegative}, where we indicate with a dashed gray line the communication energy among the stable states.
\begin{theorem}[Minimal gates for the transition $\bold 1\to\bold s\in\mathcal X^s_\text{neg}$]\label{teogatenegtarget}
If the external magnetic field is negative, then for any $\bold s\in\mathcal X^s_{\emph{neg}}$, $\mathcal W_{\emph{neg}}(\bold 1,\bold s)\equiv\mathcal W_{\emph{neg}}(\bold 1,\mathcal X^s_{\emph{neg}})$ is a minimal gate for the transition from the metastable state $\bold 1$ to $\bold s$ and
\begin{align}\label{mingatenegtarget}
\mathcal G_{\emph{neg}}(\bold 1,\bold s)\equiv\mathcal G_{\emph{neg}}(\bold 1,\mathcal X^s_{\emph{neg}}).
\end{align}
\end{theorem}
We defer to Subsection \ref{proofgates} for the proof of the theorem above. Finally, in the next corollary we prove that in both the transitions, i.e., $\bold 1\to\mathcal X^s_\text{neg}$ and $\bold 1\to\bold s\in\mathcal X^s_\text{neg}$, the process typically intersects the gates.
\begin{cor}\label{corneg}
If the external magnetic field is negative, then for any $\bold s\in\mathcal X^s_{\emph{neg}}$,
\begin{itemize}
\item[\emph{(a)}] $\lim_{\beta\to\infty} \mathbb P_\beta(\tau^\bold 1_{\mathcal W_{\emph{neg}}(\bold 1,\mathcal X^s_{\emph{neg}})}<\tau^\bold 1_{\mathcal X^s_{\emph{neg}}})=1$;
\item[\emph{(b)}] $\lim_{\beta\to\infty} \mathbb P_\beta(\tau^\bold 1_{\mathcal W_{\emph{neg}}(\bold 1,\bold s)}<\tau^\bold 1_{\bold s})=1$.
\end{itemize}
\end{cor}
\textit{Proof}. The corollary follows from Theorems \ref{teogatenegset} and \ref{teogatenegtarget} and from \cite[Theorem 5.4]{manzo2004essential}. $\qed$

\

Finally, we remark that in \cite[Theorem 4.5]{bet2021metastabilitypos} the authors identify the union of all minimal gates for the metastable transition for the $q$-Potts model with \textit{positive} external magnetic field. These minimal gates have the same geometric definition of those of our scenario, the main difference is that in the positive scenario the spins inside the quasi-square union a unit protuberance and in the sea are fixed, while in the negative case we have to take the union on all $t\in S\backslash\{1\}$, see \eqref{gate1neg}.
\subsection{Sharp estimate on the mean transition time}\label{mainresprefactor}
Exploiting the model-dependent results given in Subsections \ref{subsecenerland} and  \ref{mainresgates} and some model-independent results by \cite{bovier2002metastability,bovier2016metastability} and from \cite{baldassarri2021metastability}, in Subsection \ref{proofprefactor} we prove the following theorem in which we refine the result of Theorem \ref{timenegative}(c) by identifying the precise scaling of the prefactor multipling the exponential.
\begin{theorem}[Mean crossover time]\label{teoprefneg}
If the external magnetic field is negative, then there exists a constant $K_\emph{neg}\in(0,\infty)$ such that 
\begin{align}\label{meantimeKneg}
\lim_{\beta\to\infty} e^{-\beta\Gamma^m_\emph{neg}}\mathbb E_\bold 1(\tau_{\mathcal X^s_\emph{neg}})=K_\emph{neg}.
\end{align}
In particular, the constant $K_\emph{neg}$ is the so-called prefactor and it is given by
\begin{align}\label{Knegative}
K_\emph{neg}=\frac 3 4\frac 1 {2\ell^*-1} \frac 1 {q-1}\frac 1 {|\Lambda|}.
\end{align}
\end{theorem}
\begin{remark}
In order to prove Theorem \ref{timenegative}(c) the only model-independent inputs are the identification of $\mathcal X^m_\text{neg}$, the recurrence property given in Theorem \ref{teorecprop}, and the computation of the energy barrier $\Gamma_\text{neg}(\mathbf 1,\mathcal X^s_\text{neg})$ for the transition from the metastable state to the stable configurations, see \eqref{estimatestablevelmetaneg}. On the other hand, in order to prove Theorem \ref{teoprefneg} we need of a more accurate knowledge of the energy landscape. Indeed, it is necessary to know the geometrical identification of the critical configurations and of the configurations connected to them by a single step of the dynamics for the transition $\mathbf 1\to\mathcal X^s_\text{neg}$, that we give in Theorem \ref{teogatenegset}. 
\end{remark}

\subsection{Tube of typical trajectories of the metastable transitions}\label{mainrestubeneg}

In this subsection we give the results on the tube of typical trajectories $\mathfrak T_{\mathcal X^s_\text{neg}}(\mathbf 1)$ and $\mathfrak T_{\mathbf s}(\mathbf 1)$ for both the transitions $\mathbf 1\to\mathcal X^s_\text{neg}$ and $\mathbf1\to\mathbf s$ for any fixed $\mathbf s\in\mathcal X^s_\text{neg}$. The tube $\mathfrak T_{\mathcal X^s_\text{neg}}(\mathbf 1)$ (resp. $\mathfrak T_{\mathbf s}(\mathbf 1)$) can be characterized, and indeed identified, by only relying on the geometrical structure of the energy landscape. Once this is done it follows from standard model-independent considerations \cite{olivieri2005large, nardi2016hitting} that the dynamics leaves $\mathfrak T_{\mathcal X^s_\text{neg}}(\mathbf 1)$ (resp. $\mathfrak T_{\mathbf s}(\mathbf 1)$) through its non-principal boundary before reaching $\mathbf 1$ with exponentially small probability. In particular, the non-principal boundary are all those configurations on the boundary that do not minimize the energy. From this follows that $\mathfrak T_{\mathcal X^s_\text{neg}}(\mathbf 1)$ (resp. $\mathfrak T_{\mathbf s}(\mathbf 1)$) contains those configurations which are visited with positive probability before hitting $\mathcal X^s_\text{neg}$ (resp. $\mathbf s$) as $\beta\to\infty$. Formally, for any $\mathcal C\in\mathscr C(\mathcal X)$, we define as
\begin{align}\label{principalboundary}
\mathcal B(\mathcal C):=
\begin{cases}
\mathscr F(\partial\mathcal C) &\text{if}\ \mathcal C\ \text{is a non-trivial cycle},\\
\{\eta\in\partial\mathcal C: H(\eta)< H(\sigma)\} &\text{if}\ \mathcal C=\{\sigma\}\ \text{is a trivial cycle},
\end{cases}
\end{align}
the \textit{principal boundary} of $\mathcal C$. Furthermore, let $\partial^{np}\mathcal C$ be the \textit{non-principal boundary} of $\mathcal C$, i.e., $\partial^{np}\mathcal C:=\partial\mathcal C\backslash\mathcal B(\mathcal C).$

The tube is defined in terms of unions of $\bar B_{a,b}^l(r, s)$, defined below Theorem \ref{gate1neg}, and of the following sets.

\noindent - $\bar R_{a,b}(r, s)$ denotes the set of those configurations in which all the vertices have spins equal to $r$, except those, which have spins $s$, in a rectangle $a\times b$. Note that when either $a=L$ or $b=K$,  $\bar R_{a,b}(r, s)$ contains those configurations which have an $r$-strip and an $s$-strip. In particular, a configuration $\sigma$ has an \textit{$s$-strip} if it has a cluster of spins $s$ which is a rectangle that wraps around $\Lambda$. For any $r,s\in S$, we say that an $s$-strip is \textit{adjacent} to an $r$-strip if they are at lattice distance one from each other. For instance, in Figure \ref{figurerecurrenceproplabel}(a) there are depicted vertical adjacent strips.

\noindent - For any $s\neq 1$, we define
\begin{align}\label{setverticalstrip}
&\hspace{-18pt}\mathscr S_\text{neg}^v(1,s)\hspace{-2pt}:=\hspace{-2pt}\{\sigma\hspace{-2pt}\in\hspace{-2pt}\mathcal X(1,s)\hspace{-2pt}:\hspace{-2pt}\sigma\ \text{has a vertical $s$-strip of thickness at least $\ell^*$ with}\notag\\&\hspace{-15pt}\text{possibly a bar of length $l\hspace{-2pt}=\hspace{-2pt}1,\hspace{-2pt}...,\hspace{-2pt}K\hspace{-2pt}$ on one of the two vertical edges}\},\\ 
\label{sethorizontalstrip}
&\hspace{-18pt}\mathscr S_\text{neg}^h(1,s)\hspace{-2pt}:=\hspace{-2pt}\{\sigma\hspace{-2pt}\in\hspace{-2pt}\mathcal X(1,s)\hspace{-2pt}:\hspace{-2pt}\sigma\ \text{has a horizontal $s$-strip of thickness at least $\ell^*$}\notag\\&\hspace{-15pt}\text{with possibly a bar of length $l\hspace{-2pt}=\hspace{-2pt}1,\hspace{-2pt}...,\hspace{-2pt}L\hspace{-2pt}$ on one of the two horizontal edges}\},
\end{align} 
where $\mathcal X(r,s)=\{\sigma\in\mathcal X: \sigma(v)\in\{r,s\}\ \text{for any}\ v\in V\}$.

\begin{theorem}[Tube of typical paths for the transition $\mathbf 1\to\mathcal X^s_\text{neg}$]\label{teotubesetneg}
If the external magnetic field is negative, then the tube of typical trajectories for the transition $\mathbf 1\to\mathcal X^s_\emph{neg}$ is \begin{align}\label{tesiteotubenegxs1}
&\mathfrak T_{\mathcal X^s_\text{neg}}(\mathbf 1)\hspace{-2pt}:=\hspace{-2pt}\bigcup_{s=2}^q\hspace{-2pt}\biggr[\bigcup_{\ell=1}^{\ell^*-1}\hspace{-2pt}\bar R_{\ell-1,\ell}(1,s)\hspace{-2pt}\cup\hspace{-2pt}\bigcup_{\ell=1}^{\ell^*}\hspace{-2pt}\bar R_{\ell-1,\ell-1}(1,s)\hspace{-2pt}\cup\hspace{-3pt}\bigcup_{\ell=1}^{\ell^*-1}\bigcup_{l=1}^{\ell-1}\hspace{-2pt}\bar B^l_{\ell-1,\ell}(1,s)\hspace{-2pt}\cup\hspace{-2pt}\bigcup_{\ell=1}^{\ell^*}\bigcup_{l=1}^{\ell-2}\notag\\
&\bar B^l_{\ell-1,\ell-1}(1,s)\hspace{-2pt}\cup\hspace{-2pt}\bar B_{\ell^*-1,\ell^*}^1(1,s)\hspace{-2pt}\cup\hspace{-3pt}\bigcup_{\ell_1=\ell^*}^{K-1}\hspace{-2pt}\bigcup_{\ell_2=\ell^*}^{K-1}\hspace{-3pt}\bar R_{\ell_1,\ell_2}(1,s)\hspace{-2pt}\cup\hspace{-2pt}\bigcup_{\ell_1=\ell^*}^{K-1}\hspace{-1pt}\bigcup_{\ell_2=\ell^*}^{K-1}\hspace{-1pt}\bigcup_{l=1}^{\ell_2-1}\hspace{-3pt}\bar B_{\ell_1,\ell_2}^l(1,s)\notag\\
&\cup\hspace{-2pt}\bigcup_{\ell_1=\ell^*}^{L-1}\bigcup_{\ell_2=\ell^*}^{L-1}\hspace{-2pt}\bar R_{\ell_1,\ell_2}(1,s)\hspace{-2pt}\cup\hspace{-2pt}\bigcup_{\ell_1=\ell^*}^{L-1}\bigcup_{\ell_2=\ell^*}^{L-1}\bigcup_{l=1}^{\ell_2-1}\hspace{-2pt}\bar B_{\ell_1,\ell_2}^l(1,s)\hspace{-2pt}\cup\hspace{-2pt}\mathscr S_\text{neg}^v(1,s)\cup\mathscr S_\text{neg}^h(1,s)\biggl].
\end{align}
Furthermore, there exists $k>0$ such that for $\beta$ sufficiently large
\begin{align}\label{ristimetubenzbneg}
\mathbbm P_\beta(\tau_{\partial^{np}\mathfrak T_{\mathcal X^s_\emph{neg}}(\mathbf 1)}^\mathbf 1\le\tau_{\mathcal X^s_\emph{neg}}^\mathbf 1)\le e^{-k\beta}.
\end{align}
\end{theorem}
Note that in \cite[Theorem 4.7]{bet2021metastabilitypos} the authors identify the tube of typical trajectories for the metastable transition for the $q$-Potts model with \textit{positive} external magnetic field. This tube has the a similar geometric definition of the tube \eqref{tesiteotubenegxs1} of our scenario, the main difference is that in this negative scenario we have to take the union on all $t\in S\backslash\{1\}$. 
\begin{remark}\label{remarktube1s}
In \cite{bet2021critical} the authors study the $q$-state Potts model with zero external magnetic field. Since in this energy landscape there are $q$ stable configurations and no relevant metastable states, the authors study the transitions between stable states. More precisely, they identify the union of all minimal gates and the tube of typical paths for the transition between two fixed stable states and these results hold also in the current scenario for the transition $\mathbf r\to\mathbf s$ for any $\mathbf r,\mathbf s\in\mathcal X^s_\text{neg}$, $\mathbf r\neq\mathbf s$. Indeed, the communication height computed in Subsection \ref{commheightstable} is equal to the one given in \cite{nardi2019tunneling} and its value is strictly lower than $\Phi_\text{neg}(\mathbf 1,\mathcal X^s_\text{neg})$ as we prove in Theorem \ref{theoremcomparisonneg}. It follows that  for any $\mathbf r,\mathbf s\in\mathcal X^s_\text{neg}$, $\mathbf r\neq\mathbf s$, any optimal path for the transition $\mathbf r\to\mathbf s$ does not visit the metastable state $\mathbf 1$ and for this type of transition the identification of the union of all minimal gates and of the tube of typical trajectories is given by \cite[Theorem 3.6]{bet2021critical} and \cite[Theorem 4.3]{bet2021critical}, respectively.
\end{remark}
Using Remark \ref{remarktube1s}, the tube of typical paths for the transition from the metastable to any fixed stable state is
\begin{align}\label{tesiteotubenegstabfix}
\mathfrak T_{\mathbf s}(\mathbf 1):=\mathfrak T_{\mathcal X^s_\text{neg}}(\mathbf 1)\cup\bigcup_{\mathbf r\in\mathcal X^s_\text{neg}\backslash\{\mathbf s\}} \mathfrak T_{\mathbf s}^{\text{zero}}(\mathbf r),
\end{align}
where $\mathfrak T_{\mathbf s}^{\text{zero}}(\mathbf r)$ is given by \cite[Equation 4.23]{bet2021critical}.
\begin{theorem}[Tube of typical paths for the transition $\mathbf 1\to\mathbf s$]\label{teotubefixneg}
If the external magnetic field is negative, then for any $\mathbf s\in\mathcal X^s_\emph{neg}$ the tube of typical trajectories for the transition $\mathbf 1\to\mathbf s$ is by \eqref{tesiteotubenegstabfix}. Furthermore, there exists $k>0$ such that for $\beta$ sufficiently large
\begin{align}\label{ristimetubenzbnegfix}
\mathbbm P_\beta(\tau_{\partial^{np}\mathfrak T_{\mathbf s}(\mathbf 1)}^\mathbf 1\le\tau_{\mathbf s}^\mathbf 1)\le e^{-k\beta}.
\end{align}
\end{theorem}
We defer the explicit proof of Theorem \ref{ristimetubenzbnegfix} in Subsection \ref{prooftube}.

\section{Energy landscape analysis}\label{secenergylandscape}
We devote this section to analyse the energy landscape of the $q$-state Potts model with negative external magnetic field.
\subsection{Disagreeing edges, bridges and crosses}\label{subseclocgeo}
In the following list we introduce the notions of \textit{disagreeing edges}, \textit{bridges} and \textit{crosses} of a Potts configuration on a grid-graph $\Lambda$. These definitions are taken from \cite{nardi2019tunneling}.

\noindent - We call $e=(v,w)\in E$ a \textit{disagreeing edge} if it connects two vertices with different spin values, i.e., $\sigma(v)\neq\sigma(w)$.

\noindent - For any $i=0,\dots,K-1$, let
\begin{align}\label{diseedgesrow}
d_{ r_i}(\sigma):=\sum_{(v,w)\in r_i} \mathbbm{1}_{\{\sigma(v)\neq\sigma(w)\}}
\end{align} 
be the total number of disagreeing horizontal edges on row $r_i$. Furthermore, for any $j=0,\dots,L-1$ let
\begin{align}\label{disedgescol}
d_{c_j}(\sigma):=\sum_{(v,w)\in c_j} \mathbbm{1}_{\{\sigma(v)\neq\sigma(w)\}},
\end{align} 
be the total number of disagreeing vertical edges on column $c_j$.

\noindent - We define $d_h(\sigma)$ as the total number of disagreeing horizontal edges and $d_v(\sigma)$ as the total number of disagreeing vertical edges.

\noindent Since we may partition the edge set $E$ in the two subsets of horizontal edges $E_h$ and of vertical edges $E_v$, such that $E_h\cap E_v=\varnothing$, the total number of disagreeing edges is given by
\begin{align}\label{totnumberdisedges}
\sum_{(v,w)\in E_v} \mathbbm{1}_{\{\sigma(v)\neq\sigma(w)\}}+\sum_{(v,w)\in E_h} \mathbbm{1}_{\{\sigma(v)\neq\sigma(w)\}}=d_v(\sigma)+d_h(\sigma).
\end{align} 

\noindent - We say that $\sigma$ has a \textit{horizontal bridge} on row $r$ if $d_r(\sigma)=0$.

\noindent - We say that $\sigma$ has a \textit{vertical bridge} on column $c$ if $d_c(\sigma)=0$.

\noindent - We say that $\sigma\in\mathcal X$ has a \textit{cross} if it has at least one vertical and one horizontal bridge. If $\sigma$ has a bridge of spins $s\in S$, then we say that $\sigma$ has an $s$-bridge. Similarly, if $\sigma$ has a cross of spins $s$, we say that $\sigma$ has an $s$-cross.

\begin{figure}
\centering
\begin{tikzpicture}[scale=0.6, transform shape]
\draw [fill=gray,gray] (0,0) rectangle (0.3,0.9) (0.3,0.6) rectangle (0.6,0.9) (2.7,0) rectangle (3,0.9) (2.4,2.1) rectangle (2.7,2.4) (0.6,0.9) rectangle (0.9,0.6) (0.9,1.8) rectangle (1.8,2.1) rectangle (1.5,1.5) (0.3,1.8)rectangle(0.6,2.1);
\draw [fill=gray,lightgray] (0.9,0) rectangle (1.8,0.3) (1.2,0.3) rectangle (1.5,0.6) (2.1,1.2) rectangle (2.7,1.8) rectangle (3,2.1) (2.4,1.8) rectangle (2.7,2.1) (0.3,1.2)rectangle(0.6,1.5);
\draw [fill=gray,black!85!white] (0.6,0) rectangle (0.9,2.4) (0.9,0.3) rectangle (1.2,0.6) (2.7,0.6) rectangle (3,1.8) (2.4,0.9) rectangle (3,1.2);
\draw [fill=gray,black!12!white] (0,0.9)rectangle(0.6,1.2) (0.9,0.9)rectangle(2.4,1.2)(1.8,1.5)rectangle(2.1,2.4)(1.5,2.1)rectangle(1.8,2.4);
\draw [fill=gray,black!65!white] (0.3,0)rectangle(0.6,0.6) (0,1.2)rectangle(0.3,1.5)(0,1.5)rectangle(0.6,1.8)(2.1,0)rectangle(2.4,0.9)(2.4,0.3)rectangle(2.7,0.6);
\draw[step=0.3cm,color=black] (0,0) grid (3.3,2.4);
\draw (0.75, 2.4) node[above] {$c_2$};\draw (3.15, 2.4) node[above] {$c_{10}$};
\draw (1.65,0) node[below] {\footnotesize{(a) Two vertical bridges}};\draw(1.65,-0.4) node[below]{\footnotesize{on columns $c_2$ and $c_{10}$.}};
\end{tikzpicture}\ \ \ \ \
\begin{tikzpicture}[scale=0.6, transform shape]
\draw [fill=gray,gray] (0,0) rectangle (0.3,2.1) (0.9,1.2)rectangle(1.5,1.5)(0.6,1.5)rectangle(0.9,2.1)(0.9,0) rectangle (1.8,0.3) (1.2,0.3) rectangle (1.5,0.6) (2.1,1.2) rectangle (2.7,1.8) rectangle (3,2.1) (2.4,1.8) rectangle (2.7,2.1) (0.3,1.2)rectangle(0.6,1.5);
\draw [fill=gray,lightgray] (0.3,0) rectangle (0.6,1.5)(2.1,1.8)rectangle(2.4,2.1)(2.1,0)rectangle(3.3,0.3)(0,1.2)rectangle(0.3,1.5)(0,1.5)rectangle(0.9,1.8)(0,0) rectangle (0.3,0.9) (0.3,0.6) rectangle (0.6,0.9) (3,0) rectangle (3.3,0.9) (2.4,2.1) rectangle (2.7,2.4) (0.6,0.9) rectangle (0.9,0.6) (0.9,1.8) rectangle (1.8,2.1) rectangle (1.5,1.5) (0.3,1.8)rectangle(0.6,2.1);;
\draw [fill=gray,black!12!white] (0,0.9)rectangle(0.6,1.2) (0.9,0.9)rectangle(2.4,1.2)(1.8,1.5)rectangle(2.1,2.4)(1.5,2.1)rectangle(1.8,2.4) (0.3,0)rectangle(0.6,0.6) (0,1.2)rectangle(0.3,1.5);
\draw [fill=gray,black!65!white] (0.3,0)rectangle(0.6,0.9) (0,1.2)rectangle(0.3,1.5)(0,1.5)rectangle(0.6,1.8)(2.1,0)rectangle(2.4,0.9)(2.4,0.3)rectangle(2.7,0.6)(1.2,1.8)rectangle(1.8,2.1);
\draw [fill=gray,black!85!white] (0,0.6)rectangle(3.3,0.9)(2.7,1.2)rectangle(3.3,1.5)rectangle(3,2.4) (0.9,0)rectangle(1.2,1.2);
\draw[step=0.3cm,color=black] (0,0) grid (3.3,2.4);
\draw (0,0.75) node[left] {$r_2$};
\draw (1.65,0) node[below] {\footnotesize{(b) An horizontal bridge}};\draw(1.65,-0.4) node[below]{\footnotesize{on row $r_2$.}};
\end{tikzpicture} \ \ \ \ \
\begin{tikzpicture}[scale=0.6, transform shape]
\draw [fill=gray,lightgray] (0.3,0)rectangle(0.6,0.6) (0,1.2)rectangle(0.3,1.5)(0,1.5)rectangle(0.6,1.8)(2.1,0)rectangle(2.4,0.9)(2.4,0.3)rectangle(2.7,0.6)(1.2,1.8)rectangle(1.8,2.1);
\draw [fill=gray,black!12!white] (0,0) rectangle (0.3,2.1) (0.9,1.2)rectangle(1.5,1.5)(0.6,1.5)rectangle(0.9,2.1)(0.9,0) rectangle (1.8,0.3) (1.2,0.3) rectangle (1.5,0.6) (2.1,1.2) rectangle (2.7,1.8) rectangle (3,2.1) (2.4,1.8) rectangle (2.7,2.1) (0.3,1.2)rectangle(0.6,1.5) (0,0.6)rectangle(3.3,0.9)(2.7,1.2)rectangle(3.3,1.5)rectangle(3,2.4) (0.9,0)rectangle(1.2,1.2);
\draw [fill=gray,gray] (0,0) rectangle (0.3,2.1) (0.9,1.2)rectangle(1.5,1.5)(0.6,1.5)rectangle(0.9,2.1)(0,0.9)rectangle(0.6,1.2) (0.9,0.9)rectangle(2.4,1.2)(1.8,1.5)rectangle(2.1,2.4)(1.5,2.1)rectangle(1.8,2.4);
\draw [fill=gray,black!85!white] (1.2,0) rectangle (1.5,2.4) (0.3,0.3) rectangle (0.6,0.6) (0,1.8)rectangle(3.3,2.1)(2.7,1.2)rectangle(3.3,1.5) (0.3,0)rectangle(0.6,0.6) (2.7,0)rectangle(3.3,0.3);
\draw[step=0.3cm,color=black] (0,0) grid (3.3,2.4);
\draw (0,1.95) node[left] {$r_6$};\draw (1.35, 2.4) node[above] {$c_4$};
\draw (1.65,0) node[below] {\footnotesize{(c) A cross on column $c_4$}};\draw(1.65,-0.4) node[below]{\footnotesize{ and on row $r_6$.}};
\end{tikzpicture}
\caption{\label{esbridgescross} Example of configurations on a $8\times 11$ grid graph displaying a vertical $s$-bridge (a), a horizontal $s$-bridge (b) and a $s$-cross (c). We color black the spins $s$.}
\end{figure}
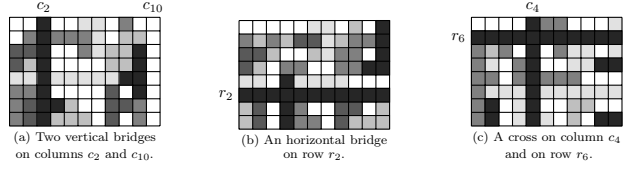\FloatBarrier

\noindent - For any $s\in S$, the total number of $s$-bridges of the configuration $\sigma$ is denoted by $B_s(\sigma)$. Note that if a configuration $\sigma\in\mathcal X$ has an $s$-cross, then $B_s(\sigma)$ is at least $2$.

\subsection{Metastable state and stability level of the metastable state}\label{submetastableneg}
In this subsection we find the unique metastable state and we compute its stability level. Furthermore, we find the set of the local minima and the set of the stable plateaux of the Hamiltonian \eqref{hamiltonianneg}. First we define a reference path from $\bold 1$ to $\bold s$, for any $\bold s\in\mathcal X^s_{\text{neg}}$. The energy value of any configuration $\sigma$ in this path, using \eqref{totnumberdisedges} can be computed as
\begin{align}\label{rewritegap1neg}
&H_{\text{neg}}(\sigma)-H_{\text{neg}}(\bold 1)
=d_v(\sigma)+d_h(\sigma)-h\sum_{u\in V} \mathbbm{1}_{\{\sigma(u)\neq1\}},
\end{align}
We say that a path $\omega\in\Omega_{\sigma,\sigma'}$ is the \textit{concatenation} of the $L$ paths $\omega^{(i)}=(\omega^{(i)}_0,\dots,\omega^{(i)}_{n_i}),$ for some $n_i\in\mathbb{N}$, $i=1,\dots,L$ if $\omega=(\omega^{(1)}_0=\sigma,\dots,\omega^{(1)}_{n_1},\omega^{(2)}_0,$ $\dots,\omega^{(2)}_{n_2},\dots,\omega^{(L)}_0,\dots,\omega^{(L)}_{n_L}=\sigma')$.
\begin{definition}[Reference path]\label{refpathmioneg} 
For any $\bold s\in\mathcal X^s_\text{neg}$, we define a reference path $\hat\omega:\bold 1\to\bold s$, $\hat\omega:=(\hat\omega_0,\dots,\hat\omega_{KL})$ as the concatenation of the two paths $\hat\omega^{(1)}:=(\bold 1=\hat\omega_0,\dots,\hat\omega_{(K-1)^2})$ and $\hat\omega^{(2)}:=(\hat\omega_{(K-1)^2+1},\dots,\bold s=\hat\omega_{KL})$. 
The path $\hat\omega^{(1)}$ is defined as follows. We set $\hat\omega_0:=\bold 1$. Then, we set $\hat\omega_1:=\hat\omega_0^{(i,j),s}$ where $(i,j)$ is any vertex in $\Lambda$. Sequentially, we flip clockwise from $1$ to $s$ all the vertices that sourround $(i,j)$ in order to construct a $3\times 3$ square. We iterate this construction until we get $\hat\omega_{(K-1)^2}\in\bar R_{K-1,K-1}(1,s)$. See Figure \ref{defrefpath}(a). The path $\hat\omega^{(2)}$ is defined as follows. Without loss of generality, assume that $\hat\omega_{(K-1)^2}\in\bar R_{K-1,K-1}(1,s)$ has the cluster of spin $s$ in the first $c_0,\dots,c_{K-2}$ columns, see Figure \ref{defrefpath}(b). We define $\hat\omega_{(K-1)^2+1},\dots,\hat\omega_{(K-1)^2+K-1}$ as a sequence of configurations in which the cluster of spins $s$ grows gradually by flipping the spins $1$ on the vertices $(K-1,j)$, for $j=0,\dots,K-2$. Thus, $\hat\omega_{(K-1)^2+K-1}\in\bar R_{K-1,K}(1,s)$ as depicted in Figure \ref{defrefpath}(c). Finally, we define the configurations $\hat\omega_{(K-1)^2+K},\dots,\hat\omega_{KL}$ as a sequence of states in which the cluster of spin $s$ grows gradually column by column. More precisely, starting from $\hat\omega_{(K-1)^2+K-1}\in\bar R_{K-1,K}(1,s)$, $\hat\omega^{(2)}$ passes through configurations in which the spins $1$ on columns $c_K,\dots,c_{L-1}$ become $s$. The procedure ends with $\hat\omega_{KL}=\bold s$. Note that the energy value of the configurations in the reference path is independent of the first flipped spin $(i,j)$.
\end{definition}
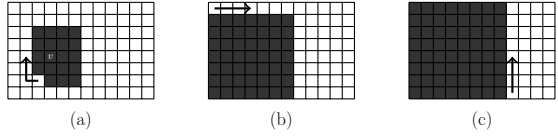
\begin{figure}[h!]
\centering 
\begin{tikzpicture}[scale=0.4, transform shape]
\draw [fill=gray,black!15!white] (1.6,1.2) rectangle (2,1.6);
\fill[black!80!white] (1.2,0.8)rectangle(1.6,2.4); 
\fill[black!80!white] (1.6,0.4) rectangle (2.8,2.4);
\draw (1.8,1.4) node[white,thick] {$v$};
\draw[<-,thick] (1,1.4) -- (1,0.585); \draw[thick] (1,0.6) -- (1.4,0.6);
\draw[step=0.4cm,color=black] (0.4,0) grid (5.2,3.2);
\draw (2.8,-0.3) node[below] {\LARGE{(a)}};
\end{tikzpicture}\ \ \ \ \ \ 
\begin{tikzpicture}[scale=0.4, transform shape]
\draw[->,thick] (0.6,3)--(1.8,3);
\fill[black!80!white] (0.4,0) rectangle (3.2,2.8);
\draw[step=0.4cm,color=black] (0.4,0) grid (5.2,3.2);
\draw (2.8,-0.3) node[below] {\LARGE{(b)}};
\end{tikzpicture}\ \ \ \ \ \ 
\begin{tikzpicture}[scale=0.4, transform shape]
\draw[->,thick] (3.8,0.2)--(3.8,1.4);
\fill[black!80!white] (0.4,0) rectangle (3.6,3.2);
\draw[step=0.4cm,color=black] (0.4,0) grid (5.2,3.2);
\draw (2.8,-0.3) node[below] {\LARGE{(c)}};
\end{tikzpicture}
\caption{\label{defrefpath} (a) First steps of path $\hat\omega^{(2)}$ on a $10\times 12$ grid $\Lambda$ starting from the vertex $v=(3,3)$. We color white the vertices with spin $1$, black those with spin $s$. The arrow indicates the order in which the spins are flipped from $1$ to $s$. (b) Illustration of $\hat\omega_{(K-1)^2}$. (c) Illustration of $\hat\omega_{(K-1)^2+K-1}$. }
\end{figure}\FloatBarrier

Next we show that any configuration belonging to $\bigcup_{t=2}^q\bar R_{\ell^*-1,\ell^*}(1,t)$ is connected to the metastable configuration $\bold 1$ by a path that does not overcome the energy value $4\ell^*-h(\ell^*(\ell^*-1)+1)+H_\text{neg}(\bold 1)$.
For any $s\in S$, we define as
\begin{align}\label{ennekappa}
N_s(\sigma):=|\{v\in V:\ \sigma(v)=s\}|
\end{align}
the number of vertices with spin $s$ in $\sigma\in\mathcal X$. 
\begin{lemma}\label{lemmaoneprefactor}
If the external magnetic field is negative, then for any $\sigma\in\bigcup_{t=2}^q\bar R_{\ell^*-1,\ell^*}(1,t)$ there exists a path $\gamma:\sigma\to\bold 1$ such that the maximum energy along $\gamma$ is bounded as
\begin{align}
\max_{\xi\in\gamma} H_\emph{neg}(\xi)<4\ell^*-h(\ell^*(\ell^*-1)+1)+H_\emph{neg}(\bold 1).
\end{align}
\end{lemma}
\textit{Proof.} Let $\sigma\in\bigcup_{t=2}^q\bar R_{\ell^*-1,\ell^*}(1,t)$. Hence, there exists $s\in\{2,\dots,q\}$ such that $\sigma\in\bar R_{\ell^*-1,\ell^*}(1,s)$. Consider the reference path of Definition \ref{refpathmioneg} and note that for any $i=0,\dots,KL$, $N_s(\hat\omega_i)=i$. The reference path may be constructed in such a way that $\hat\omega_{\ell^*(\ell^*-1)}:=\sigma$. Let $\gamma:=(\hat\omega_{\ell^*(\ell^*-1)}=\sigma,\hat\omega_{\ell^*(\ell^*-1)-1},\dots,\hat\omega_1,\hat\omega_0=\bold 1)$  be the time reversal of the subpath $(\hat\omega_0,\dots,\hat\omega_{\ell^*(\ell^*-1)})$ of $\hat\omega$. We claim that $\max_{\xi\in\gamma} H_\text{neg}(\xi)<4\ell^*-h(\ell^*(\ell^*-1)+1)+H_\text{neg}(\bold 1)$. Indeed, note that $\hat\omega_{\ell^*(\ell^*-1)}=\sigma,\dots,\hat\omega_1$ is a sequence of configurations in which all the spins are equal to $1$ except those, which are $s$, in either a quasi-square $\ell\times(\ell-1)$ or a square $(\ell-1)\times(\ell-1)$ possibly with one of the longest sides not completely filled. For any $\ell=\ell^*,\dots,2$, the path $\gamma$ moves from $\bar R_{\ell,\ell-1}(1,s)$ to $\bar R_{\ell-1,\ell-1}(1,s)$ by flipping the $\ell-1$ spins $s$ on one of the shortest sides of the $s$-cluster. In particular, $\hat\omega_{\ell(\ell-1)-1}$ is obtained by $\hat\omega_{\ell(\ell-1)}\in\bar R_{\ell,\ell-1}(1,s)$ by flipping the spin on a corner of the quasi-square from $s$ to $1$ and this increases the energy by $h$. The next $\ell-3$ steps are defined by flipping the spins on the incomplete shortest side from $s$ to $1$ where each step increases the energy by $h$. Finally, $\hat\omega_{(\ell-1)^2}\in\bar R_{\ell-1,\ell-1}(1,s)$ is defined by flipping the last spin $s$ to $1$ and this decreases the energy by $2-h$. For any $\ell=\ell^*,\dots,2$, $h(\ell-2)<2-h$. Indeed, $\ell^*=\left\lceil \frac{2}{h} \right\rceil$ and from Assumption \ref{remarkconditionneg}, we have $2-h>h(\ell^*-2)\ge h(\ell-2)$. Hence, $\max_{\xi\in\gamma} H_\text{neg}(\xi)=H_\text{neg}(\sigma)=4\ell^*-h(\ell^*(\ell^*-1)+1)-(2-h)+H_\text{neg}(\bold 1)$ and the claim is verified. $\qed$

In the next lemma we prove that any configuration in $\bigcup_{t=2}^q\bar B^2_{\ell^*-1,\ell^*}(1,t)$ is connected to the stable set $\mathcal X^s_\text{neg}$ by a path that does not overcome the energy value $4\ell^*-h(\ell^*(\ell^*-1)+1)+H_\text{neg}(\bold 1)$.
\begin{lemma}\label{lemmatwoprefactor}
If the external magnetic field is negative, then for any $\sigma\in\bar B^2_{\ell^*-1,\ell^*}(1,s)$, then there exists a path $\gamma:\sigma\to\bold s$ such that the maximum energy along $\gamma$ is bounded as
\begin{align}
\max_{\xi\in\gamma} H_\emph{neg}(\xi)<4\ell^*-h(\ell^*(\ell^*-1)+1)+H_\emph{neg}(\bold 1).
\end{align}  
\end{lemma}
\textit{Proof}. 
Consider the reference path of Definition \ref{refpathmioneg} and assume that this path is constructed in such a way that $\hat\omega_{\ell^*(\ell^*-1)+2}:=\sigma$. Let $\gamma:=(\hat\omega_{\ell^*(\ell^*-1)+2}=\sigma,\hat\omega_{\ell^*(\ell^*-1)+3},\dots,$ $\hat\omega_{KL-1},$ $\bold s)$. We claim that $\max_{\xi\in\gamma} H_\text{neg}(\xi)<4\ell^*-h(\ell^*(\ell^*-1)+1)+H_\text{neg}(\bold 1)$.
Since $\gamma$ is defined as a subpath of $\hat\omega$, we prove this claim by showing that $\max_{\xi\in\hat\omega}H_\text{neg}(\xi)=4\ell^*-h(\ell^*(\ell^*-1)+1)+H_\text{neg}(\bold 1)$ and that $\gamma$ does not intersect the unique configuration in which this maximum is reached. 
Indeed, for $\ell\le K-2$, note that the path $\hat\omega^{(1)}$ is defined by a sequence of configurations in which all the spins are equal to $1$ except those, which are $s$, in either a square $\ell\times\ell$ or a quasi-square $\ell\times(\ell-1)$ possibly with one of the longest sides not completely filled. For some $\ell\le K-2$, if $\hat\omega_{\ell(\ell-1)}\in\bar R_{\ell-1,\ell}(1,s)$ and $\hat\omega_{\ell^2}\in\bar R_{\ell,\ell}(1,s)$, then 
$\max_{\sigma\in\{\hat\omega_{\ell(\ell-1)},\hat\omega_{\ell(\ell-1)+1},\dots,\hat\omega_{\ell^2}\}} H_{\text{neg}}(\sigma)=H_{\text{neg}}(\hat\omega_{\ell(\ell-1)+1})=4\ell-h\ell^2+h\ell-h+H_{\text{neg}}(\bold 1).$
Otherwise, if $\hat\omega_j:=\hat\omega_{\ell^2}\in\bar R_{\ell,\ell}(1,s)$ and $\hat\omega_{\ell(\ell+1)}\in\bar R_{\ell,\ell+1}(1,s)$, then 
$\max_{\sigma\in\{\hat\omega_{\ell^2},\hat\omega_{\ell^2+1},\dots,\hat\omega_{\ell(\ell+1)}\}} H_{\text{neg}}(\sigma)=H_{\text{neg}}(\hat\omega_{\ell^2+1})=4\ell-h\ell^2+2-h+H_{\text{neg}}(\bold 1).$
Let $k^*:=\ell^*(\ell^*-1)+1$. By recalling the condition $\frac{2}{h}\notin\mathbb N$ of Assumption \ref{remarkconditionneg}(ii) and by studying the maxima of $H_{\text{neg}}$ as a function of $\ell$, we have $\text{arg max}_{\hat\omega^{(1)}} H_{\text{neg}}=\{\hat\omega_{k^*}\}.$

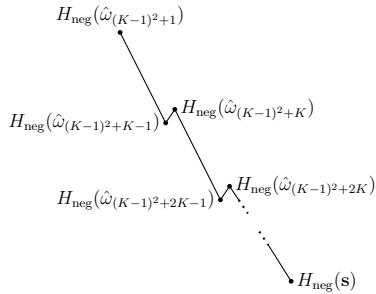
\begin{figure}[h!]
\centering
\begin{tikzpicture}[scale=0.6, transform shape]
\fill[color=black] (1,10) circle (1.5pt) node[above] {\large$H_{\text{neg}}(\hat\omega_{(K-1)^2+1})$}; 
\fill[color=black] (2,8) circle (1.5pt) node[left] {\large$H_{\text{neg}}(\hat\omega_{(K-1)^2+K-1})$}; 
\fill[color=black] (2.2,8.3) circle (1.5pt) node[above,right] {\large$H_{\text{neg}}(\hat\omega_{(K-1)^2+K})$}; 
\fill[color=black] (3.2,6.3) circle (1.5pt)  node[left] {\large$H_{\text{neg}}(\hat\omega_{(K-1)^2+2K-1})$}; 
\fill[color=black] (3.4,6.6) circle (1.5pt) node[above,right] {\large$H_{\text{neg}}(\hat\omega_{(K-1)^2+2K})$}; 
\draw[-] (1,10) -- (2,8) -- (2.2,8.3)--(3.2,6.3)--(3.4,6.6)--(3.6,6.3);
\draw[dotted,thick] (3.6,6.3)--(3.8,6) (4.05,5.6)--(4.25,5.3);
\draw [-] (4.25,5.3)--(4.75,4.5);
\fill[color=black] (4.75,4.5) circle (1.5pt) node[right] {\large$H_{\text{neg}}(\bold s)$};
\end{tikzpicture}
\caption{\label{Hmax} Qualitative illustration of the energy of the configurations belonging to $\hat\omega^{(2)}$.}
\end{figure}

\noindent Let us now study the maximum energy value reached along $\hat\omega^{(2)}$. This path is constructed by a sequence of configurations whose clusters of spins $s$ wrap around $\Lambda$. Moreover the maximum of the energy is reached by the first configuration, see Figure \ref{Hmax} for a qualitative representation of the energy of the configurations in $\hat\omega^{(2)}$. Indeed, using \eqref{energydifferenceneg}, we have
\begin{align}
&H_{\text{neg}}(\hat\omega_{(K-1)^2+j})-H_{\text{neg}}(\hat\omega_{(K-1)^2+j-1})=-2-h,\ j=2,\dots,K-1,\notag\\
&H_{\text{neg}}(\hat\omega_{(K-1)^2+K})-H_{\text{neg}}(\hat\omega_{(K-1)^2+K-1})=2-h,\notag\\
&H_{\text{neg}}(\hat\omega_{(K-1)^2+j})-H_{\text{neg}}(\hat\omega_{(K-1)^2+j-1})=-h,\ j=K+1,\dots,2K-1,\notag\\
&H_{\text{neg}}(\hat\omega_{(K-1)^2+2K})-H_{\text{neg}}(\hat\omega_{(K-1)^2+2K-1})=2-h.\notag
\end{align}
Using $K\ge 3\ell^*>3$, note that 
$H_{\text{neg}}(\hat\omega_{(K-1)^2+1})-H_{\text{neg}}(\hat\omega_{(K-1)^2+K})
=2K-6+h(K-1)>0.$
Moreover, 
$H_{\text{neg}}(\hat\omega_{(K-1)^2+K})-H_{\text{neg}}(\hat\omega_{(K-1)^2+2K})
=2K+2-h((K-1)^2+K)-(2K+2-h((K-1)^2+2K))=hK>0.$
By iterating the analysis of the energy gap between two consecutive configurations along $\hat\omega^{(2)}$, we conclude that $\text{arg max}_{\hat\omega^{(2)}} H_{\text{neg}}=\{\hat\omega_{(K-1)^2+1}\}$. In particular, 
\begin{align}\label{disequalityHneg}
H_{\text{neg}}(\hat\omega_{(K-1)^2+1})<H_{\text{neg}}(\hat\omega_{k^*})=4\ell^*-h(\ell^*(\ell^*-1)+1)+H_
\text{neg}(\bold 1)
\end{align} and, we refer to to Appendix \ref{appendixpath} for the explicit calculation. Hence, $\text{arg max}_{\hat\omega} H_{\text{neg}}=\{\hat\omega_{k^*}\}.$ Since $\gamma$ is defined as the subpath of $\hat\omega$ which goes from $\hat\omega_{\ell^*(\ell^*-1)+2}=\sigma$ to $\bold s$, $\gamma$ does not visit the configuration $\hat\omega_{k^*}$. Hence, 
$\max_{\xi\in\gamma} H_\text{neg}(\xi)<4\ell^*-h(\ell^*(\ell^*-1)+1)+H_\text{neg}(\bold 1)$ and the claim is proved.
$\qed$

\FloatBarrier

Next, we give an upper and a lower bound for $\Phi_{\text{{neg}}}(\bold 1,\mathcal X^s_{\text{{neg}}})-H_{\text{{neg}}}(\bold 1)$.
\begin{proposition}[Upper bound for the communication height]\label{refpathneg}
If the external magnetic field is negative, then  
\begin{align}\label{upperneg}
\Phi_{\text{\emph{neg}}}(\bold 1,\mathcal X^s_{\text{\emph{neg}}})-H_{\text{\emph{neg}}}(\bold 1)\le 4\ell^*-h(\ell^*(\ell^*-1)+1).
\end{align}
\end{proposition}
\textit{Proof.} The upper bound \eqref{upperneg} follows by the proof of Lemma \ref{lemmatwoprefactor}, where we proved that $\max_{\xi\in\hat\omega} H_\text{neg}(\xi)=H_\text{neg}(\hat\omega_{k^*})=4\ell^*-h(\ell^*(\ell^*-1)+1)+H_\text{neg}(\bold 1)$ where $\hat\omega$ is the reference path of Definition \ref{refpathmioneg}.
$\qed$

\begin{proposition}[Lower bound for the communication height]\label{lowerboundneg}
If the external magnetic field is negative, then  
\begin{align}\label{lowerbound2ms}
\Phi_{\text{\emph{neg}}}(\bold 1,\mathcal X^s_{\text{\emph{neg}}})-H_{\text{\emph{neg}}}(\bold 1)\ge 4\ell^*-h(\ell^*(\ell^*-1)+1).
\end{align}
\end{proposition}
\textit{Proof.} For any $\sigma\in\mathcal X$, we set $N(\sigma):=\sum_{t=2}^q N_t(\sigma)$, where $N_t(\sigma)$ is defined in \eqref{ennekappa}. Moreover, for all $k=1,\dots,|V|$, we define $\mathcal V_k:=\{\sigma\in\mathcal X: N(\sigma)=k\}$. Note that every path $\omega\in\Omega_{\bold 1,\mathcal X^s_{\text{neg}}}$ has to cross $\mathcal V_k$ for every $k=0,\dots,|V|$. 
In particular it has to intersect the set $\mathcal V_{k^*}$ with $k^*:=\ell^*(\ell^*-1)+1$. We prove the lower bound given in \eqref{lowerbound2ms} by computing that $H_\text{neg}(\mathscr F(\mathcal V_{k^*}))=4\ell^*-h(\ell^*(\ell^*-1)+1)+H_\text{neg}(\bold 1)$. Note that beacuse of the definition of $H_\text{neg}$ and of \eqref{rewritegap1neg}, the presence of disagreeing edges increases the energy. Thus, in order to describe the bottom $\mathscr F(\mathcal V_{k^*})$ we have to consider those configurations in which the $k^*$ spins different from $1$ belong to a unique $s$-cluster for some $s\neq 1$ inside a sea of spins $1$. Hence, consider the reference path $\hat\omega$ of Definition \ref{refpathmioneg} whose configurations satisfy this characterization. Note that $\hat\omega\cap\mathcal  V_{k^*}=\{\hat\omega_{k^*}\}$ with $\hat\omega\in\bar B^1_{\ell^*-1,\ell^*}(1,s)$. In particular, 
\begin{align}\label{engapconfrefpath}
H_{\text{neg}}(\hat\omega_{k^*})-H_{\text{neg}}(\bold 1)=4\ell^*-h(\ell^*(\ell^*-1)+1),
\end{align}
where $4\ell^*$ represents the perimeter of the cluster of spins different from $1$. Our goal is to prove that it is not possible to find a configuration with $k^*$ spins different from $1$ in a cluster of perimeter smaller than $4\ell^*$. Since the perimeter is an even integer, we assume that there exists a configuration belonging in $\mathcal V_{k^*}$ such that for some $s\in S\backslash\{ 1\}$ the $s$-cluster has perimeter $4\ell^*-2$. Since $4\ell^*-2<4\sqrt{k^*}$, where $\sqrt{k^*}$ is the side-length of the square $\sqrt{k^*}\times\sqrt{k^*}$ of minimal perimeter among those of area $k^*$ in $\mathbb R^2$, and since the square is the figure that minimizes the perimeter for a given area, we conclude that there does not exist a configuration with $k^*$ spins different from $1$ in a cluster with perimeter strictly smaller than $4\ell^*$. Hence, $\hat\omega_{k^*}\in\mathscr F(\mathcal V_{k^*})$ and \eqref{lowerbound2ms} is satisfied thanks to \eqref{engapconfrefpath}.
$\qed$

\begin{lemma}\label{lemmathreeprefactor}
If the external magnetic field is negative, then any $\omega\in\Omega_{\bold 1,\mathcal X^s_\emph{neg}}^{opt}$ is such that $\omega\cap\bigcup_{t=2}^q\bar R_{\ell^*-1,\ell^*}(1,t)\neq\varnothing$.
\end{lemma}
\textit{Proof}. At the beginning of the proof of Proposition \ref{lowerboundneg} we note that any path $\omega:\bold 1\to\mathcal X^s_\text{neg}$  has to visit $\mathcal V_k$ at least once for every $k=0,\dots,|V|$. Consider $\mathcal V_{\ell^*(\ell^*-1)}$. From \cite[Theorem 2.6]{alonso1996three} we get the unique configuration of minimal energy in $\mathcal V_{\ell^*(\ell^*-1)}$ is the one in which all spins are $1$ except those that are $s$, for some $s\in\{2,\dots,q\}$, in a quasi-square $\ell^*\times(\ell^*-1)$. In particular, this configuration has energy $\Phi_{\text{neg}}(\bold 1,\mathcal X^s_{\text{neg}})-(2-h)=4\ell^*-2-h\ell^*(\ell^*-1)+H_\text{neg}(\bold 1)$. Note that $4\ell^*-2$ is the perimeter of its $s$-cluster, $s\neq 1$. Since the perimeter is an even integer, we have that the other configurations belonging to $\mathcal V_{\ell^*(\ell^*-1)}$ have energy that is larger than or equal to $4\ell^*-h\ell^*(\ell^*-1)+H_\text{neg}(\bold 1)$. Thus, they are not visited by any optimal path. Indeed, $4\ell^*-h\ell^*(\ell^*-1)+H_\text{neg}(\bold 1)>\Phi_{\text{neg}}(\bold 1,\mathcal X^s_{\text{neg}}).$ Hence, we conclude that every optimal path intersects $\mathcal V_{\ell^*(\ell^*-1)}$ in a configuration belonging to $\bigcup_{t=2}^q\bar R_{\ell^*-1,\ell^*}(1,t)$.
$\qed$

Let $\sigma\in\mathcal X$ and let $v\in V$. We define the \textit{tile} centered in $v$, denoted by $v$-tile, as the set of five sites consisting of $v$ and its four nearest neighbors. See for instance Figure \ref{figmattonelleneg}.
A $v$-tile is said to be \textit{stable} for $\sigma$ if by flipping the spin on vertex $v$ from $\sigma(v)$ to any $s\in S$ the energy difference $H_\text{neg}(\sigma^{v,s})-H_\text{neg}(\sigma)$ is greater than or equal to zero. 

\noindent Our next aim is to prove a recurrence property in Proposition \ref{ricorrenzaneg}, which will be useful to prove that $\bold 1\in\mathcal X^m_{\text{neg}}$ as stated in Theorem \ref{teometastableneg}. In order to do this, in Lemma \ref{lemmastablevertices} for any configuration $\sigma\in\mathcal X$ we describe all the possible stable $v$-tiles induced by the Hamiltonian \eqref{hamiltonianneg} and we exploit this result to prove Proposition \ref{proplocalminima}.
For any $\sigma\in\mathcal X$, $v\in V$ and $s\in S$, we define $n_s(v)$ as the number of nearest neighbors to $v$ with spin $s$ in $\sigma$, i.e.,
\begin{align}\label{numberspinneighborhood}
n_s(v):=|\{w\in V:\ w\sim v,\ \sigma(w)=s\}|.
\end{align}
\begin{lemma}[Characterization of stable $v$-tiles for a configuration $\sigma$]\label{lemmastablevertices}
Let $\sigma\in\mathcal X$ and let $v\in V$. If the external magnetic field is negative, then the tile centered in $v$ is stable for $\sigma$ if and only if it satisfies one of the following conditions. 
\begin{itemize}
\item[\emph{(1)}] If $\sigma(v)=s\neq 1$, $v$ has at least two nearest neighbors with spin $s$, see Figure \ref{figmattonelleneg}\emph{(a),(c),(d),(f)--(i),(m)--(o), (q)}, or one nearest neighbor $s$ and three nearest neighbors with spin $r,t,z\in S\backslash\{s\}$, different from each other, see Figure \ref{figmattonelleneg}\emph{(r)--(s)}.
\item[\emph{(2)}] If $\sigma(v)=1$, $v$ has either at least three nearest neighbors with spin $1$ or two nearest neighbors with spin $1$ and two nearest neighbors with spin $r,s\neq 1$, $r\neq s$, see Figure \ref{figmattonelleneg}\emph{(b),(e),(l), (p)}.
\end{itemize}
\begin{figure}[h!]
\centering
\begin{tikzpicture}[scale=0.7,transform shape]
\foreach \i in {0,1.6,3.2,4.8,6.4,8,9.6,11.2,12.8}
\draw (\i,0)rectangle(\i+1.2,0.4);
\foreach \i in {0,1.6,3.2,4.8,6.4,8,9.6,11.2,12.8}
\draw (\i+0.4,-0.4)rectangle(\i+0.8,0.8);

\draw (0.6,-0.5) node[below] {(a)}(2.2,-0.5) node[below] {(b)}(3.8,-0.5) node[below] {(c)}(5.4,-0.5) node[below] {(d)}(7,-0.5) node[below] {(e)}(8.6,-0.5) node[below] {(f)}(10.2,-0.5) node[below] {(g)}(11.8,-0.5) node[below] {(h)}(13.4,-0.5) node[below] {(i)};

\foreach \i in {0.6,3.8,5.4,8.6,10.2,11.8,13.4}\draw (\i,0.2) node {$s$};
\foreach \i in {0.6,3.8,5.4}\draw (\i,0.6) node {$s$};
\foreach \i in {0.2,1,3.4,4.2,5,5.8,8.2,9,9.8,10.6,11.4,12.2,13,13.8}\draw  (\i,0.2) node {$s$}; 
\foreach \i in {1.8,2.2,2.6,6.6,7,7.4}\draw (\i,0.2) node {$1$};
\foreach \i in {2.2,7,10.2,13.4}\draw (\i,0.6) node {$1$};
\foreach \i in {2.2,5.4,10.2}\draw (\i,-0.2) node {$1$};
\draw (0.6,-0.2) node {$s$}; 
\foreach \i in {3.8,7,8.6,11.8,13.4}\draw (\i,-0.2) node {$r$};
\foreach \i in {8.6}\draw (\i,0.6) node {$r$};
\draw (11.8,0.6) node {$t$};
\end{tikzpicture}
\begin{tikzpicture}[scale=0.7,transform shape]
\foreach \i in {0,1.6,3.2,4.8,6.4,8,9.6,11.2}
\draw (\i,0)rectangle(\i+1.2,0.4);
\foreach \i in {0,1.6,3.2,4.8,6.4,8,9.6,11.2}
\draw (\i+0.4,-0.4)rectangle(\i+0.8,0.8);

\draw (0.6,-0.5) node[below] {(l)}(2.2,-0.5) node[below] {(m)}(3.8,-0.5) node[below] {(n)}(5.4,-0.5) node[below] {(o)}(7,-0.5) node[below] {(p)}(8.6,-0.5) node[below] {(q)}(10.2,-0.5) node[below] {(r)}(11.8,-0.5) node[below] {(s)};

\foreach \i in {0.2,0.6,1,2.6,4.2,6.6,7,12.2}\draw (\i,0.2) node {$1$};
\foreach \i in {7}\draw (\i,0.6) node {$1$};
\foreach \i in {2.2}\draw (\i,-0.2) node {$1$};
\foreach \i in {1.8,2.2,3.4,3.8,5,5.4,7.4,8.2,8.6,10.2,11.8}\draw (\i,0.2) node {$s$};
\foreach \i in {0.6,2.2,3.8,5.4,8.6,10.2,11.8}\draw (\i,0.6) node {$s$};
\foreach \i in {0.6,3.8,5.4,7,8.6,10.2,11.8}\draw (\i,-0.2) node {$r$};
\draw (9,0.2) node {$r$};
\foreach \i in {5.8,9.8,11.4}\draw (\i,0.2) node {$t$};
\draw (10.6,0.2) node {$z$};

\end{tikzpicture}
\caption{\label{figmattonelleneg} Stable tiles centered in any $v\in V$ for a $q$-Potts configuration $\sigma$ on $\Lambda$ for any $r,s,t,z\in S\backslash\{1\}$ different from each other. The tiles are depicted up to a rotation of $\alpha\frac\pi 2$, $\alpha\in\mathbb Z$. }
\end{figure}
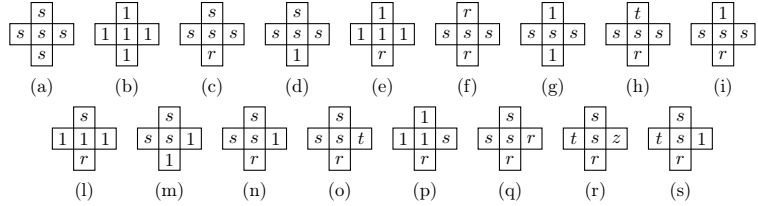\FloatBarrier
In particular, if $\sigma(v)=s$, then
\begin{align}\label{alignsummarizeenergygap}
H_{\emph{neg}}(\sigma^{v,r})-H_{\emph{neg}}(\sigma)=n_s(v)-n_r(v)-h\mathbbm{1}_{\{s=1\}}+h\mathbbm{1}_{\{r=1\}}.
\end{align}
\end{lemma}
\textit{Proof.} Let $\sigma\in\mathcal X$ and let $v\in V$. Assume that $\sigma(v)=s$, for some $s\in S$. To find if a $v$-tile is stable for $\sigma$ we reduce ourselves to flip the spin on vertex $v$ from $s$ to a spin $r$ such that
$v$ has at least one nearest neighbor $r$, i.e., $n_r(v)>1$.
Indeed, otherwise the energy difference \eqref{energydifferenceneg} is for sure strictly positive.
Let us divide the proof in several cases.\\
\noindent\textbf{Case 1.} Assume that $n_s(v)=0$ in $\sigma$. Then the corresponding $v$-tile is not stable for $\sigma$. Indeed, in view of the energy difference \eqref{energydifferenceneg}, if $r\neq 1$, by flipping the spin on vertex $v$ from $s$ to $r$ we have
\begin{align}\label{alignunomattonelleneg}
H_{\text{neg}}(\sigma^{v,r})-H_{\text{neg}}(\sigma)=-n_r(v)-h\mathbbm{1}_{\{s=1\}}.
\end{align}
Furthermore, for any $s\neq 1$, by flipping the spin on vertex $v$ from $s$ to $1$ we have 
\begin{align}
H_{\text{neg}}(\sigma^{v,1})-H_{\text{neg}}(\sigma)=-n_1(v)+h.
\end{align}
Hence, for any $s\in S$, if $v$ has spin $s$ and it has four nearest neighbors with spins different from $s$, i.e., $n_s(v)=0$, then the tile centered in $v$ is not stable for $\sigma$.\\
\noindent\textbf{Case 2.} Assume that $v\in V$ has three nearest neighbors with spin value different from $s$ in $\sigma$, i.e., $n_s(v)=1$. Then, in view of the energy difference \eqref{energydifferenceneg}, for any $s\in S$ and $r\notin\{1,s\}$, by flipping the spin on vertex $v$ from $s$ to $r$ we have
\begin{align}
H_{\text{neg}}(\sigma^{v,r})-H_{\text{neg}}(\sigma)=1-n_r(v)-h\mathbbm{1}_{\{s=1\}}.
\end{align}
Moreover, by flipping  the spin on vertex $v$ from $s\neq 1$ to $1$ we have
\begin{align}
&H_{\text{neg}}(\sigma^{v,1})-H_{\text{neg}}(\sigma)=1-n_1(v)+h.
\end{align}
Hence, for any $s\in S$, if $v$ has only one nearest neighbor with spin $s$, a tile centered in $v$ is stable for $\sigma$ only if $s\neq 1$ and $v$ has nearest neighbors with spins different from each other, see Figure \ref{figmattonelleneg}(r) and (s).\\
\noindent\textbf{Case 3.} Assume that $v\in V$ has two nearest neighbors with spin $s$, i.e., $n_s(v)=2$. Then, in view of the energy difference \eqref{energydifferenceneg}, for any $s\in S$ and $r\notin\{1,s\}$, by flipping the spin on vertex $v$ from $s$ to $r$ we have
\begin{align}
H_{\text{neg}}(\sigma^{v,r})-H_{\text{neg}}(\sigma)=2-n_r(v)-h\mathbbm{1}_{\{s=1\}}.
\end{align}
Moreover, by flipping the spin on vertex $v$ from $s\neq 1$ to $1$ we get
\begin{align}\label{minendiff}
H_{\text{neg}}(\sigma^{v,1})-H_{\text{neg}}(\sigma)=2-n_1(v)+h.
\end{align}
Hence, for any $s\in S$, if $v$ has two nearest neighbors with spin $s$ in $\sigma$, a $v$-tile is stable for $\sigma$ if $v$ has the other two nearest neighbors with different spin, see Figure \ref{figmattonelleneg}(m)--(q). Furthermore, if $s\neq 1$, a $v$-tile is stable for $\sigma$ even if $v$ has two nearest neighbors with spin $s$ and the other two nearest neighbors with the same spin, see Figure \ref{figmattonelleneg}(f)--(i).\\
\noindent\textbf{Case 4.} Assume that $v\in V$ has three nearest neighbors with spin $s$ in $\sigma$, i.e., $n_s(v)=3$, and that the fourth nearest neighbor has spin $r\neq s$. Then, for any $s\in S$ and $r\notin\{1,s\}$, we have
\begin{align}
&H_{\text{neg}}(\sigma^{v,r})-H_{\text{neg}}(\sigma)=2-h\mathbbm{1}_{\{s=1\}}.
\end{align}
Furthermore, by flipping the spin on vertex $v$ from $s\neq 1$ to $1$ we get
\begin{align}
H_{\text{neg}}(\sigma^{v,1})-H_{\text{neg}}(\sigma)=2+h.
\end{align}
\textbf{Case 5.} Assume that $v\in V$ has four nearest neighbors with spin $s$ in $\sigma$, i.e., $n_s(v)=4$. Then, for any $s\in S$ and $r\notin\{1,s\}$, we have
\begin{align}
H_{\text{neg}}(\sigma^{v,r})-H_{\text{neg}}(\sigma)=4-h\mathbbm{1}_{\{s=1\}}.
\end{align}
Furthermore, by flipping the spin on vertex $v$ from $s\neq 1$ to $1$ we get
\begin{align}\label{alignultimomattonelleneg}
H_{\text{neg}}(\sigma^{v,1})-H_{\text{neg}}(\sigma)=4+h.
\end{align}
From Case 4 and Case 5, for any $s\in S$, we get that a $v$-tile is stable for $\sigma$ if $v$ has at least three nearest neighbors with spin $s$, see Figure \ref{figmattonelleneg}(a)--(e). Finally, note that \eqref{alignsummarizeenergygap} is satisfied in all the cases $1$--$5$ above thanks to \eqref{alignunomattonelleneg}--\eqref{alignultimomattonelleneg}.
$\qed$

We define the set $C^s(\sigma)\subseteq\mathbb R^2$ as the union of unit closed squares centered at the vertices $v\in V$ such that $\sigma(v)=s$. We define $s$-\textit{clusters} the maximal connected components $C^s_1,\dots,C^s_n,\ n\in\mathbb N$, of $C^s(\sigma)$. 

\noindent For any $s\in S$, we say that a configuration $\sigma\in\mathcal X$ has an $s$-rectangle if it has a rectangular cluster in which all the vertices have spin $s$. 

\noindent Let $R_1$ an $r$-rectangle and $R_2$ an $s$-rectangle. They are said to be \textit{interacting} if either they intersect (when $r=s$) or are disjoint but there exists a site $v\notin R_1\cup R_2$ such that $\sigma(v)\neq r,s$ and $v$ has two nearest-neighbor $w, u$ lying inside $R_1, R_2$ respectively. For instance, in Figure \ref{figurerecurrenceproplabel}(b) the gray rectangles are not interacting. Furthermore, we say that $R_1$ and $R_2$ are \textit{adjacent} when they are at lattice distance one from each other, see for instance Figure \ref{figurerecurrenceproplabel}(c) and (e). 
We are now able to describe precisely the set of the local minima $\mathscr M_{\text{neg}}$ and the set of the stable plateaux $\bar{\mathscr M}_{\text{neg}}$ of the energy function \eqref{hamiltonianneg}. More precisely, the set of local minima $\mathscr M_\text{neg}$ is the set of stable points, i.e., $\mathscr M_\text{neg}:=\{\sigma\in\mathcal X:\ H_\text{neg}(\mathscr F(\partial\{\sigma\}))>H_\text{neg}(\sigma)\}$. While, a  plateau $D\subset\mathcal X$, namely a maximal connected set of equal energy states, is said to be \textit{stable} if $H_\text{neg}(\mathscr F(\partial D))>H_\text{neg}(D)$. Note that $\mathscr M_\text{neg}\cup\bar{\mathscr M}_\text{neg}\subset\hat{\mathcal X}_\text{neg}:=\{\sigma\in\mathcal X:$ for any $v\in V$ the tile centered in $v$ is stable$\}\subset\mathcal X$. In Proposition \ref{proplocalminima}, we prove that $\mathscr M_\text{neg}\cup\bar{\mathscr M}_\text{neg}$ is given by the union of the following sets. See also Figure \ref{figurerecurrenceproplabel}.

\noindent$\mathscr M^1_{\text{neg}}:=\{\bold 1,\bold2,\dots,\bold q\};$

\noindent$\mathscr M^2_{\text{neg}}:=\{\sigma\in\hat{\mathcal X}_\text{neg}:\ 
\sigma$ has strips of any spin $s\in S$ of thickness larger than or equal to one such that for any $s$ an $s$-strip of thickness one is in between strips of spins different from each other$\}$;

\noindent$\mathscr M^3_{\text{neg}}:=\{\sigma\in\hat{\mathcal X}_\text{neg}:\ \sigma$ has one or more $s$-rectangles for some $s\neq 1$, with minimum side-length larger than or equal to two, either in a sea of spins $1$ or inside a $1$-strip such that rectangles with the same spins are not interacting$\}$;

\noindent$\mathscr M^4_{\text{neg}}:=\{\sigma\in\hat{\mathcal X}_\text{neg}:\ \sigma$ has one or more $s$-rectangles for some $r,s\neq 1$, with minimum side-length larger than or equal to two, inside a $1$-strip adjacent to an $r$-strip$\}\cup\{\sigma\in\hat{\mathcal X}_\text{neg}:\ \sigma$ is covered by interacting $s$-rectangles such that each spin on the corners has outside the rectangle two nearest neighbors with different spins from each other and from the one inside the rectangle$\}$

\noindent$\bar{\mathscr M}^1_\text{neg}:=\{\sigma\in\hat{\mathcal X}_\text{neg}$: for any $r,s\neq 1$, $\sigma$ has an $s$-cluster with two consecutive sides next either to a connected $r$-cluster or to two $r$-cluster and the sides on the interfaces are of different length$\}$.

\begin{remark}
The set $\bar{\mathscr M}^1_\text{neg}$ is defined by fixing a representative configuration $\sigma$ and implicitly it includes also all the configurations connected to $\sigma$ via a path along which the energy is constant, see Figure \ref{figzeronodiscesa}.
\end{remark}

A path $\omega=(\omega_0,\dots,\omega_n)$ is said to be \textit{downhill} (\textit{strictly downhill}) if $H(\omega_{i+1})\le H(\omega_i)$ ($H(\omega_{i+1})<H(\omega_i)$) for $i=0,\dots,n-1$.

\begin{proposition}[Sets of local minima and of stable plateaux]\label{proplocalminima}
If the external magnetic field is negative, then 
\begin{align}
\mathscr M_{\emph{neg}}\cup\bar{\mathscr M}_{\emph{neg}}=\mathscr M^1_{\emph{neg}}\cup\mathscr M^2_{\emph{neg}}\cup\mathscr M^3_{\emph{neg}}\cup\mathscr M^4_{\emph{neg}}\cup\bar{\mathscr M}^1_\emph{neg}.
\end{align}
\end{proposition}
\textit{Proof.} A configuration $\sigma\in\mathcal X$ is a local minimum when, for any $v\in V$ and $s\in S$, the energy difference \eqref{energydifferenceneg} is strictly positive. On the other hand, $\sigma$ belongs to a stable plateau when, for any $v\in V$ and $s\in S$, the energy difference \eqref{energydifferenceneg} is larger than or equal to zero. Since a local minimum and a stable plateau are the union of stable tiles, we obtain \textit{all} the local minima and \textit{all} the stable plateaux by considering all the possible ways in which the stable tiles may be combined. We do this in various steps. First we consider all configurations which can be obtained from combining tiles (a)--(b). Then, we progressively add more tile types and construct all the possible resulting configurations. To refer to a tile type, we will use its corresponding lett in Figure \ref{figmattonelleneg}. 

\textbf{Step 1}. If $\sigma$ has only stable tiles as in Figure \ref{figmattonelleneg}(a) and (b), then there are no interfaces and $\sigma\in\mathscr M^1_{\text{neg}}$. 

\textbf{Step 2}. Let us assume now that the only stable tiles in $\sigma$ are (a)--(l).
Note that if $\sigma$ contains a tile of type (f), then $\sigma$ does not belong to $\mathscr M_{\emph{neg}}\cup\bar{\mathscr M}_{\emph{neg}}$. Indeed, if $\sigma$ contains at least an (f) tile, then it also contains an $s$-strip of thickness one in between two $r$-strips and there exists a downhill path of two steps. First, flip from $s$ to $r$ the central spin $s$ and this does not change the energy, then flip from $s$ to $r$ a spin, which, has now three spin $r$ neighbor. This flip reduces the energy by $2$.
On the other hand, any spin update on the central vertex of the tiles (a)--(e) and (g)--(l) strictly increases the energy of $\sigma$. By considering these, we obtain that $\sigma$ may contain horizontal (resp. vertical)  interfaces of length $L$ (resp. $K$). In particular, for any $s\in S$, an $s$-strip of thickness one must be either in between strips with different spins, using (h)--(l) tiles, or in between two $1$-strips if $s\neq 1$, using (g) tiles. 
We conclude that if $\sigma$ is obtained by a combination of the stable tiles (a)--(l), then $\sigma\in\mathscr M^1_{\text{neg}}\cup\mathscr M^2_{\text{neg}}$, see Figure \ref{figurerecurrenceproplabel}(a).

\textbf{Step 3}. Next we consider those $\sigma$ that are defined as the combination of the stable tiles (a)--(e), (g)--(p). Any spin update on the central vertex of the tiles (m)--(p). 
Since the central spin $s\neq 1$ of these tiles has at least two nearest neighbors with the same spin, the admissible shapes of an $s$-cluster are either strips or rectangles. It follows that the local minima containing only tiles (a)--(p) may additionally contain the following shapes.

\noindent (i) One or more $s$-rectangles ($s\neq 1$) with minimum side length two either in a sea of spins $1$ or inside a $1$-strip under the condition that rectangles with the same spins are not interacting, see Figure \ref{figurerecurrenceproplabel}(b). 

\noindent (ii) One or more $s$-rectangles ($s\neq 1$) with minimum side length two, inside a $1$-strip, with a side adjacent to an $r$-strip ($r\neq 1$), see Figure \ref{figurerecurrenceproplabel}(d).

\noindent (iii) Alternatively, $\sigma$ is covered by interacting rectangles under the condition that each spin on the corners is the centre of a tile of type (n)--(o), see Figure \ref{figurerecurrenceproplabel}(e).

We conclude that if $\sigma$ is defined by the combination of the stable tiles as in Figure \ref{figmattonelleneg}(a)--(p), then $\sigma\in\mathscr M^1_{\text{neg}}\cup\mathscr M^2_{\text{neg}}\cup\mathscr M^3_{\text{neg}}\cup\mathscr M^4_{\text{neg}}$.
 
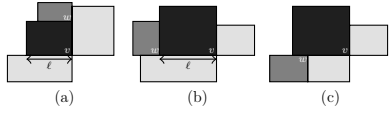
\begin{figure}
\centering
\begin{tikzpicture}[scale=0.5,transform shape]
\fill[black!12!white] (10.5,0)rectangle(12.2,0.7)(12.2,0.7)rectangle(13.3,2) (12.2,0.7)rectangle(11,1.6);
\fill[black!88!white] (12.2,0.7)rectangle(11,1.6);
\fill[black!50!white] (12.2,1.6) rectangle (11.3,2.1);
\draw (10.5,0)rectangle(12.2,0.7)(12.2,0.7)rectangle(13.3,2) (12.2,0.7)rectangle(11,1.6)(12.2,1.6) rectangle (11.3,2.1);
\draw[<->] (12.2,0.6)--(11,0.6); \draw(11.6,0.65) node[below] {\footnotesize{$\ell$}};
\draw  (12.08,1) node[white,below]{\footnotesize{$v$}};
\draw  (12.08,1.88) node[white,below]{\footnotesize{$w$}};
\draw(12,-0.1) node[below] {\large{(a)}};

\fill[black!12!white] (14,0)rectangle(16,0.7)(16,0.7)rectangle(17,1.5)(14.5,0.7)rectangle(16,2);
\fill[black!88!white] (14.5,0.7)rectangle(16,2);
\fill[black!50!white] (14.5,0.7) rectangle (13.8,1.6);
\draw (14,0)rectangle(16,0.7)(16,0.7)rectangle(17,1.5)(14.5,0.7)rectangle(16,2) (14.5,0.7) rectangle (13.8,1.6);
\draw[<->] (14.5,0.6)--(16,0.6); \draw(15.25,0.65) node[below] {\footnotesize{$\ell$}};
\draw  (15.88,1) node[white,below]{\footnotesize{$v$}};
\draw  (14.38,1) node[white,below]{\footnotesize{$w$}};
\draw(15.5,-0.1) node[below] {\large{(b)}};

\fill[black!12!white] (18.4,0)rectangle(19.5,0.7)(19.5,0.7)rectangle(20.5,1.5);
\fill[black!88!white] (18,0.7)rectangle(19.5,2);
\fill[black!50!white] (18.4,0.7) rectangle (17.4,0);
\draw (18.4,0)rectangle(19.5,0.7)(19.5,0.7)rectangle(20.5,1.5)(18,0.7)rectangle(19.5,2)(18.4,0.7) rectangle (17.4,0);
\draw  (19.38,1) node[white,below]{\footnotesize{$v$}};
\draw  (18.28,0.78) node[white,below]{\footnotesize{$w$}};
\draw(19,-0.1) node[below] {\large{(c)}};
\end{tikzpicture}
\caption{\label{figvertexv} Illustration of an $s$-rectangle, that we color black, adjacent to two $r$-rectangles, that we color light gray. Furthermore, we color gray those $t$-rectangles with $t\in S\backslash\{r,s\}$.}
\end{figure}

\begin{figure}
\centering
\begin{tikzpicture}[scale=0.45, transform shape]
\fill[black!12!white] (0,0) rectangle (3.6,2.7);
\fill[black!88!white] (0,2.7)rectangle(1.8,1.2);
\fill[black!0!white] (0,1.2)rectangle(0.6,0);
\fill[black!70!white] (0.6,0)rectangle(1.8,0.6)(3.6,2.7)rectangle(2.7,1.2);
\fill[black!30!white] (1.8,1.2)rectangle(3.6,0);
\fill[black!50!white] (1.8,2.7)rectangle(2.7,1.8);
\draw (1.65,1.35) node[white] {\footnotesize{$v$}};
\draw[step=0.3cm,color=black] (0,0) grid (3.6,2.7);
\draw(1.8,-0.1) node[below] {{$\omega_0:=\sigma$}};
\draw[->] (3.7,1.5)--(4.4,1.5); \draw [<-]  (3.7,1.2)--(4.4,1.2);
\fill[black!12!white] (4.5,0) rectangle (8.1,2.7);
\fill[black!88!white] (4.5,2.7)rectangle(6.3,1.2);
\fill[black!0!white] (4.5,1.2)rectangle(5.1,0);
\fill[black!70!white] (5.1,0)rectangle(6.3,0.6)(8.1,2.7)rectangle(7.2,1.2);
\fill[black!30!white] (6.3,1.2)rectangle(8.1,0);
\fill[black!50!white] (6.3,2.7)rectangle(7.2,1.8);
\fill[black!12!white] (6,1.2)rectangle(6.3,1.5);
\draw(6.3,-0.1) node[below] {{$\omega_1$}};
\draw (6.15,1.35) node {\footnotesize{$v$}};
\draw[step=0.3cm,color=black] (4.5,0) grid (8.1,2.7);
\draw[->] (8.2,1.5)--(8.9,1.5); \draw [<-]  (8.2,1.2)--(8.9,1.2);
\fill[black!12!white](9,0) rectangle (12.6,2.7);
\fill[black!88!white] (9,2.7)rectangle(10.8,1.2);
\fill[black!0!white] (9,1.2)rectangle(9.6,0);
\fill[black!70!white] (9.6,0)rectangle(10.8,0.6)(12.6,2.7)rectangle(11.7,1.2);
\fill[black!30!white] (10.8,1.2)rectangle(12.6,0);
\fill[black!50!white] (10.8,2.7)rectangle(11.7,1.8);
\fill[black!12!white] (10.8,1.2)rectangle(10.2,1.5);
\draw(10.8,-0.1) node[below] {{$\omega_2$}};
\draw (10.65,1.35) node {\footnotesize{$v$}};
\draw[step=0.3cm,color=black] (9,0) grid (12.6,2.7);
\end{tikzpicture}
\begin{tikzpicture}[scale=0.45, transform shape]
\draw[->] (12.7,1.5)--(12.9,1.5); \draw [<-]  (12.7,1.2)--(12.9,1.2);
\draw (13.2,1.35) node {\footnotesize$\dots$};
\draw[->] (13.7,1.5)--(13.5,1.5); \draw [<-]  (13.7,1.2)--(13.5,1.2);
\fill[black!12!white](13.8,0) rectangle (17.4,2.7);
\fill[black!88!white] (13.8,2.7)rectangle(15.6,1.2);
\fill[black!0!white] (13.8,1.2)rectangle(14.4,0);
\fill[black!70!white] (14.4,0)rectangle(15.6,0.6)(17.4,2.7)rectangle(16.5,1.2);
\fill[black!30!white] (15.6,1.2)rectangle(17.4,0);
\fill[black!50!white] (15.6,2.7)rectangle(16.5,1.8);
\fill[black!12!white] (15.6,1.2)rectangle(14.4,1.5) (15.6,1.5)rectangle(15.3,1.8);
\draw(15.6,-0.1) node[below] {{$\omega_j$}};
\draw (15.45,1.35) node {\footnotesize{$v$}};
\draw[step=0.3cm,color=black] (13.8,0) grid (17.4,2.7);
\draw[->] (17.5,1.5)--(17.7,1.5); \draw [<-]  (17.5,1.2)--(17.7,1.2);
\draw (18,1.35) node {\footnotesize$\dots$};
\draw[->] (18.5,1.5)--(18.3,1.5); \draw [<-]  (18.5,1.2)--(18.3,1.2);
\fill[black!12!white](18.6,0) rectangle (22.2,2.7);
\fill[black!88!white] (18.6,2.7)rectangle(20.4,1.2);
\fill[black!0!white] (18.6,1.2)rectangle(19.2,0);
\fill[black!70!white] (19.2,0)rectangle(20.4,0.6)(22.2,2.7)rectangle(21.3,1.2);
\fill[black!30!white] (20.4,1.2)rectangle(22.2,0);
\fill[black!50!white] (20.4,2.7)rectangle(21.3,1.8);
\fill[black!12!white] (20.4,1.2)rectangle(19.2,1.8);
\draw(20.4,-0.1) node[below] {{$\omega_n$}};
\draw (20.25,1.35) node {\footnotesize{$v$}};
\draw[step=0.3cm,color=black] (18.6,0) grid (22.2,2.7);
\end{tikzpicture}
\caption{\label{figzeronodiscesa} Example of a path $\omega:=(\omega_0,\dots,\omega_n)$ started in a configuration $\omega_0:=\sigma$ with a cluster as the one depicted in Figure \ref{figvertexv}(c) and such that $H_\text{neg}(\omega_{i})=H_\text{neg}(\omega_j)$, for any $i,j=0,\dots,n$. Since all the configurations depicted have the same energy value and they are connected by means a path, they belong to a stable plateau.}
\end{figure}
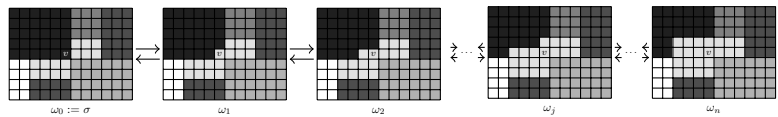

\textbf{Step 4}. Next we consider those $\sigma$ that are obtained as the combination of the stable tiles (a)--(e), (g)--(q). Combining the tiles of type (q) with all the previous ones, we obtain that for any $r,s\neq 1, r\neq s$, an $s$-cluster may have two consecutive sides adjacent either to a connected cluster or to two clusters with spins $r$ and the sides on the interface may have either the same or different length, see Figure \ref{figvertexv}. We claim that in a stable configuration $\sigma$ these are no clusters as in Figure \ref{figvertexv}(a)--(b), in which an $r$-cluster has a side longer than or equal to the side of the the $s$-cluster on the interface. Indeed, the path $(\omega_1=\sigma,\dots,\omega_\ell)$ that flips from $s$ to $r$ all the spins $s$ on the interface of length $\ell$ visits states $\omega_1,\dots,\omega_{\ell-1}$ with constant energy and $H_\text{neg}(\omega_\ell)<H_\text{neg}(\omega_{\ell-1})$.

\FloatBarrier
Let us now focus on the case in Figure \ref{figvertexv}(c), and let $\sigma$ be a configuration with such clusters. We prove that $\sigma\in\bar{\mathscr M}_\text{neg}$. In particular, all configurations connected to $\sigma$ via a path along which the energy does not change also belong to $\bar{\mathscr M}_{\text{neg}}$. In order to see this, consider the path which flips from $s$ to $r$ the spins $s$ adjacent to an $r$-rectangle (see for instance the path depicted in Figure \ref{figzeronodiscesa}). Note that at any step the energy does not change. Hence, combining all the stable tiles (a)--(q), we conclude that $\sigma\in\mathscr M^1_{\text{neg}}\cup\mathscr M^2_{\text{neg}}\cup\mathscr M^3_{\text{neg}}\cup\mathscr M^4_{\text{neg}}\cup\bar{\mathscr M}^1_\text{neg}$.

\textbf{Step 5}. Finally, assume that $\sigma$ may be obtained by a combination of tiles (a)--(s).
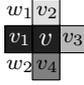
\begin{figure}
\centering
\begin{tikzpicture}[scale=0.9]
\fill[black!88!white] (-0.4,0.4)rectangle(0.4,0.8);
\fill[black!12!white] (0,0.8)rectangle(0.4,1.2);
\fill[black!50!white] (0,0)rectangle(0.4,0.4);
\fill[black!28!white] (0.4,0.4)rectangle(0.8,0.8);
\draw[step=0.4cm,color=black] (0,0) grid (0.4,1.2) (-0.4,0.4)grid(0.8,0.8);
\draw (0.2,0.6) node[white] {$v$} (-0.18,0.6) node[white] {\footnotesize{$v_1$}}(0.6,0.6) node {\footnotesize{$v_3$}} (0.22,1) node {\footnotesize{$v_2$}} (0.22,0.2) node {\footnotesize{$v_4$}} (-0.18,0.2) node{\footnotesize{$w_2$}} (-0.18,1) node{\footnotesize{$w_1$}};
\end{tikzpicture}
\caption{\label{figcolorvertx} Example of a $v$-tile equal to the one depicted in Figure \ref{figmattonelleneg}(r)--(s). We do not color the vertices $w_1$ and $w_2$ since in the proof they assume different value in different steps.}
\end{figure}

For this step, we refer to Figure \ref{figcolorvertx}, where we represent $r,s,t,z$ respectively by $\mycirclelightgray,\mycircleblack,\mycircleintgray,\mycirclegray$ and where we take $r,t\notin\{s,1\}$ and $z\neq s$. Let us assume that this type of tile belongs to a configuration $\sigma$ and consider the following cases.

 \textbf{Step 5.1}. If $n_s(v_1)=4$, then $\sigma(w_1)=\sigma(w_2)=s$. If both the $v_2$-tile and $v_4$-tile are of type (m), then $v_3$ would be the central vertex of a unstable tile. Thus, at least one of them is of type (q). Proceeding as in Step 4 we show that $\sigma$ is either unstable or it belongs to a stable plateau.

 \textbf{Step 5.2}. Assume $n_s(v_1)=3$. If $\sigma(w_1)=\sigma(w_2)=s$, then again $\sigma$ is either unstable or belongs to a stable plateau. If $\sigma(w_1)=s$ and $\sigma(w_2)\neq s$, then $v_1$ must be the central vertex of a tile of type (q) and again $\sigma$ is either unstable or belongs to a stable plateau. Otherwise, $v_3$ would be the central vertex of a unstable tile.

 \textbf{Step 5.3}. We now consider the case $n_s(v_1)=1$. This will be useful to study the case $n_s(v_1)=2$ in the next step. Along the path $\omega:=(\sigma,\sigma^{v,r},(\sigma^{v,r})^{v_1,r})$  the energy decreases. Indeed,
\begin{align}
&H_{\text{neg}}(\sigma^{v,r})-H_{\text{neg}}(\sigma)=0,\label{alignnegoneminloc} \\
&H_{\text{neg}}((\sigma^{v,r})^{v_1,r})\hspace{-2pt}-\hspace{-2pt}H_{\text{neg}}(\sigma^{v,r})\hspace{-2pt}=\hspace{-2pt}
\begin{cases} -2,\ \text{if}\ n_r(v_1)=1;\\ 
-3,\ \text{if}\ n_r(v_1)=2;\\ 
-4,\ \text{if}\ n_r(v_1)=3.
\end{cases}\label{alignnegtwominloc}
\end{align}
It follows that the tiles as (r)--(s) with $n_s(v_1)=1$ do not belong to any configuration in $\mathscr M_{\text{neg}}\cup\bar{\mathscr M}_{\text{neg}}$.

\textbf{Step 5.4}. Lastly, let us consider the case $n_s(v_1)=2$. Without loss of generality, assume that the spin $s$ nearest neighbors of $v_1$ lie on the same row. Consider the following two cases:

- $v_1$ has at least one nearest neighbor with a spin among $r,t,z\notin\{1,s\}$, say $r$, then along the path $(\sigma,\sigma^{v,r},(\sigma^{v,r})^{v_1,r})$ the energy decreases. Indeed, we have
\begin{align}
&H_{\text{neg}}(\sigma^{v,r})-H_{\text{neg}}(\sigma)=0\ \ \text{and}\ \
H_{\text{neg}}((\sigma^{v,r})^{v_1,r})-H_{\text{neg}}(\sigma^{v,r})\le-1.
\end{align}
Thus, there are no such tiles in configurations in $\mathscr M_{\text{neg}}\cup\bar{\mathscr M}_{\text{neg}}$.

- $v_1$ has two nearest neighbors with spin $s$ on vertices $v$ and $v_5$ and two nearest neighbors with spins $r'_1, r'_2\notin\{r,t,z\}$. If $r'_1=r'_2=r$, then the path $(\sigma,\sigma^{v,r},(\sigma^{v,r})^{v_1,r'})$ is downhill. Indeed, 
$H_{\text{neg}}(\sigma^{v,r})-H_{\text{neg}}(\sigma)=0$ and $H_{\text{neg}}((\sigma^{v,r})^{v_1,r'})-H_{\text{neg}}(\sigma^{v,r})=-1+h\mathbbm 1_{\{r'=1\}}$. 
\begin{figure}[h!]
\centering
\begin{tikzpicture}[scale=0.8]
\fill[black!75!white] (-0.4,0.8)rectangle(0,1.2);
\fill[black!88!white] (-0.8,0.4)rectangle(0.4,0.8);
\fill[black!12!white] (0,0.8)rectangle(0.4,1.2);
\fill[black!50!white] (0,0)rectangle(0.4,0.4);
\fill[black!28!white] (0.4,0.4)rectangle(0.8,0.8);
\draw[step=0.4cm,color=black] (-0.4,0) grid (0.4,1.2) (-0.8,0.4)grid(0.8,0.8);
\draw (0.2,0.6) node[white] {$v$} (-0.18,0.6) node[white] {\footnotesize{$v_1$}}(0.6,0.6) node {\footnotesize{$v_3$}} (0.22,1) node {\footnotesize{$v_2$}} (0.22,0.2) node {\footnotesize{$v_4$}} (-0.58,0.6) node[white] {\footnotesize{$v_5$}};
\draw (0.1,-0.2) node[below]{(a)};
\end{tikzpicture}\ \ \ \ \ \ \ \ \ \
\begin{tikzpicture}[scale=0.8]
\fill[black!75!white] (-0.4,0.8)rectangle(0,1.2);
\fill[black!88!white] (-1.2,0.4)rectangle(0,0.8);
\fill[black!12!white] (0,0.8)rectangle(0.4,1.2)(0,0.4)rectangle(0.4,0.8);
\fill[black!50!white] (0,0)rectangle(0.4,0.4);
\fill[black!28!white] (0.4,0.4)rectangle(0.8,0.8);
\draw[step=0.4cm,color=black] (-0.4,0) grid (0.4,1.2) (-1.2,0.4)grid(0.8,0.8);
\draw (0.2,0.6) node {$v$} (-0.18,0.6) node[white] {\footnotesize{$v_1$}}(0.6,0.6) node {\footnotesize{$v_3$}} (0.22,1) node {\footnotesize{$v_2$}} (0.22,0.2) node {\footnotesize{$v_4$}} (-0.58,0.6) node[white] {\footnotesize{$v_5$}} (-1,0.6) node[white] {\footnotesize{$v_6$}};
\draw (0.1,-0.2) node[below]{(b)};
\end{tikzpicture}\ \ \ \ \ \ \ \ \ \
\begin{tikzpicture}[scale=0.8]
\fill[black!88!white] (-1.2,0.4)rectangle(-0.4,0.8);
\fill[black!12!white] (0.8,0.4)rectangle(1.2,1.2);
\fill[black!90!white] (0,0.4)rectangle(0.8,0.8);
\fill[black!35!white](0.4,0)rectangle(0.8,0.4)(-1.2,0)rectangle(-0.8,0.4);
\fill[black!70!white] (-1.6,0.4)rectangle(-1.2,0.8);
\fill[black!60!white] (0.4,0.8)rectangle(0.8,1.2);
\fill[black!28!white] (1.2,0.4)rectangle(1.6,0.8);
\fill[black!75!white] (0.4,0.8)rectangle(0.8,1.2);
\draw[step=0.4cm,color=black] (0.4,0) grid (1.2,1.2) (0,0.4)grid(1.6,0.8);
\draw[step=0.4cm,color=black](-1.2,0)grid(-0.8,1.2)(-1.6,0.4)rectangle(-0.4,0.8);
\draw (1,0.6) node {$v$} (0.6,0.6) node[white] {\footnotesize{$v_1$}}(1.4,0.6) node {\footnotesize{$v_3$}} (1.02,1) node {\footnotesize{$v_2$}} (1.02,0.2) node {\footnotesize{$v_4$}} (0.2,0.6) node[white] {\footnotesize{$v_5$}} (-1,0.6)node[white] {\footnotesize{$u$}}(-0.6,0.6) node[white] {\footnotesize{$v_n$}};
\draw (-0.17,1)  node {$\dots$} (-0.17,0.2) node {$\dots$};
\draw (-0.2,-0.2) node[below]{(c)};
\end{tikzpicture}
\caption{\label{step5finale} Illustration of the Step 5.4.}
\end{figure}
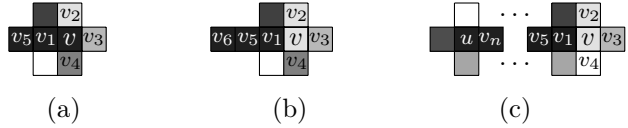\FloatBarrier
Assume now that $r'_1\neq r'_2$, as in Figure \ref{step5finale}(a) where we represent $r'_1$ by $\mycircleerreprimo$ and $r'_2$ by $\mycirclewhite$. We may repeat the discussion above by considering the tile centered in $v_1$, and performing, if possible, another zero-cost flip, and so on. This procedure necessarily ends, see, e.g., Figure \ref{step5finale}(c). Note that the vertex $u$ may coincide with $v_3$. This concludes the proof of Step 5.4.

\noindent Finally, in view of the discussion above, the stable tiles (h) and (i) belong to a stable configuration only when they belong to a strip of thickness one and the stable tiles of type (f) does not belong to any stable configuration.
$\qed$

We are now ready to prove Proposition \ref{ricorrenzaneg}. 

\emph{Proof of Proposition \ref{ricorrenzaneg}(Estimate on the stability level).} In order to prove the recurrence property it is enough to focus on the configurations belonging to $\widetilde{\mathscr M}_{\text{neg}}:=(\mathscr M_{\text{neg}}\backslash\{\bold 1,\dots,\bold q\})\cup\bar{\mathscr M_{\text{neg}}}$. For any $\eta\in\widetilde{\mathscr M}_{\text{neg}}$ we prove that $V_\eta^\text{neg}$ is smaller than or equal to $V^*:=2 <\Gamma_{\text{neg}}(\bold 1,\mathcal X^s_{\text{neg}})$. 
\begin{figure}[h!]
\centering
\begin{tikzpicture}[scale=0.45, transform shape]
\fill[black!40!white] (0.3,0)rectangle(1.2,2.7);
\fill[black!75!white] (1.2,0)rectangle(2.4,2.7);
\fill[black!15!white] (3.3,0)rectangle(3.6,2.7)(0,0)rectangle(0.3,2.7);
\draw[step=0.3cm,color=black] (0,0) grid (3.6,2.7);
\draw (1.8,-0.1) node[below] {\large(a)};
\end{tikzpicture}\ \
\begin{tikzpicture}[scale=0.45, transform shape]
\fill[black!40!white] (9.3,0.3)rectangle(10.2,0.9);
\fill[black!15!white] (10.8,0.3)rectangle(11.7,1.2)(10.2,1.8)rectangle(11.7,2.4);
\fill[black!60!white] (9,1.5)rectangle(9.6,2.4)(12.3,1.5)rectangle(12.6,2.4);
\draw[step=0.3cm,color=black] (9,0) grid (12.6,2.7);
\draw (10.8,-0.1) node[below] {\large(b)};
\end{tikzpicture}\ \
\begin{tikzpicture}[scale=0.45, transform shape]
\fill[black!40!white] (0.3,1.2)rectangle(1.8,1.8);
\fill[black!15!white] (0.9,1.8)rectangle(1.5,2.7);
\fill[black!50!white] (1.2,0.3)rectangle(3.3,0.9);
\fill[black!70!white] (2.4,01.5)rectangle(3.6,0.9);
\draw[step=0.3cm,color=black] (0,0) grid (3.6,2.7);
\draw (1.8,-0.1) node[below] {\large(c)};
\end{tikzpicture}\ \
\begin{tikzpicture}[scale=0.45, transform shape]
\fill[white] (4.5,0)rectangle(6.6,2.7);
\fill[black!40!white] (6.6,0)rectangle(8.1,2.7);
\fill[black!80!white] (5.4,0.6)rectangle(6.6,2.1);
\draw[step=0.3cm,color=black] (4.5,0) grid (8.1,2.7);
\draw (6.3,-0.1) node[below] {\large(d)};
\end{tikzpicture}\ \
\begin{tikzpicture}[scale=0.45, transform shape]
\fill[black!40!white] (4.5,0)rectangle(5.7,1.5)(6,1.5)rectangle(7.2,2.7);
\fill[black!80!white] (4.5,2.7)rectangle(6,1.5);
\fill[black!25!white] (5.7,0)rectangle(7.2,0.9);
\fill[black!60!white] (7.2,0.6)rectangle(8.1,2.1);
\fill[black!12!white] (7.2,0)rectangle(8.1,0.6)(7.2,2.1)rectangle(8.1,2.7)(5.7,0.9)rectangle(7.2,1.5);
\draw[step=0.3cm,color=black] (4.5,0) grid (8.1,2.7);
\draw (6.3,-0.1) node[below] {\large(e)};
\end{tikzpicture}\ \
\caption{\label{figurerecurrenceproplabel} Examples of local minima of the Hamiltonian \eqref{hamiltonianneg}. We color white the vertices with spin $1$ and we use the other colors to denote the other spins $2,\dots,q$.}
\end{figure}
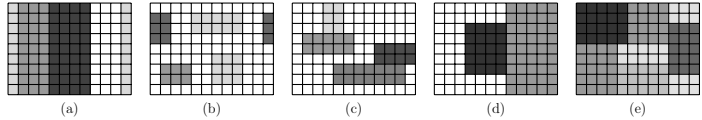\FloatBarrier
Let us first give an outline of the proof. First, we estimate of the stability level of those configurations in $\widetilde{\mathscr M}_{\text{neg}}$ that have at least two adjacent strips of different spins, see Figure \ref{figurerecurrenceproplabel}(a) and (d). 
Second, we estimate the stability level of those configurations in $\widetilde{\mathscr M}_{\text{neg}}$ that have at least an $s$-rectangle ($s\neq1$) either in a sea of spins $1$ or inside a cluster of spins $1$, as well as those configurations in which there is at least an $s$-rectangle ($s\neq1$) having a side such that on the corners there are stable tiles of type (m) and elsewhere there are stable tiles of type (d). See Figure \ref{figurerecurrenceproplabel}(b),(c) and (d).  Third, we consider those local minima in which at least an $r$-cluster has a side completely adjacent to a side of an $s$-cluster, see for instance Figure \ref{figurerecurrenceproplabel}(c)--(e) and Figure \ref{locmin4case}. Finally, we focus on those local minima that do not belong to any of the cases above, that is, those local minima with at least an $s$-rectangle with each side adjacent both to an $r$-cluster ($r\neq s$) and to a $1$-rectangle, see for instance \ref{figurerecurrenceproplabel}(e).

\textbf{Case 1.} Let us begin by assuming that $\eta$ has either at least two horizontal or vertical strips. Consider the case depicted in Figure \ref{figurerecurrenceproplabel}(a). Assume that $\eta$ has an $r$-strip $a\times K$ adjacent to an $s$-strip $b\times K$, $a,b\in\mathbb Z$, $a,b\ge 1$. Assume that $r,s\in S$, $s\neq 1$. Let $\bar\eta$ be the configuration obtained from $\eta$ by flipping from $r$ to $s$ all the spins $r$ belonging to the $r$-strip. We define a path $\omega:\eta\to\bar\eta$ as the concatenation of $a$ paths $\omega^{(1)},\dots,\omega^{(a)}$. Let $\omega:\eta\to\bar\eta$ be the path that flips the spins in the $r$-strip to $s$, column by column, starting from the column adjacent to the $s$-strip. Number the columns of the $r$-strip in order of flipping, and let $\omega^{(i)}:=(\omega^{(i)}_0=\eta_{i-1},\omega^{(i)}_1,\dots,\omega^{(i)}_{K}=\eta_{i})$ be the path that flips the $r$ spins in the $i$-th column. Then, for $i=1,\dots,a-1$,
\begin{align}\label{firsteqrecrecurrence}
H_{\text{neg}}(\omega^{(i)}_j)-H_{\text{neg}}(\omega^{(i)}_{j-1})=
\begin{cases}
2-h\mathbbm 1_{\{r=1\}},\ &\text{if}\ j=1;\\
-h\mathbbm 1_{\{r=1\}},\ &\text{if}\ j=2,\dots,K-1;\\
-2-h\mathbbm 1_{\{r=1\}},\ &\text{if}\ j=K.
\end{cases}
\end{align}
For any $i=1,\dots,a-1$, the maximum energy value along $\omega^{(i)}$ is reached at the first step. Computing the energy values along the sub-path $\omega^{(a)}$, that flips the last $r$-column, requires more care. Denoting by $v_i$ the vertex whose spin is flipping at the step $i$, 
\begin{align}
&H_{\text{neg}}(\omega^{(a)}_1)-H_{\text{neg}}(\omega^{(a-1)}_{K-1})=
\begin{cases}
1-h\mathbbm 1_{\{r=1\}},\ &\text{if $n_s(v_1)=1$,}\\ 
-h\mathbbm 1_{\{r=1\}}\ &\text{if $n_s(v_1)=2$,}
\end{cases}
\end{align}
and, if $i=2,\dots,K$, 
\begin{align}
&H_{\text{neg}}(\omega^{(a)}_i)-H_{\text{neg}}(\omega^{(a)}_{i-1})=
\begin{cases}
-1-h\mathbbm 1_{\{r=1\}},\ &\text{if $n_s(v_i)=2$,}\\ 
-2-h\mathbbm 1_{\{r=1\}},\ &\text{if $n_s(v_i)=3$.}  
\end{cases} \label{lasteqrecrecu}
\end{align}
In view of the above construction, $H_\text{neg}(\eta)>H_\text{neg}(\bar\eta)$ and, by comparing \eqref{firsteqrecrecurrence}--\eqref{lasteqrecrecu}, $V_\eta^\text{neg}\le2=V^*$.

\textbf{Case 2.} Let us now consider $\eta$ characterized by a sea of spins $1$ with some no-nteracting $s$-rectangles ($s\neq 1$). We distinguish the following cases:

(i) $\eta$ has at least a rectangle $R_{\ell_1\times\ell_2}$ of spins $s$, for some $s\in\{2,\dots,q\}$, with its minimum side of length $\ell:=\min\{\ell_1,\ell_2\}$ larger than or equal to $\ell^*$;

(ii) $\eta$ has only rectangles $R_{\ell_1\times\ell_2}$ of spins $s$, for some $s\in\{2,\dots,q\}$, with a side of length $\ell$ smaller than $\ell^*$.

\noindent In case (i), we construct a path $\omega=(\omega_0,\dots,\omega_{\ell-1})$, where $\omega_0=\eta$ and $\omega_{\ell-1}=:\tilde\eta$, that flips consecutively from $1$ to $s$ those spins adjacent to a side of length $\ell\ge\ell^*$. We have
\begin{align}
&H_{\text{neg}}(\omega_1)-H_{\text{neg}}(\eta)=2-h,\label{eqrefrecurrencerec}\\
&H_{\text{neg}}(\omega_i)-H_{\text{neg}}(\omega_{i-1})=-h, \ \text{for}\ i=2,\dots,\ell-2.\label{eqrefrecurrencerecbis}
\end{align}
It follows that $H_{\text{neg}}(\tilde\eta)-H_{\text{neg}}(\eta)=2-h\ell$. If $\ell>\ell^*=\left\lceil \frac{2}{h} \right\rceil$, then $2-h\ell<0$. Therefore the maximum energy is reached  at the first step and by \eqref{eqrefrecurrencerec} we get $V_\eta^\text{neg}=2-h<V^*$. Otherwise, if $\eta$ has only rectangles $R_{\ell^*\times\ell^*}$ of spins $s$, then $\tilde\eta$ has a rectangle $R_{\ell^*\times(\ell^*+1)}$ of spins $s$. Now, either this $s$-rectangle does not interact with the other rectangles of $\tilde\eta$ or it interacts with another rectangle $\hat R$. In the former case we conclude by arguing as previously since $\ell^*+1>\ell^*$. In the latter case, we have the following two possibilities
\begin{enumerate}
\item[(1)] $\hat R$ is an $s$-rectangle,
\item[(2)] $\hat R$ is an $r$-rectangle with $r\notin\{1,s\}$.
\end{enumerate}
In case (1), we define a configuration $\hat\eta$ from $\tilde\eta$ by flipping a spin $1$ to $s$ in the interaction interface. In particular, 
\begin{align}\label{gapenelastflip}
H_\text{neg}(\hat\eta)-H_\text{neg}(\tilde\eta)=h.
\end{align}
Hence, the maximum energy along $(\eta,\omega_1,\dots,\omega_{\ell-2},\tilde\eta,\hat\eta)$  is reached at the first step and we conclude that $V_\eta^\text{neg}=2-h < V^*$.

\noindent Let us now focus on case (2). We have to consider the two cases depicted in Figure \ref{figurarecurrencei}. Let $v_1,\dots,v_{\ell^*+1}$ be the vertices next to the side of length $\ell^*+1$ of the $s$-rectangle such that $v_1$ has two nearest neighbors with spin $1$, one nearest neighbor $s$ and one nearest neighbor inside the $r$-rectangle $\hat R$.
\begin{figure}[h!]
\centering
\begin{tikzpicture}[scale=0.7,transform shape]
\fill[black!80!white] (0,0) rectangle (1.5,1.8) (4.2,0) rectangle (5.7,1.8);
\fill[black!15!white] (1.8,0.9) rectangle (2.7,2.4) (6,1.5) rectangle (6.9,3);
\draw[step=0.3cm,color=black] (0,0) grid (1.5,1.8) (1.8,0.9) grid (2.7,2.4) (4.2,0) grid (5.7,1.8) (6,1.5) grid (6.9,3);
\draw[step=0.3cm,color=black] (1.5,0) grid (1.8,1.8) (5.7,0) grid (6,1.8);
\draw (1.35,-0.1) node[below] {\large(a)} (5.55,-0.1) node[below] {\large(b)};
\foreach \i in {1.65,5.85}\draw (\i,0.15) node[scale=0.6] {$v_6$};
\foreach \i in {1.65,5.85}\draw (\i,0.45) node[scale=0.6] {$v_5$};
\foreach \i in {1.65,5.85}\draw (\i,0.75) node[scale=0.6] {$v_4$};
\foreach \i in {1.65,5.85}\draw (\i,1.05) node[scale=0.6] {$v_3$};
\foreach \i in {1.65,5.85}\draw (\i,1.35) node[scale=0.6] {$v_2$};
\foreach \i in {1.65,5.85}\draw (\i,1.65) node[scale=0.6] {$v_1$};
\end{tikzpicture}
\caption{\label{figurarecurrencei} Examples of interacting rectangles in $\tilde\eta$ when $\ell^*=5$. We color gray the $r$-rectangle $\hat R$ and black the $s$-rectangle.}
\end{figure}
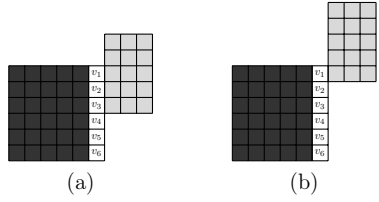\FloatBarrier
In the case depicted in Figure \ref{figurarecurrencei}(a), we define $\hat\eta_1:={\tilde\eta}^{(v_1,s)}$ and $\hat\eta_2:={\hat\eta_1}^{(v_2,s)}$. In particular,
\begin{align}
&H_{\text{neg}}(\hat\eta_1)-H_{\text{neg}}(\tilde\eta)=1-h,\label{lasttwoflipsone}\\
&H_{\text{neg}}(\hat\eta_2)-H_{\text{neg}}(\hat\eta_1)=-1-h.\label{lasttwoflipsbis}
\end{align}
Hence, from \eqref{eqrefrecurrencerec}--\eqref{eqrefrecurrencerecbis} and \eqref{lasttwoflipsone}--\eqref{lasttwoflipsbis}, we have that $H_{\text{neg}}(\hat\eta_2)-H_{\text{neg}}(\eta)=2-h(\ell^*+2)<2-h\ell^*\le0$.
Moreover, in view of \eqref{eqrefrecurrencerec} and \eqref{lasttwoflipsone}, along the path $(\eta,\omega_1,\dots,\omega_{\ell-2},\tilde\eta,\hat\eta_1,\hat\eta_2)$, we get that the maximum energy is reached at the first step. Hence, $V_\eta^\text{neg}=2-h < V^*$.

\noindent On the other hand, in the case depicted in Figure \ref{figurarecurrencei}(b) we define $\hat\eta_1:={\tilde\eta}^{(v_1,s)}$ and $\hat\eta_i:={\hat\eta_{i-1}}^{(v_i,s)}$ for any $i=2,\dots,\ell^*+1$. Note that 
\begin{align}
&H_{\text{neg}}(\hat\eta_1)-H_{\text{neg}}(\tilde\eta)=1-h,\label{lasttwoflipsonecasefigb}\\
&H_{\text{neg}}(\hat\eta_{i+1})-H_{\text{neg}}(\hat\eta_i)=-h,\ i=1,\dots,\ell^*.\label{lasttwoflipsbiscasefigb}
\end{align}
Hence, from \eqref{eqrefrecurrencerec}, \eqref{eqrefrecurrencerecbis}, \eqref{lasttwoflipsonecasefigb} and \eqref{lasttwoflipsbiscasefigb}, we have $H_{\text{neg}}(\hat\eta_{\ell^*+1})-H_{\text{neg}}(\eta)=3-h(2\ell^*+1)<0$. Moreover, by comparing \eqref{eqrefrecurrencerec} and \eqref{lasttwoflipsonecasefigb} along the path $(\eta,\omega_1,\dots,\omega_{\ell-2},\tilde\eta,\hat\eta_1,\dots,\hat\eta_{\ell^*+1})$  the maximum energy is reached at the first step. Hence, $V_\eta^\text{neg}=2-h < V^*$.

\noindent Now, we focus on the case (ii). We define a path  $\omega=(\omega_0,\dots,\omega_{\ell-1})$ that flips consecutively from $s$ to $1$ those spins $s$ next to a side of length $\ell<\ell^*$. We get: 
\begin{align}
&H_{\text{neg}}(\omega_i)-H_{\text{neg}}(\omega_{i-1})=h, \ \text{for}\ i=1,\dots,\ell-2;\\
&H_{\text{neg}}(\omega_{\ell-1})-H_{\text{neg}}(\omega_{\ell-2})=-(2-h).
\end{align}
Hence the maximum energy is achieved after $\ell-1$ steps and $H_{\text{neg}}(\omega_{\ell-1})-H_{\text{neg}}(\omega_0)=h(\ell-1)<2-h<V^*$.

\textbf{Case 3.} Let us now assume that $\eta$ has an $s$-rectangle $\bar R:=R_{a\times b}$ and an $r$-rectangle $\tilde R:=R_{c\times d}$ such that $\bar R$ has a side of length $a$ adjacent to a side of $\tilde R$ of length $c\ge a$, see for instance Figure \ref{figurerecurrenceproplabel}(e). The case $c<a$ may be studied by interchanging the role of spins $s$ and $r$. Given $\bar\eta$ the configuration obtained from $\eta$ by flipping to $r$ all the spins $s$ belonging to $\bar R$, we construct a path $\omega:\eta\to\bar\eta$ as the concatenation of $b$ paths $\omega^{(1)},\dots,\omega^{(b)}$. Let $\omega:\eta\to\bar\eta$ be the path that flips the spins in the $r$-rectangle $\tilde R$ to $s$, side by side, starting from the side adjacent to the $s$-rectangle $\bar R$. Number the sides of $\tilde R$ in order of flipping, and let $\omega^{(i)}:=(\omega^{(i)}_0=\eta_{i-1},\omega^{(i)}_1,\dots,\omega^{(i)}_{a}=\eta_{i})$ be the path that flips the $r$ spins in the $i$-th side. Then, for $i=1,\dots,b-1$,
\begin{align}\label{lastrecurrencealign}
H_{\text{neg}}(\omega^{(i)}_j)-H_{\text{neg}}(\omega^{(i)}_{j-1})=
\begin{cases}
1,\ &\text{if}\ j=1;\\
0,\ &\text{if}\ j=2,\dots,a-1;\\
-1,\ &\text{if}\ j=a.
\end{cases}
\end{align}

For any $i=1,\dots,b-1$, $H_{\text{neg}}(\eta)=H_{\text{neg}}(\eta_i)$ and the maximum energy value along $\omega^{(i)}$ is reached at the first step. Computing the energy values along the sub-path $\omega^{(b)}$, that flips the last $r$-side of the initial $\tilde R$, requires more care. Denoting by $v_i$ the vertex whose spin is flipping at the step $i$
\begin{align}
&H_{\text{neg}}(\omega^{(b)}_1)-H_{\text{neg}}(\eta_{b-1})=\begin{cases}0,\ &\text{if $n_r(v_1)=1,$ 
}\\ -1,\ &\text{if $n_r(v_1)=2,$ 
}  \end{cases}\\
&H_{\text{neg}}(\omega^{(b)}_i)-H_{\text{neg}}(\omega^{(b)}_{i-1})=\begin{cases}-1,\ &\text{if $n_r(v_i)=2,$ 
}\\ -2,\ &\text{if $n_r(v_i)=3$, 
}  \end{cases} \label{lastrecurrencealignbis}
\end{align}
for all $i=2,\dots,a$. In view of the above construction, $\Phi_{\omega}^{\text{neg}}=H_{\text{neg}}(\eta)+1$. Furthermore since $a\ge2$, $H_{\text{neg}}(\eta)>H_{\text{neg}}(\bar\eta)$ and, by comparing \eqref{lastrecurrencealign}--\eqref{lastrecurrencealignbis} we have $V_\eta^\text{neg}=1<V^*$.

\textbf{Case 4} Finally, let us consider $\eta$ with at least an $s$-rectangle, say $\hat R$, with each side adjacent both to an $r$-cluster, $r\neq s$, and to a $1$-rectangle. Let $\ell$ be the length of the interface between $\hat R$ and the $1$-rectangle and let $\omega=(\omega_0=\eta,\dots,\omega_\ell)$ be the path that flips from $1$ to $s$ all the spins $1$ on the $\ell$ vertices that lie on the interaface between $\hat R$ and the $1$-rectangle. We have that
\begin{align}
H_\text{neg}(\omega_i)-H_\text{neg}(\omega_{i-1})=
\begin{cases}
1-h,\ &\text{if}\ i=1;\\
-h,\ &\text{if}\ i=2,\dots,\ell-1;\\
-1-h,\ &\text{if}\ i=\ell.
\end{cases}
\end{align}
Since $H_\text{neg}(\omega_\ell)-H_\text{neg}(\eta)=-h\ell<0$ and
 $\Phi^\text{neg}_\omega=H_\text{neg}(\eta)+1-h$, we get $V^\text{neg}=1-h<V^*$.
 $\qed$
 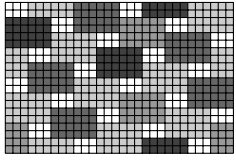
\begin{figure}
    \centering
    \begin{tikzpicture}[scale=0.4,transform shape]
    \fill[black!60!white] (1.5,0.5) rectangle (3,1.5)(0.75,2) rectangle (2.25,3)(2.25,4) rectangle (3.75,5)(5.25,3) rectangle (6.75,4)(6,1.5) rectangle (7.5,2.5);
    \fill[black!20!white] (0,1) rectangle (1.5,2)(0,2.5) rectangle (0.75,3.5)(0.75,4.5) rectangle (2.25,5)(0.75,0) rectangle (2.25,0.5)(2.25,1.5) rectangle (3.75,2.5)(1.5,3) rectangle (3,4)(3,0) rectangle (4.5,1)(5.25,0.5) rectangle (6.75,1.5)(4.5,2) rectangle (6,3)(6,4) rectangle (7.5,5)(6.75,2.5) rectangle (7.5,3.5);
    \fill[black!40!white] (0,0) rectangle (0.75,1)(6.75,0) rectangle (7.5,1)(3.75,1) rectangle (5.25,2)(3.75,3.5) rectangle (5.25,4.5);
    \fill[black!75!white] (4.5,0) rectangle (6,0.5)(4.5,4.5) rectangle (6,5)(0,3.5) rectangle (1.5,4.5)(3,2.5) rectangle (4.5,3.5);
        \draw[step=0.25cm,color=black] (0,0) grid (7.5,5);
    \end{tikzpicture}
\caption{\label{locmin4case} Local minimum on a $30\times 20$ grid graph in which there are not any $s$-rectangle with at least a side neither completely adjacent to an $r$-cluster nor completely sourrounded by spins $1$. 
}
\end{figure}
\subsection{Communication height between stable configurations}\label{commheightstable}
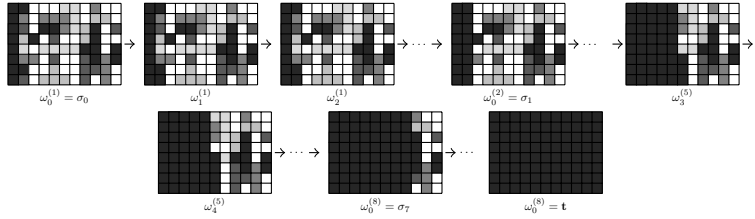
\begin{figure}[h!]
\centering
\begin{tikzpicture}[scale=0.45, transform shape]
\draw [fill=gray,gray] (0,0) rectangle (0.3,0.9) (0.3,0.6) rectangle (0.6,0.9) (2.7,0) rectangle (3,0.9) (2.4,2.1) rectangle (2.7,2.4) (0.6,0.9) rectangle (0.9,0.6) (0.9,1.8) rectangle (1.8,2.1) rectangle (1.5,1.5) (0.3,1.8)rectangle(0.6,2.1);
\draw [fill=gray,lightgray] (0.9,0) rectangle (1.8,0.3) (1.2,0.3) rectangle (1.5,0.6) (2.1,1.2) rectangle (2.7,1.8) rectangle (3,2.1) (2.4,1.8) rectangle (2.7,2.1) (0.3,1.2)rectangle(0.6,1.5);
\draw [fill=gray,black!15!white] (0,0.9)rectangle(0.6,1.2) (0.9,0.9)rectangle(2.4,1.2)(1.8,1.5)rectangle(2.1,2.4)(1.5,2.1)rectangle(1.8,2.4);
\draw [fill=gray,black!65!white] (0.3,0)rectangle(0.6,0.6) (0,1.2)rectangle(0.3,1.5)(0,1.5)rectangle(0.6,1.8)(2.1,0)rectangle(2.4,0.9)(2.4,0.3)rectangle(2.7,0.6) (1.2,0.6)rectangle(1.5,0.9)(1.2,1.2)rectangle(1.5,1.5) (3,1.5)rectangle(3.3,1.8);
\draw [fill=gray,black!85!white] (0,0) rectangle (0.3,2.4) (0.3,0.3) rectangle (0.9,0.6) (2.4,0.6) rectangle (2.7,1.8) (2.1,0.9) rectangle (2.7,1.2) (0.9,1.5)rectangle(1.5,1.8)(3,0.6)rectangle(3.3,1.2) (0.3,2.1)rectangle(0.6,2.4)(0.6,1.2)rectangle(0.9,1.5);
\fill [white] (0.6,0.3)rectangle(0.9,0.6);

\draw[step=0.3cm,color=black] (0,0) grid (3.3,2.4);
\draw (1.65,0) node[below] {$\omega^{(1)}_0=\sigma_0$};
\draw[->] (3.4,1.2)--(3.7,1.2);
\end{tikzpicture}
\begin{tikzpicture}[scale=0.45, transform shape]
\draw [fill=gray,gray] (0,0) rectangle (0.3,0.9) (0.3,0.6) rectangle (0.6,0.9) (2.7,0) rectangle (3,0.9) (2.4,2.1) rectangle (2.7,2.4) (0.6,0.9) rectangle (0.9,0.6) (0.9,1.8) rectangle (1.8,2.1) rectangle (1.5,1.5) (0.3,1.8)rectangle(0.6,2.1);
\draw [fill=gray,lightgray] (0.9,0) rectangle (1.8,0.3) (1.2,0.3) rectangle (1.5,0.6) (2.1,1.2) rectangle (2.7,1.8) rectangle (3,2.1) (2.4,1.8) rectangle (2.7,2.1) (0.3,1.2)rectangle(0.6,1.5);
\draw [fill=gray,black!15!white] (0,0.9)rectangle(0.6,1.2) (0.9,0.9)rectangle(2.4,1.2)(1.8,1.5)rectangle(2.1,2.4)(1.5,2.1)rectangle(1.8,2.4);
\draw [fill=gray,black!65!white] (0.3,0)rectangle(0.6,0.6) (0,1.2)rectangle(0.3,1.5)(0,1.5)rectangle(0.6,1.8)(2.1,0)rectangle(2.4,0.9)(2.4,0.3)rectangle(2.7,0.6) (1.2,0.6)rectangle(1.5,0.9)(1.2,1.2)rectangle(1.5,1.5) (3,1.5)rectangle(3.3,1.8);
\draw [fill=gray,black!85!white] (0,0) rectangle (0.3,2.4) (0.3,0.3) rectangle (0.9,0.6) (2.4,0.6) rectangle (2.7,1.8) (2.1,0.9) rectangle (2.7,1.2) (0.9,1.5)rectangle(1.5,1.8)(3,0.6)rectangle(3.3,1.2) (0.3,2.1)rectangle(0.6,2.4)(0.6,1.2)rectangle(0.9,1.5)(0.3,0)rectangle(0.6,0.3);
\fill [white] (0.6,0.3)rectangle(0.9,0.6);

\draw[step=0.3cm,color=black] (0,0) grid (3.3,2.4);
\draw (1.65,0) node[below] {$\omega^{(1)}_1$};
\draw[->] (3.4,1.2)--(3.7,1.2);
\end{tikzpicture}
\begin{tikzpicture}[scale=0.45, transform shape]
\draw [fill=gray,gray] (0,0) rectangle (0.3,0.9) (0.3,0.6) rectangle (0.6,0.9) (2.7,0) rectangle (3,0.9) (2.4,2.1) rectangle (2.7,2.4) (0.6,0.9) rectangle (0.9,0.6) (0.9,1.8) rectangle (1.8,2.1) rectangle (1.5,1.5) (0.3,1.8)rectangle(0.6,2.1);
\draw [fill=gray,lightgray] (0.9,0) rectangle (1.8,0.3) (1.2,0.3) rectangle (1.5,0.6) (2.1,1.2) rectangle (2.7,1.8) rectangle (3,2.1) (2.4,1.8) rectangle (2.7,2.1) (0.3,1.2)rectangle(0.6,1.5);
\draw [fill=gray,black!15!white] (0,0.9)rectangle(0.6,1.2) (0.9,0.9)rectangle(2.4,1.2)(1.8,1.5)rectangle(2.1,2.4)(1.5,2.1)rectangle(1.8,2.4);
\draw [fill=gray,black!65!white] (0.3,0)rectangle(0.6,0.6) (0,1.2)rectangle(0.3,1.5)(0,1.5)rectangle(0.6,1.8)(2.1,0)rectangle(2.4,0.9)(2.4,0.3)rectangle(2.7,0.6) (1.2,0.6)rectangle(1.5,0.9)(1.2,1.2)rectangle(1.5,1.5) (3,1.5)rectangle(3.3,1.8);
\draw [fill=gray,black!85!white] (0,0) rectangle (0.3,2.4) (0.3,0.3) rectangle (0.9,0.6) (2.4,0.6) rectangle (2.7,1.8) (2.1,0.9) rectangle (2.7,1.2) (0.9,1.5)rectangle(1.5,1.8)(3,0.6)rectangle(3.3,1.2) (0.3,2.1)rectangle(0.6,2.4)(0.6,1.2)rectangle(0.9,1.5)(0.3,0)rectangle(0.6,0.6);
\fill [white] (0.6,0.3)rectangle(0.9,0.6);

\draw[step=0.3cm,color=black] (0,0) grid (3.3,2.4);
\draw (1.65,0) node[below] {$\omega^{(1)}_2$};
\draw[->] (3.4,1.2)--(3.7,1.2) node[right] {$\dots$};
\draw[->] (4.4,1.2)--(4.7,1.2);
\end{tikzpicture}
\begin{tikzpicture}[scale=0.45, transform shape]
\draw [fill=gray,gray] (0,0) rectangle (0.3,0.9) (0.3,0.6) rectangle (0.6,0.9) (2.7,0) rectangle (3,0.9) (2.4,2.1) rectangle (2.7,2.4) (0.6,0.9) rectangle (0.9,0.6) (0.9,1.8) rectangle (1.8,2.1) rectangle (1.5,1.5) (0.3,1.8)rectangle(0.6,2.1);
\draw [fill=gray,lightgray] (0.9,0) rectangle (1.8,0.3) (1.2,0.3) rectangle (1.5,0.6) (2.1,1.2) rectangle (2.7,1.8) rectangle (3,2.1) (2.4,1.8) rectangle (2.7,2.1) (0.3,1.2)rectangle(0.6,1.5);
\draw [fill=gray,black!15!white] (0,0.9)rectangle(0.6,1.2) (0.9,0.9)rectangle(2.4,1.2)(1.8,1.5)rectangle(2.1,2.4)(1.5,2.1)rectangle(1.8,2.4);
\draw [fill=gray,black!65!white] (0.3,0)rectangle(0.6,0.6) (0,1.2)rectangle(0.3,1.5)(0,1.5)rectangle(0.6,1.8)(2.1,0)rectangle(2.4,0.9)(2.4,0.3)rectangle(2.7,0.6) (1.2,0.6)rectangle(1.5,0.9)(1.2,1.2)rectangle(1.5,1.5) (3,1.5)rectangle(3.3,1.8);
\draw [fill=gray,black!85!white] (0,0) rectangle (0.3,2.4) (0.3,0.3) rectangle (0.9,0.6) (2.4,0.6) rectangle (2.7,1.8) (2.1,0.9) rectangle (2.7,1.2) (0.9,1.5)rectangle(1.5,1.8)(3,0.6)rectangle(3.3,1.2) (0.3,2.1)rectangle(0.6,2.4)(0.6,1.2)rectangle(0.9,1.5)(0.3,0)rectangle(0.6,2.4);
\fill [white] (0.6,0.3)rectangle(0.9,0.6);

\draw[step=0.3cm,color=black] (0,0) grid (3.3,2.4);
\draw[->] (3.4,1.2)--(3.7,1.2) node[right] {$\dots$};
\draw (1.65,0) node[below] {$\omega^{(2)}_0=\sigma_1$};
\end{tikzpicture}
\begin{tikzpicture}[scale=0.45, transform shape]
\draw[->] (-0.4,1.2)--(-0.1,1.2);
\draw [fill=gray,gray] (0,0) rectangle (0.3,0.9) (0.3,0.6) rectangle (0.6,0.9) (2.7,0) rectangle (3,0.9) (2.4,2.1) rectangle (2.7,2.4) (0.6,0.9) rectangle (0.9,0.6) (0.9,1.8) rectangle (1.8,2.1) rectangle (1.5,1.5) (0.3,1.8)rectangle(0.6,2.1);
\draw [fill=gray,lightgray] (0.9,0) rectangle (1.8,0.3) (1.2,0.3) rectangle (1.5,0.6) (2.1,1.2) rectangle (2.7,1.8) rectangle (3,2.1) (2.4,1.8) rectangle (2.7,2.1) (0.3,1.2)rectangle(0.6,1.5);
\draw [fill=gray,black!15!white] (0,0.9)rectangle(0.6,1.2) (0.9,0.9)rectangle(2.4,1.2)(1.8,1.5)rectangle(2.1,2.4)(1.5,2.1)rectangle(1.8,2.4);
\draw [fill=gray,black!65!white] (0.3,0)rectangle(0.6,0.6) (0,1.2)rectangle(0.3,1.5)(0,1.5)rectangle(0.6,1.8)(2.1,0)rectangle(2.4,0.9)(2.4,0.3)rectangle(2.7,0.6) (1.2,0.6)rectangle(1.5,0.9)(1.2,1.2)rectangle(1.5,1.5) (3,1.5)rectangle(3.3,1.8);
\draw [fill=gray,black!85!white] (0,0) rectangle (0.3,2.4) (0.3,0.3) rectangle (0.9,0.6) (2.4,0.6) rectangle (2.7,1.8) (2.1,0.9) rectangle (2.7,1.2) (0.9,1.5)rectangle(1.5,1.8)(3,0.6)rectangle(3.3,1.2) (0.3,2.1)rectangle(0.6,2.4)(0.6,1.2)rectangle(0.9,1.5)(0.3,0)rectangle(0.6,2.4)(0.6,0)rectangle(1.5,2.4)(1.5,0)rectangle(1.8,0.9);

\draw[step=0.3cm,color=black] (0,0) grid (3.3,2.4);
\draw (1.65,0) node[below] {$\omega^{(5)}_3$};
\draw[->] (3.4,1.2)--(3.7,1.2);
\end{tikzpicture}\\
\begin{tikzpicture}[scale=0.45, transform shape]
\draw [fill=gray,gray] (0,0) rectangle (0.3,0.9) (0.3,0.6) rectangle (0.6,0.9) (2.7,0) rectangle (3,0.9) (2.4,2.1) rectangle (2.7,2.4) (0.6,0.9) rectangle (0.9,0.6) (0.9,1.8) rectangle (1.8,2.1) rectangle (1.5,1.5) (0.3,1.8)rectangle(0.6,2.1);
\draw [fill=gray,lightgray] (0.9,0) rectangle (1.8,0.3) (1.2,0.3) rectangle (1.5,0.6) (2.1,1.2) rectangle (2.7,1.8) rectangle (3,2.1) (2.4,1.8) rectangle (2.7,2.1) (0.3,1.2)rectangle(0.6,1.5);
\draw [fill=gray,black!15!white] (0,0.9)rectangle(0.6,1.2) (0.9,0.9)rectangle(2.4,1.2)(1.8,1.5)rectangle(2.1,2.4)(1.5,2.1)rectangle(1.8,2.4);
\draw [fill=gray,black!65!white] (0.3,0)rectangle(0.6,0.6) (0,1.2)rectangle(0.3,1.5)(0,1.5)rectangle(0.6,1.8)(2.1,0)rectangle(2.4,0.9)(2.4,0.3)rectangle(2.7,0.6) (1.2,0.6)rectangle(1.5,0.9)(1.2,1.2)rectangle(1.5,1.5) (3,1.5)rectangle(3.3,1.8);
\draw [fill=gray,black!85!white] (0,0) rectangle (0.3,2.4) (0.3,0.3) rectangle (0.9,0.6) (2.4,0.6) rectangle (2.7,1.8) (2.1,0.9) rectangle (2.7,1.2) (0.9,1.5)rectangle(1.5,1.8)(3,0.6)rectangle(3.3,1.2) (0.3,2.1)rectangle(0.6,2.4)(0.6,1.2)rectangle(0.9,1.5)(0.3,0)rectangle(0.6,2.4)(0.6,0)rectangle(1.5,2.4)(1.5,0)rectangle(1.8,1.2);

\draw[step=0.3cm,color=black] (0,0) grid (3.3,2.4);
\draw (1.65,0) node[below] {$\omega^{(5)}_4$};
\draw[->] (3.4,1.2)--(3.7,1.2) node[right] {$\dots$};
\draw[->] (4.4,1.2)--(4.7,1.2);
\end{tikzpicture}
\begin{tikzpicture}[scale=0.45, transform shape]
\draw [fill=gray,gray] (0,0) rectangle (0.3,0.9) (0.3,0.6) rectangle (0.6,0.9) (2.7,0) rectangle (3,0.9) (2.4,2.1) rectangle (2.7,2.4) (0.6,0.9) rectangle (0.9,0.6) (0.9,1.8) rectangle (1.8,2.1) rectangle (1.5,1.5) (0.3,1.8)rectangle(0.6,2.1);
\draw [fill=gray,lightgray] (0.9,0) rectangle (1.8,0.3) (1.2,0.3) rectangle (1.5,0.6) (2.1,1.2) rectangle (2.7,1.8) rectangle (3,2.1) (2.4,1.8) rectangle (2.7,2.1) (0.3,1.2)rectangle(0.6,1.5);
\draw [fill=gray,black!15!white] (0,0.9)rectangle(0.6,1.2) (0.9,0.9)rectangle(2.4,1.2)(1.8,1.5)rectangle(2.1,2.4)(1.5,2.1)rectangle(1.8,2.4);
\draw [fill=gray,black!65!white] (0.3,0)rectangle(0.6,0.6) (0,1.2)rectangle(0.3,1.5)(0,1.5)rectangle(0.6,1.8)(2.1,0)rectangle(2.4,0.9)(2.4,0.3)rectangle(2.7,0.6) (1.2,0.6)rectangle(1.5,0.9)(1.2,1.2)rectangle(1.5,1.5) (3,1.5)rectangle(3.3,1.8);
\draw [fill=gray,black!85!white] (0,0) rectangle (0.3,2.4) (0.3,0.3) rectangle (0.9,0.6) (2.4,0.6) rectangle (2.7,1.8) (2.1,0.9) rectangle (2.7,1.2) (0.9,1.5)rectangle(1.5,1.8)(3,0.6)rectangle(3.3,1.2) (0.3,2.1)rectangle(0.6,2.4)(0.6,1.2)rectangle(0.9,1.5)(0.3,0)rectangle(0.6,2.4)(0.6,0)rectangle(2.4,2.4);

\draw[step=0.3cm,color=black] (0,0) grid (3.3,2.4);
\draw (1.65,0) node[below] {$\omega^{(8)}_0=\sigma_7$};
\draw[->] (3.4,1.2)--(3.7,1.2) node[right] {$\dots$};
\end{tikzpicture}
\begin{tikzpicture}[scale=0.45, transform shape]
\draw [fill=gray,black!85!white] (0,0) rectangle (3.3,2.4); 
\draw[step=0.3cm,color=black] (0,0) grid (3.3,2.4);
\draw (1.65,0) node[below] {$\omega^{(8)}_0=\bold t$};
\end{tikzpicture}
\caption{\label{illustrationexpalg} Illustration of some particular configurations belonging to the path $\omega:\sigma\to\bold t$ of Proposition \ref{expansionalg}. We color black those vertices whose spin is $t$.}
\end{figure}

In order to study the hitting time $\tau_\bold s^\bold 1$ of a stable configuration $\bold s\in\mathcal X^s_{\text{neg}}$, we first estimate the communication height $\Phi_{\text{neg}}(\bold r,\bold s)$ between two stable configurations $\bold r,\bold s\in\mathcal X^s_{\text{neg}}$, $\bold r\neq\bold s$. Indeed, during the transition $\mathbf 1\to\mathbf s$, the process may visit a stable state $\mathbf r\neq\mathbf s$ before hitting $\mathbf s$. Using \eqref{totnumberdisedges}, the energy difference between any $\sigma\in\mathcal X$ and any $\bold s\in\mathcal X^s_{\text{neg}}$ reads
\begin{align}\label{rewritegap}
H_{\text{neg}}(\sigma)-H_{\text{neg}}(\bold s)
= d_v(\sigma)+d_h(\sigma)+h\sum_{u\in V} \mathbbm{1}_{\{\sigma(u)=1\}}.
\end{align}
In \cite[Proposition 2.4]{nardi2019tunneling} the authors define the so-called \textit{expansion algorithm}. We rewrite this procedure in the proof of the next proposition by adapting it to our scenario. Indeed, it is different from \cite{nardi2019tunneling} since in our setting there is a non-zero external magnetic field.
\begin{proposition}[Expansion algorithm]\label{expansionalg}
If the external magnetic field is negative and if $\sigma\in\mathcal X$ has a $t$-bridge for some $t\in\{2,\dots,q\}$, then there exists a path $\omega: \sigma\to\bold t$ such that
$\Phi_\omega^\emph{neg}-H_{\emph{neg}}(\sigma)\le2.$
\end{proposition}
\textit{Proof.} Without loss of generality we assume that the first column $c_0$ is the $t$-bridge. Following an iterative procedure, we define a path $\omega: \sigma\to\bold t$ that flips all spins to $t$ column-by-column starting with column $c_0$. Formally, $\omega$ is  the concatenation of $L$ paths $\omega^{(1)},\dots,\omega^{(L)}$ with $\omega^{(i)}:=(\omega^{(i)}_0=\sigma_{i-1},\dots,\omega^{(i)}_K=\sigma_{i})$ and $\omega_j^{(i)}:=(\omega_{j-1}^{(i)})^{(u,t)}$, for $u:=(i,j-1)$ and $j=1,\dots,K$.
 In particular, $\sigma_0:=\sigma$, $\sigma_{L}:=\bold t$ and the configurations $\sigma_i$, $i=0,\dots,L$, are given by
\begin{align}
\sigma_i(v):=
\begin{cases}
t \ \ \ &\text{if}\ v \in \bigcup_{j=0}^i c_j,\\
\sigma(v) &\text{if}\ v \in V \backslash \bigcup_{j=0}^i c_j.
\end{cases}
\end{align}
See Figure \ref{illustrationexpalg} for an illustration of the construction above. 
\FloatBarrier
Let us now study the energy difference $H_{\text{neg}}(\omega_j^{(i)})-H_{\text{neg}}(\omega_{j-1}^{(i)})$ for $j=1,\dots,K$. It is immediate to see that if $\sigma(u)=t$, then $H_{\text{neg}}(\omega_j^{(i)})-H_{\text{neg}}(\omega_{j-1}^{(i)})=0$. Hence, assume that $\sigma(u)\neq t$. Using \eqref{energydifferenceneg} and counting the number of spins $s$ neighbors of $u$, we get 
\begin{align}\label{upperalignone}
H_{\text{neg}}(\omega_j^{(i)})-H_{\text{neg}}(\omega_{j-1}^{(i)})\le
\begin{cases}
2-h\mathbbm 1_{\{\omega_{j-1}^{(i)}(v)=1\}},\ &\text{if}\ j=1;\\
 -h\mathbbm 1_{\{\omega_{j-1}^{(i)}(v)=1\}},\ &\text{if}\ 1<j<K;\\
-2-h\mathbbm 1_{\{\omega_{j-1}^{(i)}(v)=1\}},\ &\text{if}\ j=K.
\end{cases}
\end{align}

For every $i=1,\dots,L-1$, the inequalities \eqref{upperalignone} imply that $\Phi_{\omega^{(i)}}^\text{neg}-H_{\text{neg}}(\sigma_{i-1})\le2$. Hence, the path $\omega: \sigma\to\bold{t}$ is such that $\Phi_{\omega}^\text{neg}-H_{\text{neg}}(\sigma)\le2$. $\qed$

Thanks to Proposition \ref{expansionalg} we are able to obtain an upper bound on $\Gamma_{\text{neg}}(\bold r,\bold s):=\Phi_{\text{neg}}(\bold r,\bold s)-H_{\text{neg}}(\bold r)$, for any  $\bold r,\bold s \in \mathcal X^s_{\text{neg}}$, $\bold r\neq\bold s$. 
\begin{proposition}[Upper bound for the stability level between two stable configurations]\label{upperboundcaso2ss}
If the external magnetic field is negative, then for any $\bold r,\bold s\in\mathcal X^s_{\emph{neg}}$, $\bold r\neq\bold s$, we have
\begin{align}\label{upperboundalign}
\Phi_{\emph{neg}}(\bold r,\bold s)-H_{\emph{neg}}(\bold r)\le 2\min\{K,L\}+2.
\end{align} 
\end{proposition}
\textit{Proof.} The proof is analogous to the one of \cite[Proposition 2.5]{nardi2019tunneling} by replacing the role of \cite[Proposition 2.4]{nardi2019tunneling} with Proposition \ref{expansionalg}. For the details we refer to the Appendix \ref{appendixstablephi}. $\qed$

Now let us estimate a lower bound for $\Gamma_{\text{neg}}(\bold r,\bold s)$, for any  $\bold r,\bold s \in \mathcal X^s_{\text{neg}}$, $\bold r\neq\bold s$. The following proposition is an adaptation of \cite[Proposition 2.7]{nardi2019tunneling} to the case of Potts model with external magnetic field. Recall that $B_s(\sigma)$ denotes the total number of vertical and horizontal $s$-bridges in $\sigma\in\mathcal X$, see Subsection \ref{subseclocgeo}.
\begin{proposition}[Lower bound for the stability level between two stable configurations]\label{lowerbound2ss}
If the external magnetic field is negative, then for every $\bold r,\bold s\in\mathcal X^s_{\emph{neg}}$, the following inequality holds
\begin{align}
\Phi_{\emph{neg}}(\bold r,\bold s)-H_{\emph{neg}}(\bold r)\ge 2\min\{K,L\}+2.
\end{align} 
\end{proposition}
\textit{Proof.} We show that along every path $\omega:\bold r\to\bold s$ in $\mathcal X$ there exists a configuration $\eta$ such that $H_{\text{neg}}(\eta)-H_{\text{neg}}(\bold r)\ge 2K+2.$ Consider a path $\omega=(\omega_1,\dots,\omega_n)$ with $\omega_1=\bold r$ and $\omega_n=\bold s$. Obviously, $B_s(\bold r)=0$ and $B_s(\bold s)=K+L$. Let $\omega_{\bar k}$ be the configuration along the path $\omega$ that is the first to have at least two $s$-bridges, i.e., $\bar k:=\min\{k\le n|\ B_s(\omega_k)\ge 2\}$. We claim that the configuration $\omega_{\bar k-1}$ is such that  
\begin{align}\label{claim3}
H_{\text{neg}}(\omega_{\bar k-1})-H_{\text{neg}}(\bold r)\ge 2K+2.
\end{align}
Let us prove this claim by studying separately the following three cases:
\begin{itemize}
\item[{(i)}] $\omega_{\bar k}$ has only vertical $s$-bridges,
\item[{(ii)}] $\omega_{\bar k}$ has only horizontal $s$-bridges,
\item[{(iii)}] $\omega_{\bar k}$ has at least one $s$-cross.
\end{itemize}
We study scenarios (i) and (iii), since scenario (ii) may be studied similarly as (i). Let us begin by assuming that (i) holds. From the definition of $\bar k$, it follows that $B_s(\omega_{\bar k-1})=1$ and $B_s(\omega_{\bar k})= 2$. Otherwise $\omega_{\bar k}$ would have an $s$-cross in view of \cite[Lemma 2.6]{nardi2019tunneling} and it would be a contradiction with (i). Let us assume that $\omega_{\bar k}$ has the two vertical $s$-bridges on columns $c$ and $\hat c$ and, without loss of generality, $\omega_{\bar k-1}$ has only one $s$-bridge on column $c$. In particular, in $\omega_{\bar k-1}$ all spins in $\hat c$ are $s$, except one which is different from $s$. Thus, in view of \cite[Lemma 2.3(d)]{nardi2019tunneling} we have
\begin{align}\label{gapcolumn}
d_{\hat c}(\omega_{\bar k-1})=2.
\end{align}
Moreover, it is easy to see that there are no horizontal bridges. Thanks to this fact and to \cite[Lemma 2.3(c)]{nardi2019tunneling}, we have $d_{ r_i}(\omega_{\bar k-1})\ge 2$ for every row $r_i$, $i=0,\dots,K-1$. Then,
\begin{align}\label{gaprow}
d_h(\omega_{\bar k-1})=\sum_{i=0}^{K-1} d_{ r_i}(\omega_{\bar k-1})\ge 2K.
\end{align}
From \eqref{rewritegap}, \eqref{gapcolumn} and \eqref{gaprow} we get that
\begin{align}
H_{\text{neg}}(\omega_{\bar k-1})-H_{\text{neg}}(\bold r)
&\ge 2 + 2K + h\sum_{u\in V} \mathbbm{1}_{\{\omega_{\bar k-1}(u)=1\}}\ge 2 + 2K.
\end{align}
Let us now focus on (iii). In this case $\omega_{\bar k}$ has at least one $s$-cross and, by definition of $\bar k$, $B_s(\omega_{\bar k-1})$ is either $0$ or $1$ and we study these two cases separately.

Assume $B_s(\omega_{\bar k-1})=0$. $\omega_{\bar k-1}$ has no $s$-bridges, then, by \cite[Lemma 2.6]{nardi2019tunneling}, $B_s(\omega_{\bar k})=2$ and $\omega_{\bar k}$ has exactly one $s$-cross. 
Let us assume that this $s$-cross lies on row $\hat{r}$ and on column $\hat c$. The horizontal and vertical $s$-bridges of $\omega_{\bar k}$ must have then been created simultaneously by updating the spin on the vertex $\hat{v}:=\hat{r}\cap\hat c$. Hence, we have $\omega_{\bar k-1}(\hat v)\neq s$, $\omega_{\bar k-1}(v)=s$, for all $v\in\hat{r}\cup\hat c$, $v\neq\hat{v}$, and $\omega_{\bar k}(\hat v)=s$.
Since there is a spin equal to $s$ in every row and in every column, $\omega_{\bar k-1}$ has no $t$-bridges ($t\neq s$). Since by assumption $B_s(\omega_{\bar k-1})=0$, $\omega_{\bar k-1}$ has no bridges of any spin. Therefore, from \cite[Lemma 2.3(c)--(d)]{nardi2019tunneling} follows that
\begin{align}\label{gaprandc}
d_h(\omega_{\bar k-1})=\sum_{i=0}^{K-1} d_{ r_i}(\omega_{\bar k-1})\ge 2K\ \ \text{and}\ \ d_v(\omega_{\bar k-1})=\sum_{j=0}^{L-1} d_{c_j}(\omega_{\bar k-1})\ge 2L.
\end{align}
Plugging  \eqref{gaprandc} in \eqref{rewritegap}, we conclude that
\begin{align}
H_{\text{neg}}(\omega_{\bar k-1})-H_{\text{neg}}(\bold r)\ge 2L+2K> 2\min\{K,L\}+2=2K+2.
\end{align}

Assume now $B_s(\omega_{\bar k-1})=1$. In this case, $\omega_{\bar k-1}$ has an unique $s$-bridge and we assume that such a bridge is vertical and lies on column $\tilde c$. In view of \cite[Lemma 2.2]{nardi2019tunneling}, there are no horizontal $t$-bridges in $\omega_{\bar k-1}$ ($t\neq s$). Hence, $\omega_{\bar k-1}$ has no horizontal bridges and by \cite[Lemma 2.3(c)]{nardi2019tunneling} we get 
\begin{align}\label{gaprows2}
d_h(\omega_{\bar k-1})=\sum_{i=0}^{K-1} d_{ r_i}(\omega_{\bar k-1})\ge 2K.
\end{align}
Moreover, $\omega_{\bar k}$ has a unique horizontal $s$-bridge, say on row $\hat{r}$. Hence, if $\hat{v}$ is the vertex where $\omega_{\bar k-1}$ and $\omega_{\bar k}$ differ, $\hat v$ must lie in $\hat r$ and $\omega_{\bar k-1}(\hat{v})\neq s$ and $\omega_{\bar k-1}(v)=s$, $\forall v \in \hat{ r}$, $v \neq \hat{v}$, and $\omega_{\bar k}(\hat v)=s$. Let $\hat c$ be the column where $\hat{v}$ lies. \cite[Lemma 2.3(d)]{nardi2019tunneling} implies that $d_c(\omega_{\bar k-1})\ge2$ for any column $c\neq\tilde c,\ \hat c$. Then,
\begin{align}\label{gapcolumns2}
d_v(\omega_{\bar k-1})=\sum_{j=0}^{L-1} d_{c_j}(\omega_{\bar k-1})\ge 2L-4.
\end{align}
In view of \eqref{rewritegap}, \eqref{gaprows2} and \eqref{gapcolumns2} it follows that
\begin{align}
H_{\text{neg}}(\omega_{\bar k-1})-H_{\text{neg}}(\bold r)\ge 2L+2K-4>2\min\{K,L\}+2=2K+2,
\end{align}
where the second inequality holds because $L\ge K\ge 3\ell^*>3$.
$\qed$

\subsection{Energy landscape: proof of the main results}\label{proofenergylandscape}
We are now able to prove Theorem \ref{theoremcomparisonneg}.\\
\textit{Proof of Theorem \ref{theoremcomparisonneg}}. Let us begin by recalling that for any $\bold r,\bold s\in\mathcal X^s_{\text{neg}}$, $\bold r\neq\bold s$, from Theorem \ref{teometastableneg} we have 
\begin{align}
\Gamma_{\text{neg}}(\bold 1,\mathcal X^s_{\text{neg}})=\Phi_{\text{neg}}(\bold 1,\mathcal X^s_{\text{neg}})-H_{\text{neg}}(\bold 1)= 4\ell^*-h(\ell^*(\ell^*-1)+1),
\end{align}
and, from Proposition \ref{upperboundcaso2ss} and Proposition \ref{lowerbound2ss},
\begin{align}
\Gamma_{\text{neg}}(\bold r,\bold s)=\Phi_{\text{neg}}(\bold r,\bold s)-H_{\text{neg}}(\bold r)= 2\min\{K,L\}+2.
\end{align}
For any $\bold r,\bold s\in\mathcal X^s_{\text{neg}}$, $\bold r\neq\bold s$, first we show that $\Gamma_{\text{neg}}(\bold 1,\mathcal X^s_{\text{neg}}) < \Gamma_{\text{neg}}(\bold s, \bold r)$. Indeed, given $0<h<1$ and $L\ge K\ge3\ell^*$, we have
\begin{align}
\Gamma_{\text{neg}}&(\bold 1,\mathcal X^s_{\text{neg}}) - \Gamma_{\text{neg}}(\bold r,\bold s)=4\ell^*-h(\ell^*(\ell^*-1)+1)-(2K+2)\notag \\
&\le 4\ell^*-h(\ell^*(\ell^*-1)+1) -6\ell^*-2 
<-2\ell^*-h(\ell^*)^2+\ell^*-h-2<0\implies \eqref{comparisongammaneg}.\notag 
\end{align}
Furthermore, by Assumption \ref{remarkconditionneg}, $\Phi_{\text{neg}}(\bold r,\bold s)$ is smaller than or equal to $\Phi_{\text{neg}}(\bold 1,\bold s)$, for any $\bold r,\bold s\in\mathcal X^s_{\text{neg}}$, $\bold r\neq\bold s$. Indeed, since $|V|=KL$, 
\begin{align}\label{iffinizio}
\Phi_{\text{neg}}(\bold r,\bold s)-\Phi_{\text{neg}}(\bold 1,\bold s)&=2K+2+H_{\text{neg}}(\bold r)-4\ell^*+h(\ell^*)^2-h\ell^*+h-H_{\text{neg}}(\bold 1)\notag\\
&=2K+2-|E|-4\ell^*+h(\ell^*)^2-h\ell^*+h+|E|-h|V|\notag\\
&=2K+2-4\ell^*+h(\ell^*)^2-h\ell^*+h-hKL
\end{align}
for any $\bold r,\bold s\in\mathcal X^s_{\text{neg}}$, $\bold r\neq\bold s$. Since $\ell^*=\left\lceil \frac{2}{h} \right\rceil$, we can write $\ell^*=\frac{2}{h}+1-\delta$ where $0<\delta<1$ denotes the fractional part of $2/h$, that is not integer in view of Assumption \ref{remarkconditionneg}(ii). Assume by contradiction that \eqref{comparisonalignneg} is false, i.e., 
\begin{align}\label{hpassurdophi}
\Phi_{\text{neg}}(\bold r,\bold s)\ge\Phi_{\text{neg}}(\bold 1,\bold s). 
\end{align}
Using \eqref{iffinizio}, we have that \eqref{hpassurdophi} is verified if and only if
\begin{align}
& 2K+2-4\ell^*+h(\ell^*)^2-h\ell^*+h\ge hKL \notag \\
\iff\ & 2K+2-4(\frac 2 h+1-\delta)+h(\frac 2 h+1-\delta)^2-h(\frac 2 h+1-\delta)+h\ge hKL\notag\\
\iff\ & \frac 2 h K +\frac 2 h-\frac 4 h(\frac 2 h+1-\delta)+(\frac 2 h+1-\delta)^2-\frac 2 h-1+\delta+1\ge KL\notag\\
\iff\ & \frac 2 h K +\frac 2 h-\frac 8 {h^2}-\frac4 h+\frac 4 h \delta+\frac 4 {h^2}+1+\delta^2+\frac 4 h-\frac 4 h \delta-2\delta-\frac 2 h-1+\delta+1\ge KL\notag\\
\iff\ & \frac 2 h K -\frac 4 {h^2}+1+\delta^2-\delta\ge KL. \label{ifflast}
\end{align}
Since $L\ge K\ge 3\ell^*$ and since $0<\delta<1$, we have that 
\begin{align}\label{estimateKL}
KL\ge 3K\ell^*=3K(\frac 2 h+1-\delta)=\frac 6 h K+3K-3K\delta>\frac 6 h K.
\end{align}
Moreover, since $0<\delta<1$ implies that $\delta^2-\delta<0$, we have that
\begin{align}\label{stimadeltaaltosinistra}
\frac 2 h K -\frac 4 {h^2}+1+\delta^2-\delta<\frac 2 h K -\frac 4 {h^2}+1.
\end{align}
Combining \eqref{ifflast}, \eqref{estimateKL} and \eqref{stimadeltaaltosinistra}, since $0<\delta<1$, approximately we get that \eqref{hpassurdophi} is satisfied if and only if
\begin{align}
\frac 2 h K -\frac 4 {h^2}+1>\frac 6 h K\ \iff\ -\frac 4 h K -\frac 4 {h^2}+1>0,\label{disfinaleiff}
\end{align}
that is a contradiction. Indeed, the l.h.s. of \eqref{disfinaleiff} is strictly negative since Assumption \ref{remarkconditionneg}(i), i.e., $0<h<1$, implies that $-\frac 4 {h^2}+1<0$. Hence, we conclude that \eqref{comparisonalignneg} is satisfied. \\
Finally, we prove \eqref{gammatildenegtime}.  By \cite[Lemma 3.6]{nardi2016hitting} we get that $\widetilde{\Gamma}_\text{neg}(\mathcal X\backslash\{\bold s\})$ is the maximum energy that the process started in $\eta\in\mathcal X\backslash\{\bold s\}$ has to overcome in order to arrive in $\bold s$, i.e. 
\begin{align}
\widetilde{\Gamma}_\text{neg}(\mathcal X\backslash\{\bold s\})=\max_{\eta\in\mathcal X\backslash\{\bold s\}} \Gamma_\text{neg}(\eta,\bold s).
\end{align}
For any $\eta\in\mathcal X\backslash(\mathcal X^s_\text{neg}\cup\{\bold 1\})$ we have that
\begin{align}
\Gamma_\text{neg}(\eta,\bold s)&=\Gamma_\text{neg}(\eta,\mathcal X^s_\text{neg})=\Phi_\text{neg}(\eta,\mathcal X^s_\text{neg})-H_\text{neg}(\eta)\notag\\&\le\Phi_\text{neg}(\bold 1,\mathcal X^s_\text{neg})-H_\text{neg}(\bold 1)=\Gamma_\text{neg}(\bold 1,\mathcal X^s_\text{neg}),\notag
\end{align}
where the inequality follows by the fact that $\bold 1$ is the unique metastable configuration and this means that starting from $\eta\in\mathcal X\backslash\mathcal X^s_\text{neg}$ there are not initial cycles $\mathcal C^\eta_{\{\bold s\}}(\Gamma_\text{neg}(\eta,\bold s))$ deeper than $\mathcal C_{\{\bold s\}}^\bold 1(\Gamma^m_\text{neg})$. Note that this fact holds since we are in the metastability scenario as in the \cite[Subsection 3.5, Example 1]{nardi2016hitting}. Thus, using \eqref{comparisongammaneg}, since for any $\bold r\in\mathcal X^s_\text{neg}\backslash\{\bold s\}$ we have  $\Gamma_\text{neg}(\bold r,\bold s)=\Gamma_\text{neg}(\mathcal X^s_\text{neg}\backslash\{\bold s\},\bold s)$, we conclude that
\begin{align}
\max_{\eta\in\mathcal X\backslash\{\bold s\}}\hspace{-4pt}\Gamma_\text{neg}(\eta,\bold s)\hspace{-1.5pt}=\hspace{-1.5pt}\max\{\max_{\eta\in\mathcal X\backslash(\mathcal X^s_\text{neg}\backslash\{\bold s\})}\hspace{-8pt}\Gamma_\text{neg}(\eta,\bold s),\hspace{-4pt}\max_{\eta\in\mathcal X^s_\text{neg}\backslash\{\bold s\}}\hspace{-4pt}\Gamma_\text{neg}(\eta,\bold s)\}\hspace{-2pt}=\hspace{-2pt}\Gamma_\text{neg}(\bold r,\bold s).\notag
\end{align}
$\qed$

\section{Minimal gates and tube of typical trajectories}\label{secgatestube}
In this section we give a geometrical characterization of the critical configurations and the tube of typical paths for both metastable transitions $\mathbf 1\to\mathcal X^s_\text{neg}$ and $\mathbf 1\to\mathbf s$ for any fixed $\mathbf s\in\mathcal X^s_\text{neg}$.
\subsection{Identification of critical configurations for the transition from the metastable configuration to the set of stable states}\label{secmingatesneg}
This subsection is devoted to a more accurate study of the energy landscape $(\mathcal X,H_{\text{neg}},Q)$. From a technical point of view, the proofs are a generalization of the corresponding results for the Blume Capel model \cite[Section 6]{cirillo2013relaxation}.

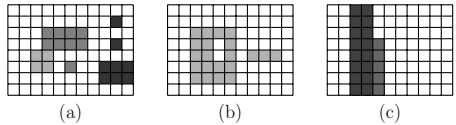
\begin{figure}[h!]
\centering
\begin{tikzpicture}[scale=0.5,transform shape]
\fill[black!30!white] (0,0.6)rectangle(0.6,1.2);
\fill[black!50!white] (0.3,1.2)rectangle(1.5,1.8)(0.9,0.6)rectangle(1.2,0.9);
\fill[black!80!white] (1.8,0.3)rectangle(2.7,0.9)(2.1,1.2)rectangle(2.4,1.5)(2.1,1.8)rectangle(2.4,2.1);
\draw[step=0.3cm,color=black] (-0.6,0) grid (2.7,2.4);
\draw (1.05,-0.1) node[below] {\Large(a)};

\fill[black!30!white](4.2,0.3)rectangle(5.4,1.8);
\fill[white] (4.8,0.9)rectangle(5.1,1.5);
\fill[black!30!white] (5.7,0.9)rectangle(6.6,1.2);
\draw[step=0.3cm,color=black] (3.6,0) grid (6.9,2.4);
\draw (5.25,-0.1) node[below] {\Large(b)};

\fill[black!75!white](8.4,0)rectangle(9,2.4);
\fill[black!60!white](9.3,0)rectangle(9,1.5);
\draw[step=0.3cm,color=black] (7.8,0) grid (11.1,2.4);

\draw (9.45,-0.1) node[below] {\Large(c)};
\end{tikzpicture}
\caption{\label{figureDnegexample} Illustration of three examples of $\sigma\in{\mathscr D}_{\text{neg}}$ when $\ell^*=5$. In (a) the $\ell^*(\ell^*-1)+1=21$ spins different from $1$ have not all the same spin value and they belong to more clusters. In (b) we consider the same number of spins with value $s\neq 1$ that belong to two different clusters.
In (c) we consider the same number of spins different from $1$ that are different between each other and that belong to two adjacent clusters.}
\end{figure}\FloatBarrier

Let $\mathscr D_{\text{neg}}\subset\mathcal X$ be the set
\begin{align}\label{setDN}
\mathscr D_{\text{neg}}:=\{\sigma\in\mathcal X: N_1(\sigma)=|\Lambda|-[\ell^*(\ell^*-1)+1]\}.
\end{align}
Furthermore, let $\mathscr D_{\text{neg}}^+$ and $\mathscr D_{\text{neg}}^-$ be the sets
\begin{align}\label{setD+N}
\mathscr D_{\text{neg}}^+:=\{\sigma\in\mathcal X: N_1(\sigma)>|\Lambda|-[\ell^*(\ell^*-1)+1]\},
\end{align} 
\begin{align}\label{setD-N}
\mathscr D_{\text{neg}}^-:=\{\sigma\in\mathcal X: N_1(\sigma)<|\Lambda|-[\ell^*(\ell^*-1)+1]\}.
\end{align}
Note that $\bold 1\in\mathscr D_{\text{neg}}^+$. For any $\sigma\in\mathscr D_{\text{neg}}$, we remark that $\sigma$ has $\ell^*(\ell^*-1)+1$ spins different from $1$ and they may have all the same spin value and may belong to one or more clusters, see Figure \ref{figureDnegexample}.

\noindent A \textit{two dimensional polyomino} on $\mathbb{Z}^2$ is a finite union of unit squares. The area of a polyomino is the number of its unit squares, while its perimeter is the cardinality of its boundary, namely, the number of interfaces on $\mathbb Z^2$ between the sites inside the polyomino and those outside. The polyominoes with minimal perimeter among those with the same area are said to be \textit{minimal polyominoes}.

\begin{lemma}\label{bottomDN}
If the external magnetic field is negative, then the minimum of the energy in $\mathscr D_{\emph{neg}}$ is achieved by those configurations in which all the spins are equal to $1$ except those, which have the same value $t\neq 1$, in a unique cluster of perimeter $4\ell^*$. 
More precisely,
\begin{align}\label{charbottomneg}
\mathscr F(\mathscr D_{\emph{neg}})=\bigcup_{t=2}^q \mathscr D^t_{\emph{neg}},
\end{align}
where 
\begin{align}\label{Dtnegdef}
\mathscr D^t_{\emph{neg}}:=&\{\sigma\in\mathscr D_{\emph{neg}}: \sigma\ \text{has all spins $1$ except those in a unique cluster $C^t(\sigma)$  }\notag \\ &\text{of spins  $t$ of perimeter}\ 4\ell^*\}.
\end{align}
Moreover, \begin{align}\label{HbottomF}
H_{\text{\emph{neg}}}(\mathscr F(\mathscr D_{\emph{neg}}))=H_{\emph{neg}}(\bold 1)+\Gamma_{\emph{neg}}(\bold 1,\mathcal X^s_{\emph{neg}})=\Phi_{\emph{neg}}(\bold 1,\mathcal X^s_{\emph{neg}}).
\end{align} 
\end{lemma}
\textit{Proof.} Since the presence of disagreeing edges increases the energy, in the configurations in $\mathscr F(\mathscr D_{\text{neg}})$, all $\ell^*(\ell^*-1)+1$ spins different from $1$ are equal to $t$ (say) and belong to a unique cluster $C^t(\sigma)$. As we have illustrated in the second part of the proof of Proposition \ref{lowerboundneg}, the minimal perimeter of a polyomino of area $\ell^*(\ell^*-1)+1$ is $4\ell^*$. Thus, \eqref{charbottomneg} is verified and we get that $\mathcal W_{\text{neg}}(\bold 1,\mathcal X^s_{\text{neg}})\subset\mathscr F(\mathscr D_{\text{neg}})$. Hence, $H_{\text{neg}}(\mathcal W_{\text{neg}}(\bold 1,\mathcal X^s_{\text{neg}}))=H_{\text{neg}}(\mathscr F(\mathscr D_{\text{neg}}))$ and, since for any $\eta\in\mathcal W_{\text{neg}}(\bold 1,\mathcal X^s_{\text{neg}})$
\begin{align}
H_{\text{neg}}(\eta)-H_{\text{neg}}(\bold 1)=4\ell^*-h(\ell^*(\ell^*-1)+1)=\Gamma_{\text{neg}}(\bold 1,\mathcal X^s_{\text{neg}}),
\end{align} 
\eqref{HbottomF} is satisfied.
$\qed$

In the next corollary we prove that $\mathscr F(\mathscr D_{\text{neg}})$ is a gate for the transition $\bold 1\to\mathcal X^s_\text{neg}$.
\begin{cor}\label{lemmarefpath}
If the external magnetic field is negative, then for any $\omega\in\Omega_{\bold 1,\mathcal X^s_{\emph{neg}}}^{opt}$, $\omega\cap\mathscr F(\mathscr D_{\emph{neg}})\neq\varnothing$. In other words, $\mathscr F(\mathscr D_{\emph{neg}})$ is a gate for the transition from $\bold 1$ to $\mathcal X^s_\emph{neg}$.
\end{cor}
\textit{Proof.} For any path $\omega\in\Omega_{\bold 1,\mathcal X^s_{\text{neg}}}$, $\omega=(\omega_0,\dots,\omega_n)$, there exists $i\in\{0,\dots,n\}$ such that $\omega_i\in\mathscr D_{\text{neg}}$. Indeed, given $N(\sigma):=\sum_{t=2}^q N_t(\sigma)$, any path has to pass through the set $\mathcal V_k:=\{\sigma\in\mathcal X:\ N(\sigma)=k\}$, for any $k=0,\dots,|V|$, at least once and $\mathcal V_{\ell^*(\ell^*-1)+1}\equiv\mathscr D_{\text{neg}}$. Since from \eqref{HbottomF} we get that the energy value of any configuration belonging to the bottom of $\mathscr D_{\text{neg}}$ is equal to the min-max reached by any optimal path from $\bold 1$ to $\mathcal X^s_{\text{neg}}$, we get that $\omega_i\in\mathscr F(\mathscr D_{\text{neg}})$. $\qed$

\begin{figure}
\centering
\begin{tikzpicture}[scale=0.5,transform shape]
\fill[black!30!white] (0.6,0.6)rectangle(1.2,1.2)(0.9,0.3)rectangle(1.5,1.8)(1.5,0.3)rectangle(1.8,2.1)(1.8,1.2)rectangle(2.1,2.1);
\draw[step=0.3cm,color=black] (0,0) grid (3.3,2.4);
\draw[thick,dotted,white] (0.6,0.3)rectangle(2.1,2.1);
\draw (1.65,-0.1) node[below] {\Large(a)};
\end{tikzpicture}\ \ \
\begin{tikzpicture}[scale=0.5,transform shape]
\fill[black!30!white](4.2,0.3)rectangle(5.4,1.8)(4.5,1.8)rectangle(4.8,2.1);
\draw[step=0.3cm,color=black] (3.6,0) grid (6.9,2.4);
\draw[thick,dotted,white] (4.2,0.3)rectangle(5.4,2.1);
\draw (5.25,-0.1) node[below] {\Large(b)};
\end{tikzpicture}\ \ \
\begin{tikzpicture}[scale=0.5,transform shape]
\fill[black!30!white](7.8,0.6)rectangle(9,2.1)(9,1.2)rectangle(9.3,1.5);
\draw[step=0.3cm,color=black] (7.2,0) grid (10.5,2.4);
\draw[thick,dotted,white] (7.8,0.6)rectangle(9.3,2.1);
\draw (8.85,-0.1) node[below] {\Large(c)};
\end{tikzpicture}\ \ \
\begin{tikzpicture}[scale=0.5,transform shape]
\fill[black!30!white](12,0.6)rectangle(12.9,2.1)(11.7,0.6)rectangle(12,1.5)(12.9,1.2)rectangle(13.2,1.5)(12,0.3)rectangle(12.6,0.6);
\draw[step=0.3cm,color=black] (10.8,0) grid (14.1,2.4);
\draw[thick,dotted,white] (11.7,0.3)rectangle(13.2,2.1);
\draw (12.45,-0.1) node[below] {\Large(d)};
\end{tikzpicture}
\caption{\label{figureexample} Examples of $\sigma\in\tilde{\mathscr D}^t_{\text{neg}}$ in (a) and of $\sigma\in\hat{\mathscr D}^t_{\text{neg}}$ in (b) and (c) when $\ell^*=5$. We associate the color gray to the spin $t$, the color white to the spin $1$. The dashed rectangle represents the smallest surrounding rectangle of $C^t(\sigma)$.  Figure (d) is an example of configuration that does not belong to $\hat{\mathscr D}^t_{\text{neg}}$.}
\end{figure}
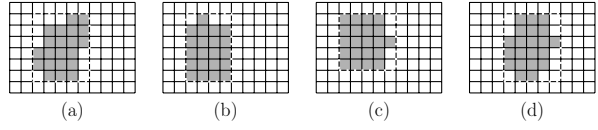

In the next proposition, we show that $\mathcal W_{\text{neg}}(\bold 1,\mathcal X^s_{\text{neg}})$ is a gate for the transition from $\bold 1$ to $\mathcal X^s_{\text{neg}}$. We define $R_{\ell_1\times\ell_2}$ be the set of the rectangles in $\mathbb R^2$ with  sides of length $\ell_1$ and $\ell_2$. 
\begin{proposition}[Gate for the transition from the metastable state to the stable set]\label{propgateneg}
If the external magnetic field is negative, then any path $\omega\in\Omega_{\bold 1,\mathcal X^s_{\emph{neg}}}^{opt}$ visits $\mathcal W_{\emph{neg}}(\bold 1,\mathcal X^s_{\emph{neg}})$. Hence,  $\mathcal W_{\emph{neg}}(\bold 1,\mathcal X^s_{\emph{neg}})$ is a gate for the transition from $\bold 1$ to $\mathcal X^s_{\emph{neg}}$.
\end{proposition}
\textit{Proof.} For any $t\neq 1$, let $\tilde{\mathscr D}^t_{\text{neg}}$ be the set of configurations $\sigma\in\mathscr D^t_{\text{neg}}$ such that the boundary of $C^t(\sigma)$ intersects each side of the boundary of its smallest surrounding rectangle $R(C^t(\sigma))$ on a set of the dual lattice $\mathbb{Z}^2+(1/2,1/2)$ made by at least two consecutive unit segments, see Figure \ref{figureexample}(a). Furthermore, let $\hat{\mathscr D}^t_{\text{neg}}$ be the set of configurations $\sigma\in\mathscr D^t_{\text{neg}}$ such that the boundary of the polyomino $C^t(\sigma)$ intersects at least one side of the boundary of $R(C^t(\sigma))$ in a single unit segment, see for instance Figure \ref{figureexample}(b) and (c).
Hence, $\mathscr F(\mathscr D_{\text{neg}})=\tilde{\mathscr D}_{\text{neg}}\cup\hat{\mathscr D}_{\text{neg}}$, where $\tilde{\mathscr D}_{\text{neg}}:=\bigcup_{t=2}^q \tilde{\mathscr D}^t_{\text{neg}}$ and $\hat{\mathscr D}_{\text{neg}}:=\bigcup_{t=2}^q \hat{\mathscr D}^t_{\text{neg}}$. The proof proceeds in five steps.

\textbf{Step 1}.
Our first aim is to prove that 
\begin{align}\label{Dneghatunion}
\hat{\mathscr D}_{\text{neg}}=\mathcal W_{\text{neg}}(\bold 1,\mathcal X^s_{\text{neg}})\cup\mathcal{W}_{\text{neg}}'(\bold 1,\mathcal X^s_{\text{neg}}).
\end{align}
From \eqref{gate1neg}, we have $\mathcal W_{\text{neg}}(\bold 1,\mathcal X^s_{\text{neg}})\cup\mathcal{W}_{\text{neg}}'(\bold 1,\mathcal X^s_{\text{neg}})\subseteq\hat{\mathscr D}_{\text{neg}}$. Thus, we are left to prove the reverse inclusion $\sigma\in\hat{\mathscr D}_{\text{neg}}$. The boundary of the cluster $C^t(\sigma)$ could intersect the other three sides of the boundary of $R(C^t(\sigma))$ in proper subsets of each side, see Figure \ref{figureexample}(d). Assume $R(C^t(\sigma))\in R_{(\ell^*+a)\times(\ell^*+b)}$ for some $a,b \in\mathbb{Z}$. 
$C^t(\sigma)$ is a minimal polyomino and so it is also convex and monotone by \cite[Lemma 6.16]{cirillo2013relaxation}. Hence, the perimeter of $C^t(\sigma)$ is equal to the perimeter of $R(C^t(\sigma))$, which implies $4\ell^*=4\ell^*+2(a+b)$, and so $a=-b$. Now, let $\hat C^t(\sigma)$ be the polyomino obtained by removing the unit protuberance from $C^t(\sigma)$ and let $\hat R$ be its smallest surrounding rectangle.
,If $C^t(\sigma)$ has the unit protuberance adjacent to a side of length $\ell^*+a$, then $\hat R$ is a rectangle $(\ell^*+a)\times(\ell^*-a-1)$. Since the area of $\hat R$ must be larger than or equal to the area of $\hat C^t(\sigma)$, we have
\[
\text{Area}(\hat R)\hspace{-2pt}=\hspace{-2pt}(\ell^*\hspace{-2pt}+\hspace{-2pt}a)(\ell^*\hspace{-2pt}-\hspace{-2pt}a\hspace{-2pt}-\hspace{-2pt}1)\hspace{-2pt}\ge\text{Area}(\hat C^t(\sigma))=\ell^*(\ell^*-1)\hspace{-4pt}\iff\hspace{-4pt}-a^2-a\ge 0.
\]
Since $a\in\mathbb{Z}$, $-a^2-a\ge 0$ is satisfied only if either $a=0$ or $a=-1$. 
On the other hand, if the unit protuberance of $C^t(\sigma)$ is adjacent to a side of length $\ell^*-a$, the same argument gives
\[
\text{Area}(\hat R)\hspace{-2pt}=\hspace{-2pt}(\ell^*\hspace{-2pt}+\hspace{-2pt}a\hspace{-2pt}-\hspace{-2pt}1)(\ell^*\hspace{-2pt}-\hspace{-2pt}a)\hspace{-2pt}\ge\text{Area}(\hat C^t(\sigma))=\ell^*(\ell^*\hspace{-2pt}-\hspace{-2pt}1) \iff -a^2+a\ge 0.
\]
Again this is satisfied only if either $a=0$ or $a=1$. In both cases we get that $\hat R\in R_{\ell^*\times(\ell^*-1)}$. Thus, if the protuberance is attached to one of the longest sides of $\hat R$, then $\sigma\in\mathcal W_{\text{neg}}(\bold 1,\mathcal X^s_{\text{neg}})$, otherwise $\sigma\in\mathcal W'_{\text{neg}}(\bold 1,\mathcal X^s_{\text{neg}})$. Then, \eqref{Dneghatunion} is verified.

\textbf{Step 2.}  For any $\omega=(\omega_0,\dots,\omega_n)\in\Omega_{\bold 1,\mathcal X^s_{\text{neg}}}^{opt}$ and any $t\in\{2,\dots,q\}$, let
\begin{align}
f_t(\omega):=\{k\in\mathbb{N}: \omega_k\in\mathscr F(\mathscr D_{\text{neg}}),\ N_1(\omega_{k-1})=|\Lambda|-\ell^*(\ell^*-1),\ N_t(\omega_{k-1})\hspace{-2pt}=\hspace{-2pt}\ell^*(\ell^*-1)\}.
\end{align}
We claim that the set $f(\omega):=\bigcup_{t=2}^q f_t(\omega)$ is not empty. Let $\omega=(\omega_0,\dots,\omega_n)\in\Omega_{\bold 1,\mathcal X^s_{\text{neg}}}^{opt}$ and let $\bar k\le n$ be the smallest integer such that $(\omega_{\bar k},\dots,\omega_n)\cap\mathscr D_\text{neg}^+=\varnothing$. Since $\omega_{\bar k-1}$ is the last configuration in $\mathscr D_\text{neg}^+$ along $\omega$, it follows that $\omega_{\bar k}\in\mathscr D_\text{neg}$ and, by the proof of Corollary \ref{lemmarefpath} we have that $\omega_{\bar k}\in\mathscr F(\mathscr D_{\text{neg}})$. Thus, there exists $t\neq 1$ such that $\omega_{\bar k}\in\mathscr D^t_\text{neg}$. Furthermore, $N_1(\omega_{\bar k-1})=|\Lambda|-\ell^*(\ell^*-1)$ and $\omega_{\bar k}$ is obtained from $\omega_{\bar k-1}$ by flipping a spin $1$ to $s\neq 1$. Note that $N_1(\omega_{\bar k-1})=|\Lambda|-\ell^*(\ell^*-1)$ implies that $N_t(\omega_{\bar k-1})\le\ell^*(\ell^*-1)$. Since by Lemma \ref{bottomDN} we have that $N_t(\omega_{\bar k})=\ell^*(\ell^*-1)+1$, we conclude that $N_t(\omega_{\bar k-1})=\ell^*(\ell^*-1)$ since in a single spin flip the number of spins $t$ changes by at most one. Thus, $N_t(\omega_{\bar k-1})=\ell^*(\ell^*-1)$ and $\bar k\in f(\omega)$. 

\textbf{Step 3}. We claim that for any path $\omega\in\Omega_{\bold 1,\mathcal X^s_{\text{neg}}}^{opt}$ one has $\omega_i\in\hat{\mathscr D}_{\text{neg}}$  for all $i\in f(\omega)$. We argue by contradiction. Assume that there exists $i\in f(\omega)$ such that $\omega_i\notin\hat{\mathscr D}_{\text{neg}}$ and $\omega_i\in\tilde{\mathscr D}_{\text{neg}}$. Since $i\in f(\omega)$, there exists $t\neq 1$, such that $i\in f_t(\omega)$. Furthermore, $\omega_{i-1}$ is obtained from $\omega_i$ by flipping a spin $t$ from $t$ to $1$. In view of the definition of $\tilde{\mathscr D}_{\text{neg}}$, every spin equal to $t\neq1$ has at least two nearest neighbors with spin $t$. Hence,
\begin{align}\label{weget}
H_{\text{neg}}(\omega_{i-1})-H_{\text{neg}}(\omega_i)\ge(2-2)+h=h>0.
\end{align} 
From \eqref{weget} we get a contradiction since 
\begin{align}
\Phi_\omega^\text{neg}\ge H_{\text{neg}}(\omega_{i-1})>H_{\text{neg}}(\omega_i)=H_{\text{neg}}(\bold 1)+\Gamma_{\text{neg}}(\bold 1,\mathcal X^s_\text{neg})=\Phi_{\text{neg}}(\bold 1,\mathcal X^s_{\text{neg}}),\notag
\end{align}
where the first identity follows from Lemma \ref{bottomDN}. Then the claim is proved. 

\textbf{Step 4}. 
Now we claim that for any path $\omega\in\Omega_{\bold 1,\mathcal X^s_{\text{neg}}}^{opt}$, $\omega_i\in\mathscr F(\mathscr D_{\text{neg}})$ implies $\omega_{i-1}, \omega_{i+1}\notin\mathscr D_{\text{neg}}.$
Using Corollary \ref{lemmarefpath}, there exists a positive integer $i$ such that $\omega_i\in\mathscr F(\mathscr D_{\text{neg}})$. Thus, there exists $t\neq 1$ such that $\omega_i\in\mathscr D^t_{\text{neg}}$. Assume by contradiction that $\omega_{i+1}\in\mathscr D_{\text{neg}}$. Then $\omega_{i+1}$ must be obtained by $\omega_i$ by flipping a spin $t$ to $s\neq t$, since $N_1(\omega_{i})= N_1(\omega_{i+1})$. In particular, this spin-update increases the energy and so, using Lemma \ref{bottomDN}, we obtain
$\Phi_\omega^{\text{neg}}\ge H_{\text{neg}}(\omega_{i+1})>H_{\text{neg}}(\omega_i)=H_{\text{neg}}(\bold 1)+\Gamma_{\text{neg}}(\bold 1,\mathcal X^s_{\text{neg}})=\Phi_{\text{neg}}(\bold 1,\mathcal X^s_{\text{neg}}),$
which is a contradiction. Hence $\omega_{i+1}\notin\mathscr D_{\text{neg}}$ and similarly we may also prove that $\omega_{i-1}\notin\mathscr D_{\text{neg}}$.

\textbf{Step 5}. 
Our final aim is to show that for any path $\omega\in\Omega_{\bold 1,\mathcal X^s_{\text{neg}}}^{opt}$, we have that $\omega\cap\mathcal W_{\text{neg}}(\bold 1,\mathcal X^s_{\text{neg}})\neq\varnothing$. 
Given a path $\omega=(\omega_0,\dots,\omega_n)\in\Omega_{\bold 1,\mathcal X^s_{\text{neg}}}^{opt}$, assume by contradiction that $\omega\cap\mathcal W_{\text{neg}}(\bold 1,\mathcal X^s_{\text{neg}})=\varnothing$. From step 4 we know that along $\omega$ the configurations which belong to $\mathscr F(\mathscr D_{\text{neg}})$ are not consecutive and they are separated by a subpath which belongs either to $\mathscr D_{\text{neg}}^+$ or to $\mathscr D_{\text{neg}}^-$. Let $j\in\{1,\dots,n\}$ be the smallest integer such that $\omega_j\in\mathscr F(\mathscr D_{\text{neg}})$ and such that $(\omega_j,\dots,\omega_n)\cap\mathscr D_{\text{neg}}^+=\varnothing$.  In particular, $j\in f(\omega)$ since $j$ plays the same role of $\bar k$ in Step 2.  Note that using \eqref{Dneghatunion}, Step 2 and the assumption $\omega\cap\mathcal W_{\text{neg}}(\bold 1,\mathcal X^s_{\text{neg}})=\varnothing$, we have $\omega_j\in\mathcal W'_{\text{neg}}(\bold 1,\mathcal X^s_{\text{neg}})$. Furthermore, by \eqref{HbottomF} the energy along the path from $\omega_j\in\mathscr F(\mathscr D_{\text{neg}})$ to $\omega_n$ decreases. Let $t\neq 1$ be such that $\omega_j\in\mathscr D^t_\text{neg}$. Then the only moves that decrease the energy are

(i ) flipping the spin in the unit protuberance from $t$ to $1$,

(ii) flipping a spin $1$ with two nearest neighbors with spin $t$ from $1$ to $t$.

\noindent Since $\omega_{j+1}\notin\mathscr D_{\text{neg}}^+$, (i) is not feasible. Hence, necessarily $H_{\text{neg}}(\omega_{j+1})=H_{\text{neg}}(\bold 1)+\Gamma_{\text{neg}}(\bold 1,\mathcal X^s_{\text{neg}})-h$ and starting from $\omega_{j+1}$ we consider a spin-update that either decreases the energy or increases the energy of at most $h$. Hence the only feasible moves are

(iii) flipping a spin $1$, with two nearest neighbors with spin $t$, from $1$  to $t$,

(iv) flipping a spin $t$, with two nearest neighbors with spin $1$, from $t$ to $1$.

\noindent Note that by (iii) and (iv), the process reaches a configuration $\sigma$ with all spins equal to $1$ except those, which are $t$, in a polyomino $C^t(\sigma)$ that is convex and such that $R(C^t(\sigma))=R_{(\ell^*+1)\times(\ell^*-1)}$. Note that we cannot iterate move (iv) since otherwise we would find a configuration that does not belong to $\mathscr D_\text{neg}$. On the other hand, applying once (iv) and iteratively (iii), until we fill the rectangle $R_{(\ell^*+1)\times(\ell^*-1)}$ with spins $t$, we find a set of configurations in which the one with the smallest energy is $\sigma$ such that $C^t(\sigma)\equiv R(C^t(\sigma))$. Starting from any configuration of this set, the smallest energy increase is $2-h$ and it is achieved by flipping from $1$ to $t$ a spin $1$ with three nearest neighbors with spin $1$ and a neighbor of spin $t$ inside $C^t(\sigma)$. It follows that
\begin{align}\label{absurdfinally}
\Phi_\omega^\text{neg}-H_{\text{neg}}(\bold 1)&\ge 4\ell^*-h(\ell^*+1)(\ell^*-1)+2-h
>\Gamma_{\text{neg}}(\bold 1,\mathcal X^s_{\text{neg}}),
\end{align}
where the last inequality holds because $2>h(\ell^*-1)$ since $0<h<1$ and $\ell^*:=\left\lceil \frac{2}{h}\right\rceil$, see Assumption \ref{remarkconditionneg}.
Since in \eqref{absurdfinally} we obtained a contradiction, we conclude that any path $\omega\in\Omega_{\bold 1,\mathcal X^s_{\text{neg}}}^{opt}$ must visit $\mathcal W_{\text{neg}}(\bold 1,\mathcal X^s_{\text{neg}})$.
$\qed$

\subsection{Minimal gates: proof of the main results}\label{proofgates}
We are now able to prove Theorems \ref{teogatenegset} and \ref{teogatenegtarget}.

\textit{Proof of Theorem \ref{teogatenegset}}. By Proposition \ref{propgateneg} we get that $\mathcal W_{\text{neg}}(\bold 1,\mathcal X^s_{\text{neg}})$ is a gate for the transition from the metastable state $\bold 1$ to $\mathcal X^s_{\text{neg}}$.  In order to show that $\mathcal W_{\text{neg}}(\bold 1,\mathcal X^s_{\text{neg}})$ is a minimal gate, we exploit \cite[Theorem 5.1]{manzo2004essential} and we show that any $\eta\in\mathcal W_{\text{neg}}(\bold 1,\mathcal X^s_{\text{neg}})$ is an essential saddle. In order to do this, in view of the definition of an essential saddle given in Subsection \ref{mainresgates}, for any $\eta\in\mathcal W_{\text{neg}}(\bold 1,\mathcal X^s_{\text{neg}})$ we construct an optimal path from $\bold 1$ to $\mathcal X^s_{\text{neg}}$ passing through $\eta$ and reaching its maximum energy only there. Since $\eta\in\mathcal W_{\text{neg}}(\bold 1,\mathcal X^s_{\text{neg}})$, there exists $s\neq 1$ such that $\eta\in\bar B_{\ell^*-1,\ell^*}(s,1)$ and the optimal path above is defined by modifying
the reference path $\hat\omega$ of Definition \ref{refpathmioneg} in a such a way that $\hat\omega_{\ell^*(\ell^*-1)+1}=\eta$ in which $C^s(\eta)$ is a quasi-square $\ell^*\times(\ell^*-1)$ with a unit protuberance. It follows that $\hat\omega\cap\mathcal W_{\text{neg}}(\bold 1,\mathcal X^s_{\text{neg}})=\{\eta\}$ and $\text{arg max}_{\hat\omega}H_\text{neg}=\{\eta\}$ by the proof of Lemma \ref{lemmatwoprefactor}.
To conclude, we prove that $\mathcal W_{\text{neg}}(\bold 1,\mathcal X^s_{\text{neg}})$ is the only minimal gate. Note that the above reference path $\hat\omega$ reaches the energy $\Phi_\text{neg}(\bold 1,\mathcal X^s_{\text{neg}})$ only in $\mathcal W_{\text{neg}}(\bold 1,\mathcal X^s_{\text{neg}})$. It follows that, for any $\eta_1\in\mathcal W_{\text{neg}}(\bold 1,\mathcal X^s_{\text{neg}})$, the set $\mathcal W_{\text{neg}}(\bold 1,\mathcal X^s_{\text{neg}})\backslash\{\eta_1\}$ is not a gate for the transition $\bold 1\to\mathcal X^s_{\text{neg}}$. Indeed, from the above discussion we get that there exists an optimal path $\hat\omega$ such that $\hat\omega\cap\mathcal W_{\text{neg}}(\bold 1,\mathcal X^s_{\text{neg}})\backslash\{\eta_1\}=\varnothing$. Note that the uniqueness of the minimal gate follows by the condition $2/h\notin\mathbb N$, see Assumption \ref{remarkconditionneg}.
$\qed$

\textit{Proof of Theorem \ref{teogatenegtarget}}.
For any $\bold s\in\mathcal X^s_{\text{neg}}$, the min-max energy value that is reached by any path $\omega:\bold 1\to\bold s$ is $\Phi_{\text{neg}}(\bold 1,\bold s)\equiv\Phi_{\text{neg}}(\bold 1,\mathcal X^s_{\text{neg}})$. Furthermore, Theorem \ref{theoremcomparisonneg} implies that when a path $\omega:\bold 1\to\bold s$ visits some $\bold r\in\mathcal X^s_{\text{neg}}\backslash\{\bold s\}$, the min-max energy value that the path reaches is still $\Phi_{\text{neg}}(\bold 1,\mathcal X^s_{\text{neg}})$. Indeed, for instance in the case in which the path $\omega$ may be decomposed in two paths $\omega_1:\bold1\to\bold r$ and $\omega_2:\bold r\to\bold s$, we have $\Phi_{\omega}^\text{neg}=\max\{\Phi_{\omega_1}^\text{neg},\Phi_{\omega_2}^\text{neg}\}=\Phi_{\text{neg}}(\bold 1,\mathcal X^s_{\text{neg}})$ where we used \eqref{comparisonalignneg}.
Hence, the saddles visited by the process are only the ones crossed during the transition between $\bold 1$ and the first stable state. This fact, together with Theorem \ref{teogatenegset}, allows us to state that the set $\mathcal W_{\text{neg}}(\bold 1,\mathcal X^s_{\text{neg}})$ is the unique minimal gate for the transition from $\bold 1$ to $\bold s$, for any fixed $\bold s\in\mathcal X^s_\text{neg}$. Thus, \eqref{mingatenegtarget} is satisfied. $\qed$

\subsection{Tube of typical trajectories: proof of the main results}\label{prooftube}
In order to give the proofs of Theorems \ref{teotubesetneg} and \ref{teotubefixneg}, first we prove the following lemmas. 

\begin{lemma}\label{lemmapbsottocritical}
 Let $\mathcal C(\eta)$ and $\mathcal C(\zeta)$ be the non-trivial cycles whose bottom are $\eta\in\bar R_{\ell,\ell-1}(1,s)$ and $\zeta\in\bar R_{\ell,\ell}(1,s)$ with $\ell\le\ell^*-1$ and $s\neq 1$, respectively. Then,  
\begin{align}
\mathcal B(\mathcal C(\eta))&=\bar B^1_{\ell-1,\ell-1}(1,s);\\
\label{lemmatubeone}
\mathcal B(\mathcal C(\zeta))&=\bar B^1_{\ell-1,\ell}(1,s).
\end{align}
\end{lemma}
\textit{Proof.} For any $s\neq 1$, let $\eta_1\in\bar R_{\ell,\ell-1}(1,s)$ with $\ell\le\ell^*$. By Proposition \ref{proplocalminima}, $\eta_1\in\mathscr M^3_\text{neg}$ is a local minimum for the Hamiltonian $H_\text{neg}$. Using \eqref{principalboundary}, our aim is to prove the following
\begin{align}\label{claimonetubeuno}
\bar B^1_{\ell-1,\ell-1}(1,s)=\mathscr F(\partial\mathcal C(\eta_1)).
\end{align}
In $\eta_1$, for any $v\in V$ the corresponding $v$-tile (see before Lemma \ref{lemmastablevertices} for the definition) is of type (a), (b), (d), (e) and (h), see Figure \ref{figmattonelleneg}. Starting from $\eta_1$, by flipping to $1$ (resp.~$s$) the spin $s$ (resp.~$1$) on a vertex whose tile is of type (a), (d) (resp. (b), (e)), the process visits a configuration $\sigma_1$ such that
\begin{align}\label{maxincreases}
H_\text{neg}(\sigma_1)-H_\text{neg}(\eta_1)\ge2-h.
\end{align}
Thus, the smallest energy increase is given by $h$ by flipping to $1$ a spin $s$ on a vertex $v_1$ centered in a tile of type (h). Let $\eta_2:=\eta_1^{v_1,1}\in\bar B_{\ell-1,\ell-1}^{\ell-2}(1,s)$.
 In $\eta_2$, for any $v\in V$ the corresponding $v$-tile is one among those depicted in Figure \ref{figmattonelleneg}(a), (b), (d), (e), (h) and (p) with $r=s$.
 Since $H_\text{neg}(\eta_2)=H_\text{neg}(\eta_1)+h$, the spin flips on a vertex whose tile is of type (a), (b), (d) and (e)
 lead to $H_\text{neg}(\sigma_2)-H_\text{neg}(\eta_1)\ge2$. Thus, as in the previous case, the smallest energy increase is given by flipping to $1$ a spin $s$ on a vertex $v_1$ centered in a tile of type (h). Note that starting from $\eta_2$ the only spin flip which decreases the energy leads to the bottom of $\mathcal C(\eta_1)$, namely in $\eta_1$.

\noindent Let us now note that
\begin{align}\label{etaelleuno}
H_\text{neg}(\eta_{\ell-1})-H_\text{neg}(\eta_1)=h(\ell-2).
\end{align}
Since $\ell\le\ell^*$, comparing \eqref{maxincreases} with \eqref{etaelleuno}, we get that $\eta_{\ell-1}\in\mathscr F(\partial\mathcal C(\eta_1))$, and \eqref{claimonetubeuno} is verified.

\noindent Let us now consider for any $s\neq 1$ the local minimum $\zeta_1\in\bar R_{\ell,\ell}(1,s)\subset\mathscr M^3_\text{neg}$ with $\ell\le\ell^*-1$.
Arguing similarly to the previous case, \eqref{lemmatubeone} may be verified by proving that
%
$\bar B^1_{\ell-1,\ell}(1,s)=\mathscr F(\partial\mathcal C(\zeta_1)).$
$\qed$
\begin{lemma}\label{lemmapbsupercritical}
 Let $\mathcal C(\eta)$ be the non-trivial cycle whose bottom is $\eta\in\bar R_{\ell_1,\ell_2}(1,s)$ with $\min\{\ell_1,\ell_2\}\ge\ell^*$ and $s\neq 1$. Then,  
\begin{align}
\mathcal B(\mathcal C(\eta))&=\bar B^1_{\ell_1,\ell_2}(1,s)\cup\bar B^1_{\ell_2,\ell_1}(1,s).
\end{align}
\end{lemma}
\textit{Proof.} For any $s\neq 1$, let $\eta_1\in\bar R_{\ell_1,\ell_2}(1,s)$ with $\ell^*\le\ell_1\le\ell_2$. By Proposition \ref{proplocalminima}, $\eta_1\in\mathscr M^3_\text{neg}$ is a local minimum for the Hamiltonian $H_\text{neg}$. Using \eqref{principalboundary}, our aim is to prove the following
\begin{align}\label{claimonetubethree}
\bar B^1_{\ell_1,\ell_2}(1,s)\cup\bar B^1_{\ell_2,\ell_1}(1,s)=\mathscr F(\partial\mathcal C(\eta_1)).
\end{align}
In $\eta_1$, for any $v\in V$ the corresponding $v$-tile is of type (a), (b), (d), (e) and (h). Let $v_1\in V$ such that the $v_1$-tile is of type (e) with $r=s$, and let $\eta_2:=\eta_1^{v_1,s}$. Note that if $v_1$ is adjacent to a side of length $\ell_2$, then $\eta_2\in\bar B^1_{\ell_1,\ell_2}(1,s)$, otherwise $\eta_2\in\bar B^1_{\ell_2,\ell_1}(1,s)$. Without loss of generality, let us assume that $\eta_2\in\bar B^1_{\ell_1,\ell_2}(1,s)$. By simple algebraic calculation we obtain that
\begin{align}
H_\text{neg}(\eta_2)-H_\text{neg}(\eta_1)=2-h.
\end{align}
In $\eta_2$ for any $v\in V$ the corresponding $v$-tile is of type (a), (b), (d), (e), (h) and (p) with $t=r=1$. By flipping to $s$ a spin $1$ on a vertex $w$ whose tile is of type (p) with $r=s$ the energy decreases by $h$ and the process enters a cycle different from the previous one that is either the cycle $\bar{\mathcal C}$ whose bottom is a local minimum belonging to $\bar R_{\ell_1+1,\ell_2}(m,1)$, or a trivial cycle for which iterating this procedure the process enters $\bar{\mathcal C}$. Thus, $\bar B^1_{\ell_1,\ell_2}(1,s)\subseteq\partial\mathcal C(\eta_1)$. Similarly we prove that $\bar B^1_{\ell_2,\ell_1}(1,s)\subseteq\partial\mathcal C(\eta_1)$. 

Let us now note that starting from $\eta_1$ the smallest energy increase is $h$, and it is given by flipping to $1$ a spin $s$ on a vertex whose tile is of type (h). Let us consider the uphill path $\omega$ started in $\eta_1$ and constructed by flipping to $1$ all the spins $s$ along a side of the rectangular $\ell_1\times\ell_2$ $s$-cluster, say one of length $\ell_1$. Using the discussion given in the proof of Lemma \ref{lemmapbsottocritical} and the construction of $\omega$, we get that the process intersects $\partial\mathcal C(\eta)$ in a configuration $\sigma$ belonging to $\bar B^1_{\ell_2-1,\ell_1}(1,s)$. By algebraic computations, we obtain the following
\begin{align}\label{etaelleunobis}
H_\text{neg}(\sigma)-H_\text{neg}(\eta_1)=h(\ell_2-1).
\end{align}
Since $\ell_2\ge\ell^*$, it follows that $H_\text{neg}(\sigma)>H_\text{neg}(\eta_2)$.

\noindent Since by flipping to $1$ (resp.~$s$) the vertex centered in a tile of type (a), (d) (resp.~(e)), the energy increase is largest than or equal to $2+h$, it follows that \eqref{claimonetubethree} is satisfied. $\qed$

In order to prove Theorems \ref{teotubesetneg} and \ref{teotubefixneg}, we need some further definitions that are taken from \cite{nardi2016hitting, cirillo2015metastability, olivieri2005large}. Our goal is to give an equivalent definition of the tube that only relies on the energy landscape data. 
We call \textit{cycle-path} a finite sequence $(\mathcal C_1,\dots,\mathcal C_m)$ of trivial or non-trivial cycles $\mathcal C_1,\dots,\mathcal C_m\in\mathscr C(\mathcal X)$, such that $\mathcal C_i\cap\mathcal C_{i+1}=\varnothing\ \ \text{and}\ \ \partial\mathcal C_i\cap\mathcal C_{i+1}\neq\varnothing,\ \text{for every}\ i=1,\dots,m-1.$ 

\noindent A cycle-path  $(\mathcal C_1,\dots,\mathcal C_m)$ is said to be \textit{downhill} (\textit{strictly downhill}) if the cycles $\mathcal C_1,\dots,\mathcal C_m$ are pairwise connected with decreasing height, i.e., when $H(\mathscr F(\partial\mathcal C_i))\ge H(\mathscr F(\partial\mathcal C_{i+1}))$ ($H(\mathscr F(\partial\mathcal C_i))> H(\mathscr F(\partial\mathcal C_{i+1}))$) for any $i=0,\dots,m-1$.

\noindent We denote the set of cycle-paths that lead from $\sigma$ to $\mathcal A$ and consist of maximal cycles  in $\mathcal X\backslash\mathcal A$ as 
\begin{align}
\mathcal {P}_{\sigma,\mathcal A}\hspace{-3pt}:=\hspace{-3pt}\{\text{cycle-path} (\mathcal C_1,...,\mathcal C_m)|\mathcal C_1,...,\mathcal C_m\hspace{-3pt}\in\hspace{-3pt}\mathcal M(\mathcal C_{\mathcal{A}}^+(\sigma)\backslash A),\sigma\hspace{-3pt}\in\hspace{-3pt}\mathcal C_1, \partial\mathcal C_m\hspace{-2pt}\cap\hspace{-2pt}\mathcal A\neq\hspace{-2pt}\varnothing\}.\notag
\end{align}

 Given a non-empty set $\mathcal{A}\subset\mathcal X$ and $\sigma\in\mathcal{X}$, we constructively define a mapping $G: \Omega_{\sigma,A}\to\mathcal P_{\sigma,\mathcal A}$ in the following way. Given $\omega=(\omega_1,\dots,\omega_n)\in\Omega_{\sigma,A}$, we set $m_0=1$, $\mathcal C_1=\mathcal C_{\mathcal A}(\sigma)$ and define recursively $m_i:=\min\{k>m_{i-1}|\ \omega_k\notin\mathcal C_i\}$ and $\mathcal C_{i+1}:=\mathcal C_{\mathcal A}(\omega_{m_i})$. We note that $\omega$ is a finite sequence and $\omega_n\in\mathcal A$, so there exists an index $n(\omega)\in\mathbb N$ such that $\omega_{m_{n(\omega)}}=\omega_n\in\mathcal A$ and there the procedure stops. By  $(\mathcal C_1,\dots,\mathcal C_{m_{n(\omega)}})$ is a cycle-path with $\mathcal C_1,\dots,\mathcal C_{m_{n(\omega)}}\subset\mathcal M(\mathcal X\backslash\mathcal A)$. Moreover, the fact that $ \omega\in\Omega_{\sigma,A}$ implies that $\sigma\in\mathcal C_1$ and that $\partial\mathcal C_{n(\omega)}\cap\mathcal A\neq\varnothing$, hence $G(\omega)\in\mathcal P_{\sigma,\mathcal A}$ and the mapping is well-defined.

\noindent We say that a cycle-path $(\mathcal C_1,\dots,\mathcal C_m)$ is \textit{connected via typical jumps} to $\mathcal A\subset\mathcal X$ or simply $vtj-$\textit{connected} to $\mathcal A$ if
\begin{align}\label{cyclepathvtj}
\mathcal B(\mathcal C_i)\cap\mathcal C_{i+1}\neq\varnothing,\ \ \forall i=1,\dots,m-1,\ \ \text{and}\ \ \mathcal B(\mathcal C_m)\cap\mathcal A\neq\varnothing.\end{align}
Let $J_{\mathcal C,\mathcal A}$ be the collection of all cycle-paths $(\mathcal C_1,\dots,\mathcal C_m)$ that are vtj-connected to $\mathcal A$ and such that $\mathcal C_1=\mathcal C$. Given a non-empty set $\mathcal{A}$ and $\sigma\in\mathcal{X}$, we define $\omega\in\Omega_{\sigma,A}$ as a \textit{typical path} from $\sigma$ to $\mathcal A$ if its corresponding cycle-path $G(\omega)$ is vtj-connected to $\mathcal A$ and we denote by $\Omega_{\sigma,A}^{\text{vtj}}$ the collection of all typical paths from $\sigma$ to $\mathcal A$, i.e., 
\begin{align}\label{defvtjpaths}
\Omega_{\sigma,\mathcal A}^{\text{vtj}}:=\{\omega\in\Omega_{\sigma, \mathcal A}|\ G(\omega)\in J_{\mathcal C_{\mathcal A}(\sigma),\mathcal A}\}.
\end{align}

\noindent Finally, we define the \textit{tube of typical paths} $T_{\mathcal A}(\sigma)$ from $\sigma$ to $\mathcal A$ as the subset of states $\eta\in\mathcal X$ that can be reached from $\sigma$ by means of a typical path which does not enter $\mathcal A$ before visiting $\eta$, i.e.,
\begin{align}
T_{\mathcal A}(\sigma):=\{\eta\in\mathcal X|\ \exists\omega\in\Omega_{\sigma, \mathcal A}^{\text{vtj}}:\ \eta\in\omega\}.
\end{align}
Finally, we define $\mathfrak T_{\mathcal A}(\sigma)$ as the set of all maximal cycles that belong to at least one vtj-connected path from $\mathcal C_{\mathcal A}^\sigma(\Gamma)$ to $\mathcal A$, i.e.,
\begin{align}\label{tubostorto}
\mathfrak T_{\mathcal A}(\sigma)\hspace{-2pt}:=\hspace{-2pt}\{\mathcal C\hspace{-2pt}\in\hspace{-2pt}\mathcal M(\mathcal C_{\mathcal A}^+(\sigma)\backslash\mathcal A)|\exists(\mathcal C_1,\dots,\mathcal C_n)\hspace{-2pt}\in\hspace{-2pt} J_{\mathcal C_{\mathcal A}^\sigma(\Gamma),\mathcal A}, \exists j\in\{1,...,n\}\hspace{-2pt}:\hspace{-2pt}\mathcal C_j\hspace{-2pt}=\hspace{-2pt}\mathcal C\}.
\end{align}
Note that 
\begin{align}\label{reltubodrittoetorto}
\mathfrak T_{\mathcal A}(\sigma)=\mathcal M(T_{\mathcal A}(\sigma)\backslash\mathcal A)
\end{align}
and that the boundary of $T_{\mathcal A}(\sigma)$ consists of states either in $\mathcal A$ or in the non-principal part of the boundary of some $\mathcal C\in\mathfrak T_{\mathcal A}(\sigma)$, i.e, $\partial T_{\mathcal A}(\sigma)\backslash\mathcal A\subseteq \bigcup_{\mathcal C\in\mathfrak T_{\mathcal A}(\sigma)} (\partial\mathcal C\backslash\mathcal B(\mathcal C))=:\partial^{np}\mathfrak T_{\mathcal A}(\sigma).$

\textit{Proof of Theorem \ref{teotubesetneg}}. Following the same approach as \cite[Section 6.7]{olivieri2005large}, we characterize the tube of typical trajectories using the so-called ``standard cascades''. See \cite[Figure 6.3]{olivieri2005large} for an example of a standard cascade. We describe these in terms of the paths that are started in $\mathbf 1$ and are vtj-connected to $\mathcal X^s_\text{neg}$. See \eqref{defvtjpaths} for the formal definition and see \cite[Lemma 3.12]{nardi2016hitting} for an equivalent characterization of these paths. We remark that  any typical path from $\mathbf 1$ to $\mathcal X^s_\text{neg}$ is also an optimal path for the same transition.

In order to give a geometrical description of these typical paths, we proceed similarly to \cite[Section 7.4]{olivieri2005large}, where the authors apply the model-independent results given in Section 6.7 to identify the tube of typical paths in the context of the Ising model. We define a vtj-connected cycle-path that is the concatenation of both trivial and non-trivial cycles. Let $\eta_1$ be a configuration belonging to one of the minimal gates for the transition $\mathbf 1\to\mathcal X^s_\text{neg}$, see Theorem \ref{teogatenegset} . We begin by studying the first descent from $\eta_1$ both to $\mathbf 1$ and to $\mathcal X^s_\text{neg}$. Then, we complete the description of $\mathfrak T_{\mathcal X^s_\text{neg}}(\mathbf 1)$ by joining the time reversal of the first descent from $\eta_1$ to $\mathbf 1$ with the first descent from $\eta_1$ to $\mathcal X^s_\text{neg}$. In view of \eqref{gate1neg} we have that $\eta_1\in\bar B^1_{\ell^*-1,\ell^*}(1,s)$ for some $s\neq 1$, and for the sake of semplicity we describe a vtj-connected path from $\mathbf 1$ to $\mathcal X^s_\text{neg}$ conditioned to hit $\mathcal X^s_\text{neg}$ for the first time in $\mathbf s$.

Let us begin by studying the standard cascades from $\eta_1$ to $\mathbf 1$. Since a spin flip from $s$ to $t\notin\{1,s\}$ implies an increase of the energy value equal to the increase of the number of the disagreeing edges, we consider only the splin-flips from $s$ to $1$ on those vertices belonging to the $s$-cluster.
Thus, starting from $\eta_1$ and given $v_1$ a vertex such that $\eta_1(v_1)=s$, since $H_\text{neg}(\eta_1)=\Phi_\text{neg}(\mathbf 1,\mathcal X^s_\text{neg})$, we get
\begin{align}
H_\text{neg}(\eta_1^{v_1,1})=\Phi_\text{neg}(\mathbf 1,\mathcal X^s_\text{neg})+n_s(v_1)-n_1(v_1)+h.
\end{align}
It follows that the only possibility for the path to be optimal is $n_s(v_1)=1$ and $n_1(v_1)=3$. Thus, along the first descent from $\eta_1$ to $\mathbf 1$ the process visits $\eta_2$ in which all the vertices have spin $1$ except those, which are $s$, in a rectangular cluster $\ell^*\times(\ell^*-1)$, i.e., $\eta_2\in\bar R_{\ell^*-1,\ell^*}(1,s)$. By Proposition \ref{proplocalminima} $\eta_2\in\mathscr M^3_\text{neg}$ is a local minimum, thus according to \eqref{cyclepathvtj} we have to describe its non-trivial cycle and its principal boundary. Starting from $\eta_2$, the next configuration along a typical path is defined by flipping to $1$ a spin $s$ on a vertex $v_2$ on one of the four corners of the rectangular $s$-cluster. Indeed, since $H_\text{neg}(\eta_2)=\Phi_\text{neg}(\mathbf 1,\mathcal X^s_\text{neg})-2+h$, we have
\begin{align}
H_\text{neg}(\eta_2^{v_2,1})
=\Phi_\text{neg}(\mathbf 1,\mathcal X^s_\text{neg})-2+2h+n_s(v_2)-n_1(v_2),
\end{align}
and for the path to be optimal, we must have $n_s(v_2)=2$ and $n_1(v_2)=2$. 
By \eqref{minendiff}, the smallest energy increase for any single step of the dynamics is equal to $h$. Thus, a typical path towards $\mathbf 1$ proceeds by eroding the $\ell^*-2$ unit squares with spin $s$ belonging to a side of length $\ell^*-1$ that are corners of the $s$-cluster and that belong to the same side of $v_2$. Each of the first $\ell^*-3$ spin flips increases the energy by $h$, and these uphill steps are necessary in order to exit from the cycle whose bottom is the local minimum $\eta_2$. After these $\ell^*-3$ steps, the process hits the bottom of the boundary of this cycle in a configuration $\eta_{\ell^*}\in\bar B_{\ell^*-1,\ell^*-1}^1(1,s)$, see Lemma \ref{lemmapbsottocritical}. 
The last spin-update, that flips from $s$ to $1$ the spin $s$ on the unit protuberance of the $s$-cluster, decreases the energy by $2-h$. Thus, the typical path arrives in a local minimum $\eta_{\ell^*+1}\in\bar R_{\ell^*-1,\ell^*-1}(1,s)$, i.e., it enters a new cycle whose bottom is a configuration in which all the vertices have spin $1$, except those, which are $s$, in a square $(\ell^*-1)\times(\ell^*-1)$ $s$-cluster. Summarizing the construction above, we have the following sequence of vtj-connected cycles
\begin{align}
\{\eta_1\},\mathcal C^{\eta_2}_{\mathbf 1}(h(\ell^*-2)),\{\eta_{\ell^*}\},\mathcal C^{\eta_{\ell^*+1}}_{\mathbf 1}(h(\ell^*-2)).
\end{align}
Iterating this argument, we obtain that the first descent from $\eta_1\in\mathcal W_\text{neg}(\mathbf 1,\mathcal X^s_\text{neg})$ to $\mathbf 1$ is characterized by the concatenation of those vtj-connected cycle-subpaths between the cycles whose bottom is the local minima in which all the vertices have spin equal to $1$, except those, which are $s$, in either a quasi-square $(\ell-1)\times\ell$ or a square $(\ell-1)\times(\ell-1)$ for any $\ell=\ell^*,\dots,1$, and whose depth is given by $h(\ell -2)$. More precisely, from a quasi-square to a square, a typical path proceeds by flipping to $1$ those spins $s$ on one of the shortest sides of the $s$-cluster. On the other hand, from a square to a quasi-square, it proceeds by flipping to $1$ those spins $s$ belonging to one of the four sides of the square. 
Thus, a standard cascade from $\eta_1$ to $\mathbf 1$ is characterized by the sequence of those configurations that belong to
\begin{align}\label{firstdescenteeta1}
\bigcup_{\ell=1}^{\ell^*}\biggr[\bigcup_{l=1}^{\ell-1}\hspace{-1.2pt}\bar B^l_{\ell-1,\ell}(1,s)\hspace{-1.2pt}\cup\hspace{-1.2pt}\bar R_{\ell-1,\ell}(1,s)\hspace{-1.2pt}\cup\hspace{-1.2pt}\bigcup_{l=1}^{\ell-2}\hspace{-1.2pt}\bar B^l_{\ell-1,\ell-1}(1,s)\hspace{-1.2pt}\cup\hspace{-1.2pt}\bar R_{\ell-1,\ell-1}(1,s)\biggl]. 
\end{align}

\noindent Let us now consider the first descent from $\eta_1\in\bar B^1_{\ell^*-1,\ell^*}(1,s)$ to $\mathbf s\in\mathcal X^s_\text{neg}$. Since the path is optimal, we only consider flips from $1$ to $s$. Thus, let $w_1$ be a vertex such that $\eta_1(w_1)=1$. Flipping the spin $1$ on the vertex $w_1$, we get
\begin{align}
H_\text{neg}(\eta_1^{w_1,s})=\Phi_\text{neg}(\mathbf 1,\mathcal X^s_\text{neg})+n_1(w_1)-n_s(w_1)-h,
\end{align}
and the only feasible choice is $n_1(w_1)=2$ and $n_s(w_1)=2$. Thus, $\eta_1^{w_1,s}\in\bar B^2_{\ell^*-1,\ell^*}(1,s)$, namely the bar is now of length two. Arguing similarly, we get that along the descent to $\mathbf s$ a typical path proceeds by flipping from $1$ to $s$ the spins $1$ with two nearest-neighbors with spin $s$ and two nearest-neighbors with spin $1$ belonging to the incomplete side of the $s$-cluster. More precisely, it proceeds downhill visiting $\bar\eta_i\in\bar B_{\ell^*-1,\ell^*}^{i}(1,s)$ for any $i=2,\dots,\ell^*-1$ and $\bar\eta_{\ell^*}\in\bar R_{\ell^*,\ell^*}(1,s)$, which is a local minimum by Proposition \ref{proplocalminima}. In order to exit from the cycle whose bottom is $\bar\eta_{\ell^*}$, the process crosses the bottom of its boundary by creating a unit protuberance of spin $s$ adjacent to one of the four edges of the $s$-square, i.e., visits $\{\bar\eta_{\ell^*+1}\}$ where $\bar\eta_{\ell^*+1}\in\bar B_{\ell^*,\ell^*}^1(1,s)$, see Lemma \ref{lemmapbsupercritical}. Starting from $\{\bar\eta_{\ell^*+1}\}$, a typical path towards $\mathbf s$ proceeds by enlarging the protuberance to a bar of length two to $\ell^*-1$, thus it visits $\bar\eta_{\ell^*+i}\in\bar B_{\ell^*,\ell^*}^i(1,s)$ for any $i=2,\dots,\ell^*-1$. Each of these steps decreases the energy by $h$, and eventually the bottom of the cycle is reached, i.e., in the local minimum $\bar\eta_{2\ell^*}\in\bar R_{\ell^*,\ell^*+1}(1,s)$. Then, the process exits from this cycle through the bottom of its boundary by adding a unit protuberance of spin $s$ on \textit{any} one of the four edges of the rectangular $\ell^*\times(\ell^*+1)$ $s$-cluster in $\bar\eta_{2\ell^*}$. Thus, it visits the trivial cycle $\{\bar\eta_{2\ell^*+1}\}$, where $\bar\eta_{2\ell^*+1}\in\bar B^1_{\ell^*,\ell^*+1}(1,s)\cup\bar B^1_{\ell^*+1,\ell^*}(1,s)$. Note that  the resulting standard cascade is different from the one towards $\mathbf 1$. Thus, summarizing the construction above, we have defined the following sequence of vtj-connected cycles
\begin{align}
\{\eta_1\},\mathcal C^{\bar\eta_{\ell^*}}_{\mathbf s}(h(\ell^*-1)),\{\bar\eta_{\ell^*+1}\},\mathcal C^{\bar\eta_{2\ell^*}}_{\mathbf s}(h(\ell^*-1)),\{\bar\eta_{2\ell^*+1}\}. 
\end{align}
Note that if $\bar\eta_{2\ell^*}\in\bar B^1_{\ell^*,\ell^*+1}(1,s)$, then the process enters the cycle whose bottom is a configuration belonging to $\bar R_{\ell^*+1,\ell^*+1}(1,s)$. On the other hand, if $\bar\eta_{2\ell^*}\in\bar B^1_{\ell^*+1,\ell^*}(1,s)$, then the standard cascade enters the cycle whose bottom is a configuration belonging to $\bar R_{\ell^*,\ell^*+2}(1,s)$. In the first case the cycle has depth $h\ell^*$, in the second case the cycle has depth $h(\ell^*-1)$. Iterating this argument, 
we get that the first descent from $\eta_1$ to $\mathbf s$ is characterized by vtj-connected cycle-subpaths from $\bar R_{\ell_1,\ell_2}(1,s)$ to $\bar R_{\ell_1,\ell_2+1}(1,s)$ 
defined as the sequence of those configurations belonging to $\bar B^l_{\ell_1,\ell_2}(1,s)$ for any $l=1,\dots,\ell_2-1$. Eventually, a configuration in which this cluster is either a vertical or a horizontal strip is reached, i.e., it intersects one of the two sets defined in \eqref{setverticalstrip}--\eqref{sethorizontalstrip}. If the descent arrives in $\mathscr S_\text{neg}^v(1,s)$, then it proceeds by enlarging the vertical strip column by column. Otherwise, if it arrives in $\mathscr S_\text{neg}^h(1,s)$, then it enlarges the horizontal strip row by row. In both cases, starting from a configuration with an $s$-strip, i.e., a local minimum in $\mathscr M^2_\text{neg}$ by Proposition \ref{proplocalminima}, the path exits from its cycle by adding a unit protuberance with a spin $s$ adjacent to one of the two vertical (resp. horizontal) edges and increasing the energy by $2-h$. Starting from this trivial cycle, the standard cascade proceeds downhill in a new cycle by filling the column (resp.~row) with spins $s$. More precisely, the standard cascade visits $K-1$ (resp.~$L-1$) configurations such that each of them is defined by the previous one flipping from $1$ to $s$ a spin $1$ with two nearest-neighbors with spin $1$ and two nearest-neighbors with spin $s$. Each of these spin-updates decreases the energy by $h$. The process arrives in this way to the bottom of the cycle, i.e., in a configuration in which the thickness of the $s$-strip has been enlarged by a column (resp.~row). Starting from this state with the new $s$-strip, we repeat the same arguments above until the standard cascade arrives in the trivial cycle of a configuration $\sigma$ with an $s$-strip of thickness $L-2$ (resp.~$K-2$) and with a unit protuberance. Starting from $\{\sigma\}$, the process enters the cycle whose bottom is $\mathbf s$ and it proceeds downhill either by flipping from $1$ to $s$ those spins $1$ with two nearest-neighbors with spin $1$ and two nearest-neighbors with spin $s$, or by flipping to $s$ all the spins $1$ with three nearest-neighbors with spin $s$ and one nearest-neighbor with spin $1$. The last step flips from $1$ to $s$ the last spin $1$ with four nearest-neighbors with spin $s$. 
Thus, the first descent from $\eta_1$ to $\mathcal X^s_\text{neg}$ conditioning to hit this set in $\mathbf s$ is characterized by the sequence of those configurations that belong to 
\begin{align}\label{firstdescenteeta2}
&\bigcup_{\ell_1=\ell^*}^{K-1}\bigcup_{\ell_2=\ell^*}^{K-1} \bar R_{\ell_1,\ell_2}(1,s)\cup\bigcup_{\ell_1=\ell^*}^{K-1}\bigcup_{\ell_2=\ell^*}^{K-1}\bigcup_{l=1}^{\ell_2-1}\bar B_{\ell_1,\ell_2}^l(1,s)\cup\bigcup_{\ell_1=\ell^*}^{L-1}\bigcup_{\ell_2=\ell^*}^{L-1} \bar R_{\ell_1,\ell_2}(1,s)\notag\\
&\cup\bigcup_{\ell_1=\ell^*}^{L-1}\bigcup_{\ell_2=\ell^*}^{L-1}\bigcup_{l=1}^{\ell_2-1}\bar B_{\ell_1,\ell_2}^l(1,s)\cup\mathscr S_\text{neg}^v(1,s)\cup\mathscr S_\text{neg}^h(1,s).
\end{align}
To conclude we need to find the standard cascade from $\mathbf 1$ to $\mathcal X^s_\text{neg}$. Using Theorem \ref{teogatenegset} and the symmetry of the energy landscape with respect to the $q-1$ stable states, we complete the proof by taking the union of the standard cascades from $\mathbf 1$ to all possible $\mathbf s\in\mathcal X^s_\text{neg}$ given by \eqref{firstdescenteeta1}--\eqref{firstdescenteeta2}. Finally, \eqref{ristimetubenzbneg} follows by \cite[Lemma 3.13]{nardi2016hitting}.
$\qed$ 

\textit{Proof of Theorem \ref{teotubefixneg}} Let us assume $q>2$, otherwise the result is proven in \cite[Section 7.4]{olivieri2005large}. Starting from the metastable state $\mathbf 1$, the process hits $\mathcal X^s_\text{neg}$ in any stable state $\mathbf r$ with the same probability $\frac 1 {q-1}$. The set of typical paths $\Omega_{\mathbf 1,\mathbf s}^{\text{vtj}}$ may be partitioned in two subsets $\Omega_{\mathbf 1,\mathbf s}^{\text{vtj},1}:=\{\omega\in\Omega_{\mathbf 1,\mathbf s}^{\text{vtj}}:\ \omega\cap\mathcal X^s_\text{neg}\backslash\{\mathbf s\}=\varnothing\}$ and $\Omega_{\mathbf 1,\mathbf s}^{\text{vtj},2}:=\{\omega\in\Omega_{\mathbf 1,\mathbf s}^{\text{vtj}}:\ \omega\cap\mathcal X^s_\text{neg}\backslash\{\mathbf s\}\neq\varnothing\}$. Since the process follows a path belonging to $\Omega_{\mathbf 1,\mathbf s}^{\text{vtj},2}$ with probability $\frac{q-2}{q-1}>0$, these trajectories also belong to the tube of typical paths. Thus, the tube $\mathfrak T_{\mathbf s}(\mathbf 1)$ is comprised of those configurations that belong to all the typical paths that go from $\mathbf 1$ to $\mathcal X^s_\text{neg}$, i.e., those states belonging to $\mathfrak T_{\mathcal X^s_\text{neg}}(\mathbf 1)$, and of those configurations that belong to all typical paths from any $\mathbf r\in\mathcal X^s_\text{neg}\backslash\{\mathbf s\}$ to $\mathbf s$. Using Remark \ref{remarktube1s}, these last configurations belong to the tube $\mathfrak T_{\mathbf s}^{\text{zero}}(\mathbf r)$ given by \cite[Equation 4.25, Theorem 4.3]{bet2021critical}. Finally, we apply \cite[Lemma 3.13]{nardi2016hitting} to prove \eqref{ristimetubenzbnegfix}. $\qed$

\section{Sharp estimate on the mean transition time from the metastable state to the set of the stable states}\label{secprefactor}
In order to prove our main results on the computation of the prefactor and on the estimate of the expected value of the transition time from a metastable state to the stable set, we adopt the \textit{potential theoretic approach}. In order to apply this method, let us give some further definitions and some known results taken from \cite{bovier2002metastability,bovier2016metastability} and from \cite{baldassarri2021metastability}. 

We begin by introducing some further model-independent definitions and results. Consider any energy landscape $(\mathcal X,H,Q)$ and let $h:\mathcal X\to\mathbb R$. We define \textit{Dirichlet form} as
\begin{align}\label{dirichletform}
\mathfrak E_\beta(h)&:=\frac 1 2 \sum_{\sigma,\eta\in\mathcal X}\mu_\beta(\sigma)P_\beta(\sigma,\eta)[h(\sigma)-h(\eta)]^2\notag\\
&= \frac 1 2 \sum_{\sigma,\eta\in\mathcal X} \frac{e^{-\beta H(\sigma)}}{Z}\frac{e^{-\beta[H(\eta)-H(\sigma)]_+}}{|\Lambda|}[h(\sigma)-h(\eta)]^2.
\end{align}

\noindent Given two non-empty disjoint sets $\mathcal A_1,\mathcal A_2\subset\mathcal X$, the \textit{capacity} of the pair $\mathcal A_1,\mathcal A_2$ is defined by
\begin{align}\label{defcapacity}
\text{CAP}(\mathcal A_1,\mathcal A_2):=\min_{\substack{h:\mathcal X\to[0,1]\\ h_{|_{\mathcal A_1}}=1, h_{|_{\mathcal A_2}}=0}}\mathfrak E_\beta(h).
\end{align}
Note that from \eqref{defcapacity} it follows immediately that the capacity is symmetric in $\mathcal A_1$ and $\mathcal A_2$. In particular, the right hand side of \eqref{defcapacity} has a unique minimizer $h^*_{\mathcal A_1,\mathcal A_2}$ known as \textit{equilibrium potential} of $\mathcal A_1,\mathcal A_2$ and given by
$h^*_{\mathcal A_1,\mathcal A_2}(\eta)=\mathbb P(\tau^\eta_{\mathcal A_1}<\tau^\eta_{\mathcal A_2}),$
for any $\eta\in\mathcal X$. Finally, using what we have just defined, consider the following.
\begin{definition}\label{ptametastableset}
A set $\mathcal A\subset\mathcal X$ is said to be \textit{p.t.a.-metastable} if 
\begin{align}\label{alignptametastableset}
\lim_{\beta\to\infty} \frac{\max_{\sigma\notin\mathcal A}\mu_\beta(\sigma)[\text{CAP}_\beta(\sigma,\mathcal A)]^{-1}}{\min_{\sigma\in\mathcal A}\mu_\beta(\sigma)[\text{CAP}_\beta(\sigma,\mathcal A\backslash\{\sigma\})]^{-1}}=0.
\end{align}
The prefix p.t.a. stands for potential theoretic approach and it is used for distinguishing the Definition \ref{ptametastableset} from that of the metastable set $\mathcal X^m$. We remark that the idea of defining a set as in Definition \ref{ptametastableset} was introduced in \cite{bovier2002metastability}, where the authors refer to it as \textit{set of metastable points}. We refer to \cite{bovier2002metastability} and to \cite[Chapter 8]{bovier2016metastability} for the study of the main properties of this set. 

Since the identification of a p.t.a.-metastable set is quite difficult if one starts from the Definition \ref{ptametastableset}, we recall \cite[Theorem 3.6]{cirillo2013relaxation} where the authors give a constructive method for defining any p.t.a.-metastable set. In particular, for any $\sigma,\eta\in\mathcal X$, the authors introduced the following equivalence relation 
\begin{align}\label{equivalencerelation}
\sigma\sim\eta\ \text{if and only if}\ \Phi(\sigma,\eta)-H(\sigma)<\Gamma^m\ \text{and}\ \Phi(\eta,\sigma)-H(\eta)<\Gamma^m.
\end{align}
Assumed $\mathcal X\backslash\mathcal X^s\neq\varnothing$, let $\mathcal X^m_{(1)},\dots,\mathcal X^m_{(k_m)}$ and $\mathcal X^s_{(1)},\dots,\mathcal X^s_{(k_s)}$ be the equivalence classes in which $\mathcal X^m$ and $\mathcal X^s$ are partitioned with respect to the relation $\sim$, respectively. 
\begin{theorem}\emph{\cite[Theorem 3.6]{cirillo2013relaxation}}\label{teocirillonardi}
Assume that $\mathcal X\backslash\mathcal X^s\neq\varnothing$ and $\mathcal X\backslash(\mathcal X^s\cup\mathcal X^m)\neq\varnothing$. Choose arbitrarily $\sigma_{s,i}\in\mathcal X^s_{(i)}$ for any $i=1,\dots,k_s$ and $\sigma_{m,j}\in\mathcal X^m_{(j)}$ for any $j=1,\dots,k_m$. The set $\{\sigma_{s,1},\dots,\sigma_{s,k_s},\sigma_{m,1},\dots,\sigma_{m,k_m}\}$ is a p.t.a.-metastable.
\end{theorem}
\begin{remark}\label{remarklibrometpair}
In \cite[Chapters 8 and 16]{bovier2016metastability} the authors state the main metastability theorems for those energy landscapes in which the stable set $\mathcal X^s=\{\bold s\}$ and the metastable set $\mathcal X^m=\{\bold m\}$ are singletons. In particular, in \cite[Lemma 16.13]{bovier2016metastability} the authors prove that the pair $\mathcal A = \{\bold m,\bold s\}$ is a p.t.a.-metastable set.  
\end{remark}

\end{definition}
\subsection{Mean crossover time and computation of prefactor: proof of main results}\label{proofprefactor}
In this subsection we prove Theorem \ref{teoprefneg} by using the model independent results given in \cite{baldassarri2021metastability} and \cite[Chapter 16]{bovier2016metastability}, by exploiting the discussion given in \cite[Subsection 3.1]{cirillo2013relaxation} and also by using some results given in \cite{bovier2002metastability},\cite{bashiri2019on}. Let us begin by giving the following list of definitions and notations.
\begin{itemize}
\item[-] With an abuse of notation we consider $\mathcal X$ as a graph whose vertices are the configurations. Given two configurations $\sigma,\eta\in\mathcal X$ there is an edge between the corresponding vertices if it is possible to move from $\sigma$ to $\eta$ (resp.~$\eta$ to $\sigma$) in one step of the dynamics.
\item[-] 
Let $\mathcal X^*_\text{neg}\subset\mathcal X$ be the subgraph obtained by removing all the vertices corresponding to configurations $\sigma\in\mathcal X$ such that  $H_\text{neg}(\sigma)>\Gamma_\text{neg}^m+H_\text{neg}(\bold 1)$ and also removing all edges incident to these configurations. 
\item[-] Let $\mathcal X^{**}_\text{neg}\subset\mathcal X^*_\text{neg}$ be the subgraph obtained by removing all the vertices corresponding to configurations $\sigma$ such that $H_\text{neg}(\sigma)=\Gamma_\text{neg}^m+H_\text{neg}(\bold 1)$ and also removing all edges incident to these configurations. 

\item[-] Let $\mathscr P^*_\text{PTA}(\bold 1,\mathcal X^s_\text{neg})$ be the \textit{protocritical set} and let $\mathscr C^*_\text{PTA}(\bold 1,\mathcal X^s_\text{neg})$ be the \textit{critical set}. More precisely, we exploit \cite[Definition 16.3]{bovier2016metastability} and define $(\mathscr C^*_\text{PTA}(\bold 1,\mathcal X^s_\text{neg}),\mathscr P^*_\text{PTA}(\bold 1,\mathcal X^s_\text{neg}))$ as the maximal subset of $\mathcal X\times\mathcal X$ such that:
\noindent(1) for any $\sigma\in\mathscr P^*_\text{PTA}(\bold 1,\mathcal X^s_\text{neg})$ there exists $\eta\in\mathscr C^*_\text{PTA}(\bold 1,\mathcal X^s_\text{neg})$ such that $\sigma\sim\eta$ and for any $\eta\in\mathscr C^*_\text{PTA}(\bold 1,\mathcal X^s_\text{neg})$ there exists $\sigma\in\mathscr P^*_\text{PTA}(\bold 1,\mathcal X^s_\text{neg})$ such that $\eta\sim\sigma$;

\noindent(2) for any $\sigma\in\mathscr P^*_\text{PTA}(\bold 1,\mathcal X^s_\text{neg})$, $\Phi_\text{neg}(\sigma,\bold 1)<\Phi_\text{neg}(\sigma,\mathcal X^s_\text{neg})$;

\noindent(3) for any $\eta\in\mathscr C^*_\text{PTA}(\bold 1,\mathcal X^s_\text{neg})$ there exists a path $\omega:\eta\to\mathcal X^s_\text{neg}$ such that $\max_{\zeta\in\omega}H_\text{neg}(\zeta)-H_\text{neg}(\bold 1)\le\Gamma_\text{neg}^m$ and $\omega\cap\{\zeta\in\mathcal X:\ \Phi_\text{neg}(\zeta,\bold 1)<\Phi_\text{neg}(\zeta,\mathcal X^s_\text{neg})\}=\varnothing$. 
\end{itemize}
Next, consider $\mathcal W_{\text{neg}}(\bold 1,\mathcal X^s_{\text{neg}})=\mathcal G^1_\text{neg}\cup\mathcal G^2_\text{neg}$ where $\mathcal G^1_\text{neg}$ and $\mathcal G^2_\text{neg}$ are defined as follows.

\noindent- $\mathcal G^1_\text{neg}:=\{\sigma\in\mathcal W_{\text{neg}}(\bold 1,\mathcal X^s_{\text{neg}}):$ the cluster of spins different from $1$ has the unit protuberance on a corner of one of the longest sides of the quasi-square $\ell^*\times(\ell^*-1)\}.$

\noindent- $\mathcal G^2_\text{neg}:=\{\sigma\in\mathcal W_{\text{neg}}(\bold 1,\mathcal X^s_{\text{neg}})$: the cluster of spins different from $1$ has the unit protuberance on one of the $\ell^*-2$ vertices different from the corners of one of the longest sides of the quasi-square $\ell^*\times(\ell^*-1)\}.$
Following the same strategy given in \cite{baldassarri2021metastability}, let us consider the set
\begin{align}
\mathcal X^{**}_\text{neg}\backslash(\mathcal C^\mathbf 1_{\mathcal X^s_\text{neg}}(\Gamma^m_\text{neg})\cup\mathcal C^{\mathcal X^s_\text{neg}}_{\mathbf 1}(\Gamma_\text{neg}(\mathcal X^s_\text{neg},\mathbf 1)))=\bigcup_{i=1}^I\mathcal X(i),
\end{align}
where each $\mathcal X(i)$ is a set of communicating states with energy strictly lower than $\Phi_\text{neg}(\mathbf 1,\mathcal X^s_\text{neg})$ and with communication energy $\Phi_\text{neg}(\mathbf 1,\mathcal X^s_\text{neg})$ with respect to both $\mathbf 1$ and $\mathcal X^s_\text{neg}$. Among these sets we find also the wells $\mathcal Z^\mathbf 1_j$ (resp. $\mathcal Z^{\mathcal X^s_\text{neg}}_j$) that are connected by one step of the dynamics with the unessential saddles that in \cite[Definitions 3.2 and 3.4]{baldassarri2021metastability} are said to be ``of the first type'' (resp. ``of the second type'') and that are denoted by $\sigma_j$ (resp. $\zeta_j$). In view of the above discussion, let us define the following subsets of $\mathcal X^*_\text{neg}$.
\begin{itemize}
\item[-] $A_\text{neg}:=\mathcal C^\bold 1_{\mathcal X^s_\text{neg}}(\Gamma^m_\text{neg})\cup\bigcup_{j=1}^{J_{meta}}(\{\sigma_j\}\cup\mathcal Z^{\mathbf 1}_j)$.
\item[-] $B_\text{neg}:=\mathcal C^{\mathcal X^s_\text{neg}}_{\mathbf 1}(\Gamma_\text{neg}(\mathcal X^s_\text{neg},\mathbf 1))\cup\bigcup_{j=1}^{J_{stab}}(\{\zeta_j\}\cup\mathcal Z^{\mathcal X^s_\text{neg}}_j)$.
\end{itemize}

\noindent Before of the proof of Theorem \ref{teoprefneg}, it is useful to state the following results.

\begin{lemma}\label{remarkcardinalityG}
The cardinality of $\mathcal G^1_\emph{neg}$ and $\mathcal G^2_\emph{neg}$ are 
$|\mathcal G^1_\emph{neg}|=8|\Lambda|(q-1)$ and $|\mathcal G^2_\emph{neg}|=4|\Lambda|(\ell^*-2)(q-1)$, respectively.
\end{lemma}
\textit{Proof.} In $\mathcal G^1_\text{neg}$ the protuberance lies at one of the two extreme ends of one of the side of length $\ell^*$, hence there are four possible positions. On the other hand, in $\mathcal G^2_\text{neg}$ there are $2(\ell^*-2)$ sites in which can place the unit protuberance. In both cases, the quantity $2|\Lambda|$ counts the number of locations and rotations of the cluster with spins different from $1$. Indeed, the quasi-square with the unit protuberance may be located anywhere in $\Lambda$ in two possible orientations. Furthermore, the factor $(q-1)$ counts the number of possible spins that may characterize this homogenous cluster.
$\qed$

\begin{lemma}\label{lemmaptametaneg}
If the external magnetic field is negative, then the set $\{\bold 1,\mathcal X^s_\emph{neg}\}$ is p.t.a-metastable.
\end{lemma}
\textit{Proof.} Consider the equivalence relation $\sim$ given in \eqref{equivalencerelation}. From Theorem \ref{teometastableneg}, we get that in the energy landscape $(\mathcal X,H_\text{neg},Q)$ the metastable set is a singleton. Hence, there exists only one equivalence class with respect to $\sim$ given by $\mathcal X^m_\text{neg}$ itself. On the other hand, $\mathcal X^s_\text{neg}=\{\bold 2,\dots,\bold q\}$ and from Equation \eqref{comparisongammaneg} of Theorem \ref{theoremcomparisonneg} we get that $\mathcal X^s_{(1)}:=\{\bold 2\},\dots,\mathcal X^s_{(q-1)}:=\{\bold q\}$ are the equivalence classes with respect to the relation $\sim$ that partition $\mathcal X^s_\text{neg}$. Thus, by Theorem \ref{teocirillonardi} we conclude that the set $\{\bold 1,\bold 2,\dots,\bold q\}=\{\bold 1,\mathcal X^s_\text{neg}\}$ is p.t.a.-metastable. $\qed$

\begin{proposition}\label{identificationCstar}
If the external magnetic field is negative, then 
\begin{align}\label{alignCstarW}
\mathscr C^*_\emph{PTA}(\bold 1,\mathcal X^s_\emph{neg})=\mathcal W_{\emph{neg}}(\bold 1,\mathcal X^s_{\emph{neg}}).
\end{align}
\end{proposition}
\textit{Proof.} Following the same strategy of the proof of \cite[Theorem 17.3]{bovier2016metastability}, \eqref{alignCstarW} follows by the definition of $\mathscr C^*_\text{PTA}(\bold 1,\mathcal X^s_\text{neg})$, by Lemmas \ref{lemmaoneprefactor}--\ref{lemmathreeprefactor} and by Proposition \ref{propgateneg}. 
$\qed$
%

\begin{lemma}\label{lemmastepsellaneg}
Let $\eta\in\mathcal W_\emph{neg}(\mathbf 1,\mathcal X^s_\emph{neg})$ and let $\bar\eta\in\mathcal X$ such that $\bar\eta:=\eta^{v,t}$ for some $v\in V$ and $t\in S$, $t\neq\eta(v)$. If the external magnetic field is negative, then either $H_{\text{\emph{neg}}}(\eta)<H_{\text{\emph{neg}}}(\bar\eta)$ or $H_{\text{\emph{neg}}}(\eta)>H_{\text{\emph{neg}}}(\bar\eta)$.
\end{lemma}
\textit{Proof.} Since $\eta\in\mathcal W_\text{neg}(\mathbf 1,\mathcal X^s_\text{neg})=\bigcup_{t=2}^q\bar B_{\ell^*-1,\ell^*}^1(1,t)$, there exists $s\neq 1$ such that $\eta\in\bar B_{\ell^*-1,\ell^*}^1(1,s)$. This implies that $\eta$ is characterized by all spins $1$ except those, which are $s$, in a quasi-square $(\ell^*-1)\times\ell^*$ with a unit protuberance on one of the longest sides. In particular, for any $u\in V$, either $\eta(u)=1$ or $\eta(u)=s$. If $\eta(u)=1$, then for any $t\in S\backslash\{1\}$, depending on the distance between the vertex $u$ and the $s$-cluster, we have 
\begin{align}
&H_{\text{neg}}(\bar\eta)-H_{\text{neg}}(\eta)=
\begin{cases}
4-h\mathbbm 1_{\{t=s\}},\ &\text{if $n_1(u)=4$ };\\
3-\mathbbm 1_{\{t=s\}}-h\mathbbm 1_{\{t=s\}},\ &\text{if $n_1(u)=3$, $n_s(u)=1$};\\
2-2\mathbbm 1_{\{t=s\}}-h\mathbbm 1_{\{t=s\}},\ &\text{if $n_1(u)=2$, $n_s(u)=2$}.
\end{cases}
\end{align}
Otherwise, if $\eta(u)=1$, for any $t\in S\backslash\{1\}$, depending on the distance between the vertex $u$ and the boundary of the $s$-cluster,  we get
\begin{align}
&H_{\text{neg}}(\bar\eta)-H_{\text{neg}}(\eta)=
\begin{cases}
4+h,\ &\text{if $n_s(u)=4$};\\
3-\mathbbm 1_{\{t=1\}}+h,\ &\text{if $n_1(u)=1$, $n_s(u)=3$};\\
2-2\mathbbm 1_{\{t=1\}}+h,\ &\text{if $n_1(u)=2$, $n_s(u)=2$};\\
1-3\mathbbm 1_{\{t=1\}}+h,\ &\text{if $n_1(u)=3$, $n_s(u)=1$}.
\end{cases}
\end{align}
We conclude that $H_{\text{neg}}(\eta)\neq H_{\text{neg}}(\bar\eta)$. 
$\qed$

In \cite[Definitions  3.2 and 3.4]{baldassarri2021metastability} the authors define two subsets of unessential saddles for the metastable transition and they call them respectively unessential saddles of the first type'' and of the second type and in \cite[Equations (3.16)--(3.17)]{baldassarri2021metastability} they define the sets $K$ and $\tilde K$. Using these definitions and Lemma \ref{lemmastepsellaneg}, we are now able to prove the following. 
\begin{lemma}\label{lemmasimone}
If the external magnetic field is negative, then the following properties are verified.
\begin{itemize}
\item[\emph{(a)}] $K=\varnothing$, $\tilde K=\varnothing$.
\item[\emph{(b)}] Any $\sigma\in\mathcal W'_\emph{neg}(\mathbf 1,\mathcal X^s_\emph{neg})$ is such that $\sigma\in\bigcup_{j=1}^{J_{meta}}(\{\sigma_j\}\cup\mathcal Z^{\mathbf 1}_j)$, namely there exist at least a unessential saddle $\sigma_i$ ``of the first type'' and its well $\mathcal Z^\mathbf 1_i$ is not empty. 
\item[\emph{(c)}] The set $\bigcup_{j=1}^{J_{stab}}(\{\zeta_j\}\cup\mathcal Z^{\mathcal X^s_\emph{neg}}_j)$ is not empty, namely there exists at least a unessential saddle $\zeta_i$ ``of the second type''. 
\end{itemize}
\end{lemma}
\textit{Proof}. By Lemma \ref{lemmastepsellaneg} we have that any $\eta\in\mathcal W_\text{neg}(\mathbf 1,\mathcal X^s_\text{neg})$ that communicates with configurations in the cycles $\mathcal C^\bold 1_{\mathcal X^s_\text{neg}}(\Gamma^m_\text{neg})\cup\mathcal C^{\mathcal X^s_\text{neg}}_{\mathbf 1}(\Gamma_\text{neg}(\mathcal X^s_\text{neg},\mathbf 1))$, in $\mathcal X\backslash\mathcal X^*_\text{neg}$, and it does not communicate by a single step of the dynamics with another saddle.
This implies that for any $\bar\eta\in\mathcal S_\text{neg}(\mathbf 1,\mathcal X^s_\text{neg})\backslash\mathcal W_\text{neg}(\mathbf 1,\mathcal X^s_\text{neg})$, visited by the process before visiting the gate $\mathcal W_\text{neg}(\mathbf 1,\mathcal X^s_\text{neg})$, it does not exist a path $\omega_1:\eta\to\bar\eta$ such that $\omega_1\cap\mathcal C^\bold 1_{\mathcal X^s_\text{neg}}(\Gamma^m_\text{neg})=\varnothing$, $\omega_1\cap\mathcal W_\text{neg}(\mathbf 1,\mathcal X^s_\text{neg})=\{\eta\}$, and $\max_{\sigma\in\omega_1}H_\text{neg}(\sigma)\le\Phi_\text{neg}(\mathbf 1,\mathcal X^s_\text{neg})$. This concludes that $K=\varnothing$. Furthermore, for any $\bar\eta\in\mathcal S_\text{neg}(\mathbf 1,\mathcal X^s_\text{neg})\backslash\mathcal W_\text{neg}(\mathbf 1,\mathcal X^s_\text{neg})$, visited by the process after visiting the gate $\mathcal W_\text{neg}(\mathbf 1,\mathcal X^s_\text{neg})$, there does not exist $\omega_1:\eta\to\bar\eta$ such that $\omega_1\cap\mathcal C^{\mathcal X^s_\text{neg}}_{\mathbf 1}(\Gamma_\text{neg}(\mathcal X^s_\text{neg},\mathbf 1))=\varnothing$, $\omega_1\cap\mathcal W_\text{neg}(\mathbf 1,\mathcal X^s_\text{neg})=\{\eta\}$, and $\max_{\sigma\in\omega_1}H_\text{neg}(\sigma)\le\Phi_\text{neg}(\mathbf 1,\mathcal X^s_\text{neg})$. This concludes that $\tilde K=\varnothing$ and the proof of item (a).

\noindent Let us now prove item (b). Using Theorem \ref{teogatenegset}, we get that any saddle in which the protuberance is on one of the shortest sides: $\sigma_i\in\mathcal W'_\text{neg}(\mathbf 1,\mathcal X^s_\text{neg})$, is an unessential saddle. Thus, $\sigma_i$ satisfies \cite[Definition 3.2]{baldassarri2021metastability} and it belongs to $\bigcup_{j=1}^{J_{meta}}(\{\sigma_j\}\cup\mathcal Z^{\mathbf 1}_j)$. Moreover, if $\sigma_i\in\bar B^1_{\ell^*,\ell^*-1}(1,s)$, and without loss of generality the protuberance is on the shortest side that is north, then it communicates by one step of the dynamics with a configuration in $\bar B^2_{\ell^*,\ell^*-1}(1,s)$ with a bar of length two on the north side. This belongs to $\mathcal Z^\mathbf 1_i$ together with those configurations with a bar of length $l$ on the north side belonging to $\bar B^l_{\ell^*,\ell^*-1}(1,s)$ for any $l=3,\dots,\ell^*-2$ and its bottom is a configuration belonging to $\bar R_{\ell^*-1,\ell^*+1}(1,s)$ with the shortest sides that are north and south. The same arguments hold by replacing north with south, east, west. 

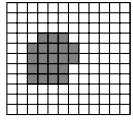
\begin{figure}
\centering
\begin{tikzpicture}[scale=0.45,transform shape]
\fill[gray] (-0.3,0) rectangle (0,1.2) (0,0) rectangle (0.9,1.5) (0.9,0.6) rectangle (1.2,1.2);
\draw[step=0.3cm,color=black] (-0.3,0) grid (0,1.2) (0,0) grid (0.9,1.5) (0.9,0.6) grid (1.2,1.2);
\draw[step=0.3cm,color=black] (-0.9,-0.9) grid (2.7,2.4);
\draw (-0.9,-0.9) rectangle (2.7,2.4);
\end{tikzpicture}
\caption{\label{figuraunessentialsecond} Example of a unessential saddle $\zeta$ ``of the second type'' defined in \cite{baldassarri2021metastability} when $\ell^*=5$. We color white the vertices with spin $1$ and gray the vertices with spin $s\neq 1$.}
\end{figure}

\noindent Let us now prove item (c) by illustrating an example of unessential saddle ''of the second type''. We choose this unessential saddles as the configuration $\zeta\in\partial\mathcal C^{\mathcal X^s_\text{neg}}_{\mathbf 1}(\Gamma_\text{neg}(\mathcal X^s_\text{neg},\mathbf 1))\cap(\mathcal S_\text{neg}(\mathbf 1,\mathcal X^s_\text{neg})\backslash\mathcal W_\text{neg}(\mathbf 1,\mathcal X^s_\text{neg}))$ in which all the vertices have spin equal to $1$ except those, which are all equal to $s$ for some $s\neq 1$, in a cluster that is a square $(\ell^*-1)\times(\ell^*-1)$ with a bar of length two on one of the four sides and a bar of length $\ell^*-2$ on one of the two consecutive sides, see Figure \ref{figuraunessentialsecond}. 
\FloatBarrier
Note that $\zeta\in\mathcal S_\text{neg}(\mathbf 1,\mathcal X^s_\text{neg})\backslash\mathcal W_\text{neg}(\mathbf 1,\mathcal X^s_\text{neg})$ since the perimeter of the $s$-cluster is $4\ell^*$ and since its area is equal to $\ell^*(\ell^*-1)+1$, and so by \eqref{rewritegap1neg} we get that
%
$H_\text{neg}(\zeta)=H_\text{neg}(\mathbf 1)+4\ell^*-h(\ell^*(\ell^*-1)+1))=\Phi_\text{neg}(\mathbf 1,\mathcal X^s_\text{neg}).$
%
Furthermore, $\zeta\in\partial\mathcal C^{\mathcal X^s_\text{neg}}_{\mathbf 1}(\Gamma_\text{neg}(\mathcal X^s_\text{neg},\mathbf 1))$. Indeed, by flipping to $s$ the spin $1$ adjacent to the bar of length $\ell^*-2$, the process intersects a configuration belonging to $\bar B^2_{\ell^*-1,\ell^*}(1,s)\subset\mathcal C^{\mathcal X^s_\text{neg}}_{\mathbf 1}(\Gamma_\text{neg}(\mathcal X^s_\text{neg},\mathbf 1))$.
 $\qed$

Now we are able to give the proof of Theorem \ref{teoprefneg}. Since our model is under Glauber dynamics, we exploit the proof of \cite[Theorem 17.4]{bovier2016metastability}.

\textit{Proof of Theorem \ref{teoprefneg}}. Let us begin to compute the prefactor \eqref{Knegative} by exploiting the variational formula for $\Theta_\text{neg}=1/{K_\text{neg}}$ given in \cite[Lemma 10.7]{baldassarri2021metastability}. This variational problem is simplified because of our Glauber dynamics. Indeed, from the definition of $A_\text{neg}$ and $B_\text{neg}$ and from Proposition \ref{identificationCstar}, we get that $\mathcal X^*_\text{neg}\backslash(A_\text{neg}\cup B_\text{neg})=\mathscr C^*_\text{PTA}(\bold 1,\mathcal X^s_\text{neg})$. It follows that there are no wells inside $\mathscr C^*_\text{PTA}(\bold 1,\mathcal X^s_\text{neg})$ and any critical configuration may not transform into each other via single spin-update. We proceed by computing a lower and un upper bound for $\Theta_\text{neg}$ as follows.

\noindent\textit{Upper bound.} In order to estimate un upper bound for the capacity we choose a test function $h:\mathcal X^*_\text{neg}\to\mathbb R$ defined as
\begin{align} h(\sigma):=
\begin{cases}
1,\ &\text{if}\ \sigma\in A_\text{neg},\\
0,\ &\text{if}\ \sigma\in B_\text{neg},\\
c_i,\ &\text{if}\ \sigma\in\mathcal G_\text{neg}^i, i=1,2,
\end{cases}
\end{align}
where $c_1,c_2$ are two constants, see \cite[Equation (10.17)]{baldassarri2021metastability}. Thus, we get
\begin{align}
\Theta_\text{neg}&\le(1+o(1))\min_{c_1,c_2\in [0,1]}\min_{\substack{h:\mathcal X^*_\text{neg}\to[0,1]\\ h_{|_{A_\text{neg}}}=1, h_{|_{B_\text{neg}}}=0\\ h_{|_{\mathcal G_\text{neg}^i}}=c_i, i=1,2}} \frac 1 2 \sum_{\sigma,\eta\in\mathcal X^*_\text{neg}} \mathbbm 1_{\{\sigma\sim\eta\}}[h(\sigma)-h(\eta)]^2\notag \\
&=(1+o(1))\min_{c_1,c_2\in [0,1]}[\sum_{\substack{\sigma\in A_\text{neg}\\ \eta\in\mathcal G_\text{neg}^i, i=1,2\\\sigma\sim\eta}} (1-h(\eta))^2+\sum_{\substack{\sigma\in B_\text{neg}\\ \eta\in\mathcal G_\text{neg}^i, i=1,2\\\sigma\sim\eta}} h(\eta)^2]\notag \\
&=(1+o(1))\min_{c_1,c_2\in [0,1]}[\sum_{\substack{\eta\in\mathcal G_\text{neg}^i, i=1,2\\\sigma\sim\eta}}N^-(\eta)(1-c_i)^2+\sum_{\substack{\eta\in\mathcal G_\text{neg}^i, i=1,2\\\sigma\sim\eta}}N^+(\eta)c_i^2]
\end{align}
where $N^-(\eta) :=|\{\xi\in\bigcup_{t=2}^q\bar R_{\ell^*-1,\ell^*}(1,t):\ \xi\sim\eta\}|$, and $N^+(\eta):= |\{\xi\in\bigcup_{t=2}^q\bar B_{\ell^*-1,\ell^*}^2(1,t):\ \xi\sim\eta\}|$. Let us note that
\begin{align}\label{Nmenonum}
N^-(\eta)&=1,\ \text{if}\ \eta\in\mathcal G^1_\text{neg}\cup\mathcal G^2_\text{neg},\ \text{and}\ 
N^+(\eta)&=\begin{cases}1,\ &\text{if}\ \eta\in\mathcal G^1_\text{neg},\\ 2,\ &\text{if}\ \eta\in\mathcal G^2_\text{neg}.\end{cases}
\end{align}
Thus, we have
\begin{align}
\Theta_\text{neg}&\le(1+o(1))\min_{c_1,c_2\in [0,1]}[\sum_{\eta\in\mathcal G^1_\text{neg}} (1-c_1)^2+c_1^2+\sum_{\eta\in\mathcal G^2_\text{neg}} (1-c_2)^2+2c_2^2]\notag \\
&=(1+o(1))\min_{c_1,c_2\in [0,1]}[|\mathcal G^1_\text{neg}|(2c_1^2-2c_1+1)+|\mathcal G^2_\text{neg}|(3c^2_2-2c_2+1)],\notag
\end{align}
where the equality follows by the fact that the sums are independent from $\eta\in\mathcal G_\text{neg}^i$, $i=1,2$. Furthermore, since the minimum value of the function $g_1(c_1):=2c_1^2-2c_1+1$ is $\frac 1 2$ and the minimum value of the function $g_2(c_2):=3c^2_2-2c_2+1$ is $\frac 2 3$, we have
\begin{align}
\Theta_\text{neg}\hspace{-2pt}=\hspace{-2pt}|\mathcal G^1_\text{neg}|\frac 1 2\hspace{-2pt}+\hspace{-2pt}|\mathcal G^2_\text{neg}|\frac 2 3\hspace{-2pt}=\hspace{-2pt}\frac 1 2 8|\Lambda|(q\hspace{-2pt}-\hspace{-2pt}1)\hspace{-2pt}+\hspace{-2pt}\frac 2 3 4|\Lambda|(\ell^*\hspace{-2pt}-\hspace{-2pt}2)(q\hspace{-2pt}-\hspace{-2pt}1)\hspace{-2pt}=\hspace{-2pt}\frac 4 3 |\Lambda|(2\ell^*\hspace{-2pt}-\hspace{-2pt}1)(q\hspace{-2pt}-\hspace{-2pt}1),\notag
\end{align}
where the second equality follows by Lemma \ref{remarkcardinalityG}. 

\noindent\textit{Lower bound.} Since the variational formula for $\Theta_\text{neg}=1/{K_\text{neg}}$ given in \cite[Lemma 10.7]{baldassarri2021metastability} is defined by a sum with only non-negative summands, we obtain a lower bound for $\Theta_{\text{neg}}$ as follows 
\begin{align}
\Theta_\text{neg}\ge\min_{c_1,c_2\in [0,1]}\min_{\substack{h:\mathcal X^*_\text{neg}\to[0,1]\\ h_{|_{A_\text{neg}}}=1, h_{|_{B_\text{neg}}}=0\\ h_{|_{\mathcal G_\text{neg}^i}=c_i,} i=1,2}} \frac 1 2 \sum_{\sigma,\eta\in(\mathscr C^*_\text{PTA}(\mathbf 1,\mathcal X^s_\text{neg}))^+} \mathbbm 1_{\{\sigma\sim\eta\}}[h(\sigma)-h(\eta)]^2\notag
\end{align}
where $(\mathscr C^*_\text{PTA}(\mathbf 1,\mathcal X^s_\text{neg}))^+:=\mathscr C^*_\text{PTA}(\mathbf 1,\mathcal X^s_\text{neg})\cup\partial\mathscr C^*_\text{PTA}(\mathbf 1,\mathcal X^s_\text{neg})$.

\noindent Note that $\partial\mathscr C^*_\text{PTA}(\mathbf 1,\mathcal X^s_\text{neg})\cap\mathcal X^*_\text{neg}=\bigcup_{s=2}^q (\bar R_{\ell^*,\ell^*-1}(1,s)\cup\bar B^2_{\ell^*,\ell^*-1}(1,s))$, with $\bigcup_{s=2}^{q} \bar R_{\ell^*,\ell^*-1}(1,s)\subset\mathcal C^\mathbf 1_{\mathcal X^s_\text{neg}}(\Gamma^m_\text{neg})$ and $\bigcup_{s=2}^{q}\bar B^2_{\ell^*,\ell^*-1}(1,s)\subset\mathcal C^{\mathcal X^s_\text{neg}}_{\mathbf 1}(\Gamma_\text{neg}(\mathcal X^s_\text{neg},\mathbf 1))$. Thus, we have
\begin{align}
\Theta_\text{neg}&\ge\min_{h:\mathcal X^*_\text{neg}\to[0,1]} \sum_{\eta\in\mathscr C^*_\text{PTA}(\mathbf 1,\mathcal X^s_\text{neg})} \biggl(\sum_{\substack{\sigma\in\bigcup_{s=2}^{q}\bar R_{\ell^*,\ell^*-1}(1,s),\\ \sigma\sim\eta}}\hspace{-22pt}[1-h(\eta)]^2+\hspace{-12pt}\sum_{\substack{\sigma\in\bigcup_{s=2}^{q} \bar B^2_{\ell^*,\ell^*-1}(1,s),\\\sigma\sim\eta}}\hspace{-22pt}h(\eta)^2\biggr)\notag\\
&=\sum_{\sigma,\eta\in\mathscr C^*_\text{PTA}(\mathbf 1,\mathcal X^s_\text{neg})} \min_{h\in[0,1]}\biggl(N^-(\eta)[1-h]^2+N^+(\eta)h^2\biggr).
\end{align}
Since the minimizer of the function $f(h):=N^-(\eta)[1-h]^2+N^+(\eta)h^2$ is $h_\text{min}=\frac{N^-(\eta)}{N^-(\eta)+N^+(\eta)}$, we obtain
\begin{align}
\Theta_\text{neg}&\ge\sum_{\sigma,\eta\in\mathscr C^*_\text{PTA}(\mathbf 1,\mathcal X^s_\text{neg})} \frac{N^-(\eta)N^+(\eta)}{N^-(\eta)+N^+(\eta)}
%
=\frac 4 3 |\Lambda|(2\ell^*-1)(q-1),
\end{align}
where the first equality follows by \eqref{Nmenonum}. Finally, \eqref{meantimeKneg} is proven following the strategy given in \cite[Subsection 16.3.2]{bovier2016metastability} by taking into account the metastable set $\{\bold 1,\mathcal X^s_\text{neg}\}$ by replacing the role of Lemma 16.17 with \cite[Lemma 10.7]{baldassarri2021metastability}, see Remark \ref{remarklibrometpair} and Lemma \ref{lemmaptametaneg}.
$\qed$ 
\appendix

\section{Appendix}
\subsection{Additional material for Subsection \ref{submetastableneg}}\label{appendixpath}
\subsubsection{Explicit calculation of the inequality \eqref{disequalityHneg}}
We have
\begin{align}
&H_{\text{neg}}(\hat\omega_{k^*})-H_{\text{neg}}(\bold 1)=4\ell^*-h(\ell^*(\ell^*-1)+1),\notag\\
&H_{\text{neg}}(\hat\omega_{(K-1)^2+1})-H_{\text{neg}}(\bold 1)=4K-4-h(K-1)^2-h.\notag
\end{align}
Note that
\begin{align}\label{alignusingHcomp}
H_{\text{neg}}(\hat\omega_{k^*})-H_{\text{neg}}(\hat\omega_{(K-1)^2+1})
&=4\ell^*-h(\ell^*)^2+h\ell^*-4K+4+hK^2-2hK+h.
\end{align}
Using the constraints of Assumption \ref{remarkconditionneg} it follows that, we may write $\ell^*=\frac{2}{h}+1-\delta$ where $0<\delta<1$ denotes the fractional part of $2/h$. Hence, using \eqref{alignusingHcomp}, we get
\begin{align}\label{comparisonHpath}
&H_{\text{neg}}(\hat\omega_{k^*})\le H_{\text{neg}}(\hat\omega_{(K-1)^2+1})\\
\iff&\ 4\ell^*-h(\ell^*)^2+h\ell^*-4K+4+hK^2-2hK+h \le 0\notag \\
\iff&-\frac 4 h(\frac{2}{h}+1-\delta)+(\frac{2}{h}+1-\delta)^2-(\frac{2}{h}+1-\delta)+\frac 4 h K-\frac 4 h-K^2+2K-1\ge0\notag\\
\iff& -\frac 8{h^2}-\frac 4 h+\frac 4 h \delta+\frac4 {h^2}+1+\delta^2+\frac 4 h-\frac 4 h \delta-2\delta-\frac 2 h -1 +\delta+\frac 4 h K-\frac 4 h-1\ge K^2-2K\notag\\
\iff&-\frac 4 {h^2}-\frac 6 h+\frac 4 h K+\delta^2-\delta-1\ge K^2-2K.\notag
\end{align}
Since $K\ge 3\ell^*=3(\frac{2}{h}+1-\delta)$ and since $0<\delta<1$, it follows that \[K^2-2K\ge K(3\ell^*)-2K=3K(\frac{2}{h}+1-\delta)-2K=\frac 6 h K+K-3K\delta>\frac 6 h K-2K.\] Moreover, since $0<\delta<1$ implies that $\delta^2-\delta<0$, we have that
\begin{align}
-\frac 4 {h^2}-\frac 6 h+\frac 4 h K+\delta^2-\delta-1<-\frac 4 {h^2}-\frac 6 h+\frac 4 h K.
\end{align}
Hence, approximately we get that \eqref{comparisonHpath} is verified if and only if
\begin{align}\label{absurdappendifinal}
-\frac 4 {h^2}-\frac 6 h+\frac 4 h K> \frac 6 h K-2K\iff -\frac4 {h^2}-\frac 6 h-\frac 2 h K+2K> 0,\notag
\end{align}
that is an absurd because of the l.h.s. is strictly negative. Indeed, Assumption \ref{remarkconditionneg}(ii), i.e., $0<h<1$, implies that $-\frac 2 h K+2K=2K(1-\frac1 h)<0$. Thus, \eqref{comparisonHpath} is not verified and 
\begin{align}
H_{\text{neg}}(\hat\omega_{k^*})> H_{\text{neg}}(\hat\omega_{(K-1)^2+1}).
\end{align}
\subsection{Additional material for Subsection \ref{commheightstable}}\label{appendixstablephi}
\subsubsection{Proof of Proposition \ref{upperboundcaso2ss}}
\textit{Proof.} Our aim is to prove \eqref{upperboundalign} by constructing a path $\omega:\bold r\to\bold s$ such that 
\begin{align}
\Phi_\omega^{\text{neg}}-H_{\text{neg}}(\bold r)= 2\min\{K,L\}+2=2K+2,
\end{align}
where the last equality follows by our assumption $L\ge K$. Let $\sigma^*\in\mathcal X$ be the configuration defined as 
\begin{align}
\sigma^*(v):=
\begin{cases}
s,\ &\text{if}\ v\in c_0,\\
r,\ &\text{otherwise}.
\end{cases}
\end{align}
We define the path $\omega$ as the concatenation of the two paths $\omega^{(1)}:\bold r\to\sigma^*$ and $\omega^{(2)}:\sigma^*\to\bold s$ such that
$\Phi_{\omega^{(1)}}^{\text{neg}}=H_{\text{neg}}(\bold r)+2K\ \ \text{and}\ \ \Phi_{\omega^{(2)}}^{\text{neg}}=H_{\text{neg}}(\bold r)+2K+2.$
We define $\omega^{(1)}:=(\omega^{(1)}_0,\dots,\omega^{(1)}_K)$ where $\omega^{(1)}_0=\bold r$ and where for any $i=1,\dots,K$ the state $\omega^{(1)}_i$ is obtained by flipping the spin on the vertex $(i-1,0)$ from $r$ to $s$. The energy difference at each step of the path is
\begin{align}\label{energydifferencerefpathNZ}
H_{\text{neg}}(\omega^{(1)}_i)-H_{\text{neg}}(\omega^{(1)}_{i-1})=
\begin{cases}
4, &\text{if}\ i=1,\\
2, &\text{if}\ i=2,\dots,K-1,\\
0, &\text{if}\ i=K.
\end{cases}
\end{align}
Hence, $\text{arg max}_{\omega^{(1)}}=\{\omega^{(1)}_{K-1}, \omega^{(1)}_K=\sigma^*\}$. Indeed, in view of the periodic boundary conditions and of the \eqref{energydifferencerefpathNZ}, we have
\begin{align}\label{diffstablealign}
H_{\text{neg}}(\omega^{(1)}_{K-1})-H_{\text{neg}}(\bold r)=2K=H_{\text{neg}}(\omega^{(1)}_{K})-H_{\text{neg}}(\bold r).
\end{align}
Therefore, $\Phi_{\omega^{(1)}}^{\text{neg}}=H_{\text{neg}}(\bold r)+2K$. 
Let us now define the path $\omega^{(2)}$. We note that $\sigma^*$ has an $s$-bridge on column $c_0$ and so we apply to it the expansion algorithm introduced in Proposition \ref{expansionalg}. The algorithm gives a path $\omega^{(2)}:\sigma^*\to\bold s$ such that $\Phi_{\omega^{(2)}}^{\text{neg}}=H_{\text{neg}}(\sigma^*)+2=H_{\text{neg}}(\bold r)+2K+2$, where the last equality follows by \eqref{diffstablealign}.
$\qed$

\bibliographystyle{abbrv}
\bibliography{mybib}
\end{document}